\documentclass{article}

\usepackage[margin=1.5in]{geometry}

\usepackage[T1]{fontenc}
\usepackage[latin9]{inputenc}
\usepackage{array}
\usepackage{longtable}
\usepackage{float}
\usepackage{amsmath}
\usepackage{amsthm}
\usepackage{mathtools}
\usepackage{amssymb}
\usepackage{graphicx}
\usepackage{esint}
\usepackage{dsfont}
\usepackage{hyperref}
\usepackage{psfrag,epsf}
\usepackage{enumerate}
\usepackage[numbers]{natbib}

\newtheorem{theorem}{Theorem}[section]

\newtheorem*{remark}{Remark}

\allowdisplaybreaks

\newcommand{\RN}[1]{%
  \textup{\uppercase\expandafter{\romannumeral#1}}%
}

\title{\bf Extreme Value Analysis Without the Largest Values: What Can Be Done?}
  \author{Jingjing Zou \thanks{
    The authors would like to thank Zhi-Li Zhang for providing the Google+ data. This research is funded by ARO MURI grant W911NF-12-1-0385.} \\
    Department of Statistics, Columbia University
	\and
    Richard A. Davis \\
    Department of Statistics, Columbia University 
	\and
	Gennady Samorodnitsky \\
	School of Operations Research and Information Engineering\\
	Cornell University
}

\date{}

\begin{document}
\maketitle

\begin{abstract}

In this paper we are concerned with the analysis of heavy-tailed data when a portion of the extreme values is unavailable. 
This research was motivated by an analysis of the degree distributions in a large social network. The degree distributions of such networks tend to have power law behavior in the tails. 
We focus on the Hill estimator, which plays a starring role in heavy-tailed modeling. 
The Hill estimator for this data exhibited a smooth and increasing ``sample path'' as a function of the number of upper order statistics used in constructing the estimator.  
This behavior became more apparent as we artificially removed more of the upper order statistics.
Building on this observation we introduce a new version of the Hill estimator. It is a function of the number of the upper order statistics used in the estimation, but also depends on the number of unavailable extreme values.
We establish functional convergence of the normalized Hill estimator to a Gaussian process. 
An estimation procedure is developed based on the limit theory to estimate the number of missing extremes and extreme value parameters including the tail index and the bias of Hill's estimator.
We illustrate how this approach works in both simulations and real data examples.
  
\end{abstract}

\noindent%
\vfill

\newpage

\section{Introduction} \label{sec: intro}

In studying data exhibiting heavy-tailed behavior, a widely used model is the family of distributions that are regular varying. A distribution $ F $ is regularly varying of index $\alpha$ if
\begin{equation} \label{reg varying}
\frac{\bar{F}(tx)}{\bar{F}(t)} \rightarrow x^{-\alpha}  
\end{equation} 
as $ t \rightarrow \infty $ for all $ x > 0 $, where $ \alpha > 0 $ and $ \bar{F}(t) = 1- F(t) $ is the survival function.
The parameter $ \alpha $ is called the tail index or the extreme value index, and it controls the heaviness of the tail of the distribution. 
This is perhaps the most important parameter in extreme value theory and a great deal of research has been devoted to its estimation. The most used and studied estimate of $ \alpha $ is based on the Hill estimator for its reciprocal $ \gamma = 1/\alpha $ (see \citet{Hill:1975iu}, \citet{Drees:2000ic} and Section 2.1 of \citet{deHaan:2006bz} for further discussion on this estimator). The Hill estimator is defined by
\begin{equation*}
H_n(k) =  \frac{1}{k} \sum_{i = 1}^{k} \log X_{(i)} - \log X_{(k+1)},
\end{equation*}
for $k = 1, 2, \dots, n-1$, 
where $ X_{(1)} \ge X_{(2)} \ge \dots \ge X_{(n)}$ are the order statistics of an independent and identically distributed (iid) sample $ X_1, X_2, \dots, X_n \sim F(x) $. 
The Hill estimator is consistent in estimating $\gamma$: $H_n(k) \stackrel{P}{\rightarrow} \gamma$ as $n \rightarrow \infty$, $k = k(n) \rightarrow \infty$ and $k/n \rightarrow 0$ (see, for example, Section 3.2 of \citet{deHaan:2006bz}).

\begin{figure}[H]
  \centering
  \begin{minipage}[b]{0.45\textwidth} 
  \centering
    \includegraphics[width=\textwidth]{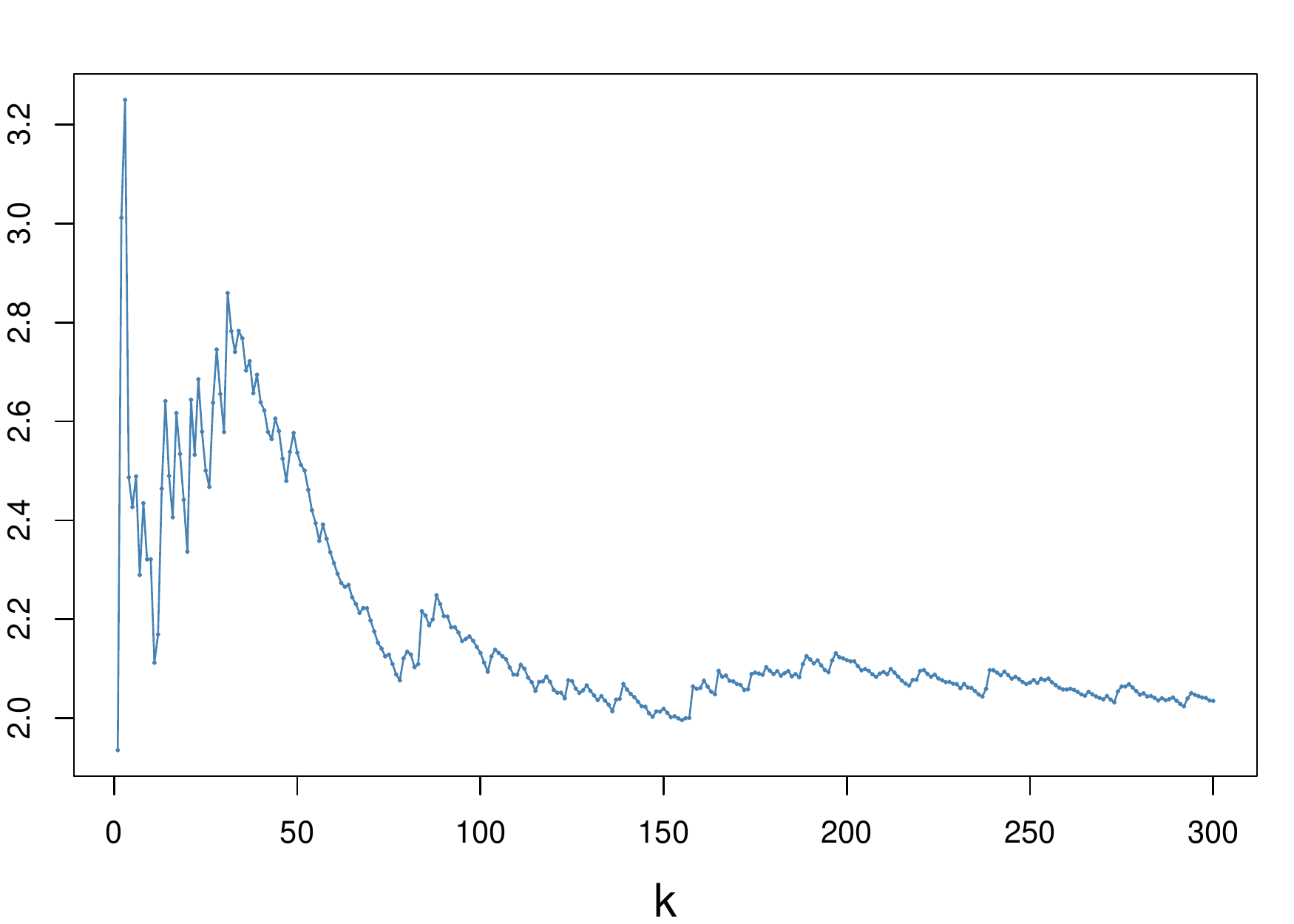}
  \end{minipage}
  \hfill
  \begin{minipage}[b]{0.45\textwidth}  
  \centering
    \includegraphics[width=\textwidth]{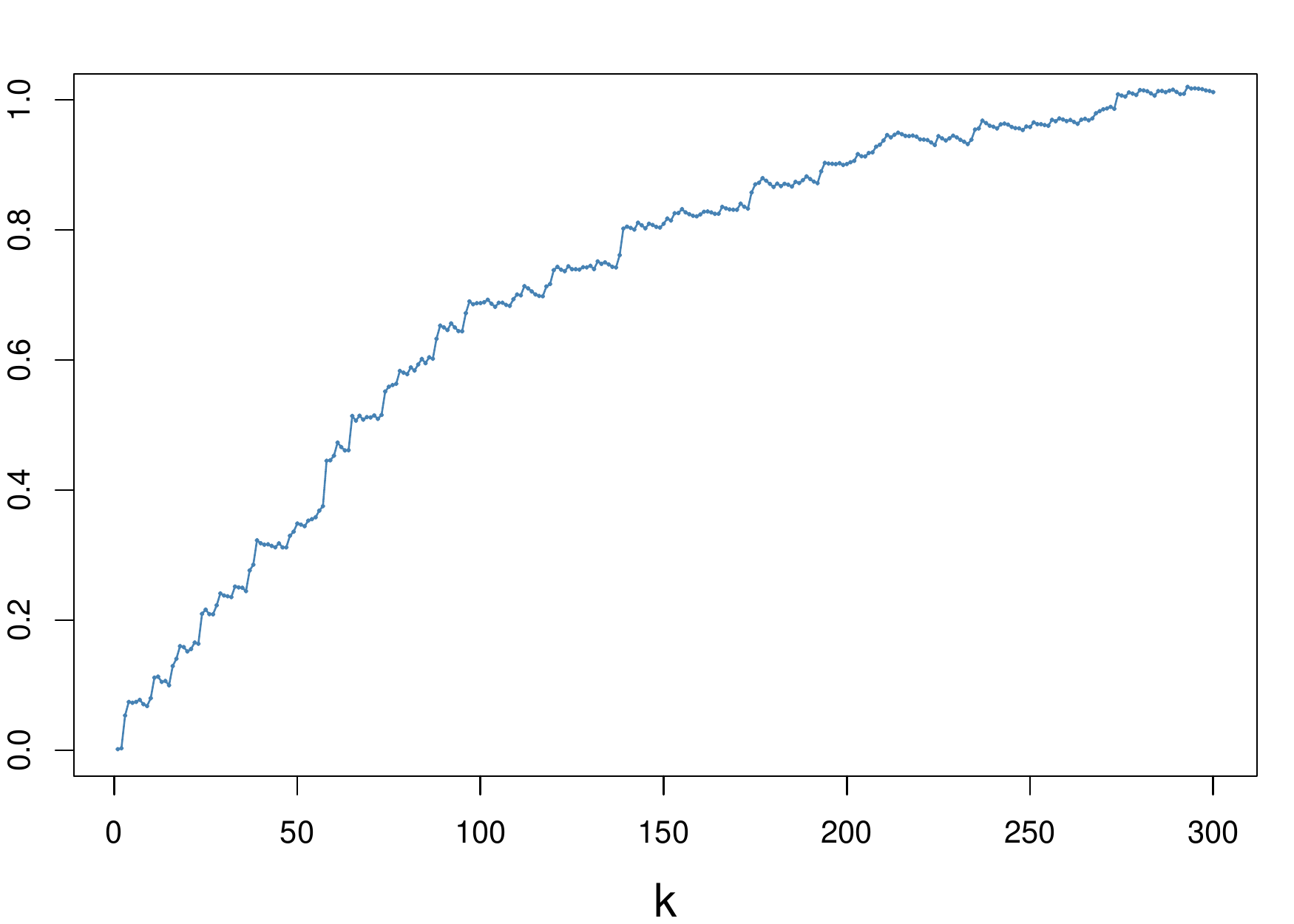}
  \end{minipage} 
\caption{Hill plot of iid Pareto ($\gamma = 2$) variables ($ n = 1000 $). 
$ x $-axis: number $ k $ of upper order statistics used in the calculation. $ y $-axis: $ H_n(k) $.
Left: without removal. Right: top $ 100 $ removed} \label{hill_pareto}
\end{figure}

As an illustration, the left panel of Figure \ref{hill_pareto} shows the Hill plot of $ 1000 $ iid observations from a Pareto distribution with $ \gamma = 2 $ ($ F(x) = 1- x^{-0.5} $ for $ x \ge 1 $ and $ 0 $ otherwise). 
In general, one chooses $ k $ for which the Hill plot remains relatively horizontal and uses the corresponding value of $ H_n(k) $ as the estimate for $ \gamma $. 

If the largest several observations in the data are removed, the Hill curve behaves very differently. For example, when the $ 100 $ largest observations of the previous Pareto sample have been removed, the Hill plot renders a much smoother curve that is generally increasing.
As seen in the right panel of Figure \ref{hill_pareto}, the Hill plot has no region in which it is horizontal.
Hence the choice of $ k $ in the Hill estimator is problematic in the presence of missing extremes, especially if the number of missing is unknown. The principle objective of this paper is to estimate the number of missing extremes simultaneously with other relevant parameters, including the tail index $ \alpha $.

\begin{figure}[H] 
\begin{center}
\includegraphics[width= 0.98 \textwidth]{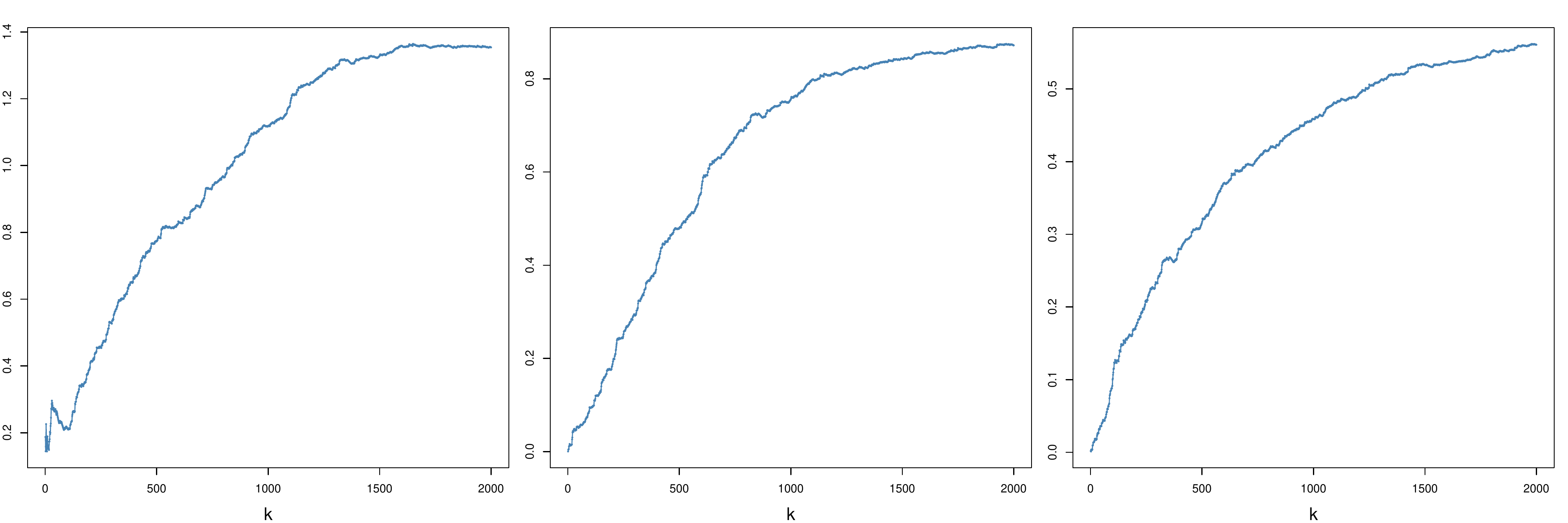}
\end{center}
\caption{Hill plots of in-degrees of the Google+ network. Left: without removal. Middle: $ 500 $ largest values removed. Right: $ 1000 $ largest values removed} \label{figure: GPlus 3plots}
\end{figure}

As a real-world example, a similar phenomenon is observed when we study the tail behavior of the in- and out-degrees in a large social network. We looked at data from a snapshot of Google+, the social network owned and operated by Google, taken on October 19, 2012. The data contain 76,438,791 nodes (registered users) and 1,442,504,499 edges (directed connections). 
The in-degree of each user is the number of other users following the user and the out-degree is the number of others followed by the user. 
The degree distributions in natural and social networks are often heavy-tailed (see Chapter 8 of \citet{Newman:2010ur}).
The resulting Hill plot for the in-degrees of the Google+ data (the first plot in Figure \ref{figure: GPlus 3plots}) resembles the curve of the Hill plot for the Pareto observations with the largest extremes removed. This raises the question of whether some extreme in-degrees of the Google+ data are also unobserved. For example, some users with extremely large in-degrees may have been excluded from the data. 
This pattern of a smooth curve becomes even more pronounced when we apply an additional removal of the top 500 and 1000 values of the in-degree (the second and the third plots in Figure \ref{figure: GPlus 3plots}).

In addition to detecting possible manipulation of data in the tail, modeling and analyzing data in the presence of missing extremes can also be applied to a variety of fields. For example, in studying natural disasters such as earthquakes, forest fires and floods, extreme values might be missing due to difficulty in data collection. In actuarial sciences, claims of extremely large amounts might be covered by a reinsurance company and not included in the claims (Section 8.7 of \citet{Embrechts:1997gj}, \citet{Benchaira:2016dm}).

In order to understand the behavior of the Hill curves of samples in which some of the top extreme values have been removed, we introduce a new parametrization to the Hill estimator. 
Let $k=k_n$ be an intermediate sequence. We denote the number of upper order statistics used in the Hill estimator by $\theta k$ and the number of missing extremes by $\delta k$, and define a functional version of the Hill estimator without extremes (HEWE) as a function of $\theta$ and $\delta$. 
This new parametrization allows one to explore missing extremes both visually and theoretically. The Hill estimator curve of the data without the top extremes exhibits a strikingly smooth and increasing pattern, in contrast to the fluctuating shapes when no extremes are missing. And the differences in the shape of the curves are explained by the functional properties of the limiting process of the HEWE.
Under a second-order regular varying condition, we show that the HEWE, suitably normalized, converges in distribution to a continuous Gaussian process with mean zero and covariance depending on $ \delta $ and parameters of the distribution $ F $ including the tail index $ \alpha $. 

Based on the likelihood function of the limiting process, an estimation procedure is developed for $ \delta $ and the parameters of the distribution, in particular, the tail index $ \alpha $. 
The proposed approach may also have value in assessing the fidelity of the data to the heavy-tailed assumptions. Specifically, one would expect consistency of the estimation of the tail index when more extremes are artificially removed from the data.

A natural question is whether the observed phenomenon, such as those illustrated in the Hill plots in Figure \ref{figure: GPlus 3plots}, is an artifact of the data coming from a light-tailed distribution. In fact, our method is robust to the light-tailed case and can differentiate between the case of heavy-tailed data with missing extremes and light-tailed data. A theoretical justification can be found in \citet{Davis:1984kc}, in which the consistency of the Hill estimator when $ \alpha = \infty $ was established. We also include an example in the simulation section to demonstrate the good performance of the proposed method when applied to light-tailed data.

There has been recent work (\citet{Aban:2006eh}, \citet{Beirlant:2016gs,2016arXiv160602090B}) that involves adapting classical extreme value theory to the case of truncated Pareto distributions. The truncation is modeled via an unknown threshold parameter and the probability of an observation exceeding the threshold is zero. Maximum likelihood estimators (MLE) are derived for the threshold and the tail index. 

Our focus here is to study the path behavior of the HEWE if any arbitrary number of largest values are unavailable. 
Moreover, the estimation procedure we propose has a built-in mechanism to compensate for the bias introduced by non-Pareto heavy-tailed distributions. 
Ultimately, the HEWE provides a graphical and theoretical method for estimation and assessment of modeling assumptions.
An R Shiny web application has been built to interactively estimate and evaluate results from user uploaded data (see the supplementary material for details). 

In addition, we feel the proposed approach may shed some useful insight on classical extreme value theory even when extreme values are not missing in the observed data. 
It is possible to remove a number of top extreme values artificially and study the effect of the artificial removal on the estimation of the tail index. In this case we know the true value of $ \delta $.

This paper is organized as follows. 
Section \ref{sect: conv} introduces the HEWE process and states the main result of this paper dealing with the functional convergence of the HEWE to a continuous Gaussian process. 
Section \ref{section estimation} explains the details of the estimation procedure based on the asymptotic results. 
Section \ref{sect: sim} demonstrates how our estimation procedure works on simulated data from both Pareto and non-Pareto distributions. 
We also illustrate this procedure on a light-tailed distribution.
Section \ref{sect: data} applies our procedure to several interesting real data sets. All the proofs are postponed to the Appendix.

\section{Functional Convergence of HEWE} \label{sect: conv}
In this section we set up the framework for studying the reparametrized Hill estimator. To start,
let $ X_1, X_2, \dots $ be iid random variables with distribution function $ F $ satisfying the regular varying condition Eq. \ref{reg varying}. 
Let $ X_{(1)} \ge X_{(2)} \ge \dots \ge X_{(n)}$ denote the order statistics of $ X_1,\dots, X_n $. 
Let integer $ k \in \{1, \dots, n-1\} $. For fixed $\delta \ge 0$, the HEWE process is defined to be the function of $\theta>0$ given by
\begin{equation} \label{eq: Hill}
H_{k, n}(\theta; \delta) = \begin{cases}
\frac{1}{\lfloor \theta k \rfloor} \sum_{i = 1 }^{\lfloor \theta k \rfloor} \log X_{(\lfloor \delta k \rfloor + i)} - \log X_{(\lfloor \delta k \rfloor + \lfloor \theta k \rfloor + 1)}, & \theta \ge 1/k, \\
0 ,& \theta < 1/k. 
\end{cases}
\end{equation}
Strictly speaking, the process in Eq. \ref{eq: Hill} is defined only when $\lfloor \delta k \rfloor + \lfloor \theta k \rfloor < n$. Asymptotically, we will assume that $k/n \rightarrow 0$, so the process will be defined for all $\delta\ge 0$ and $\theta>0$.

To see the idea behind this definition, imagine that the top $ \lfloor \delta k \rfloor $ observations are not available in the data set and the Hill estimator is computed based on $ \lfloor \theta k \rfloor $ extreme order statistics of the remaining observations. 
Viewed as a function of the observable part of the sample, $ H_{k, n} $ is the usual Hill estimator based on the $\lfloor\theta k \rfloor$ upper order statistics. 
A special case is when $ \delta = 0 $ and no extreme values are missing, then $ H_{k, n}(\theta; \delta = 0) $ corresponds to the usual Hill estimator based on the upper $ \lfloor \theta k \rfloor $ observations. 

Here we treat $ \delta $ as a fixed unknown parameter and Eq. \ref{eq: Hill} a single-parameter process $ H_{k, n}(\theta; \delta) $ indexed by $ \theta $.  $ H_{k, n}(\theta; \delta)$ will play a key role in estimating relevant parameters such as $ \delta $ and $ \alpha $. The estimation is based on the asymptotic distribution of $ H_{k, n}(\theta; \delta) $ and is described in detail in Section \ref{section estimation}. 

In order to obtain the functional convergence of $ H_{k, n}(\theta; \delta) $, a second-order regular variation condition, which provides a rate of convergence in Eq. \ref{reg varying} is needed. This condition can be found, for example, in Section 2.4 of \citet{deHaan:2006bz}, and it states that for $ x > 0 $,
\begin{equation} \label{eq: second order cond}
\lim_{t\rightarrow \infty} \frac{\log U(tx) - \log U(t) - \alpha^{-1} \log x}{A(t)} = \frac{x^\rho - 1}{\rho},  
\end{equation}
where $ \rho \le 0 $, $ U(t) = F^{\leftarrow}(1- 1/t)  $  and $ A $ is a positive or negative function with $ \lim_{t\rightarrow \infty} A(t) = 0 $. Assume
that the sequence $k=k_n\to\infty $ satisfies
\begin{equation} \label{cond lambda}
\lim_{n \rightarrow \infty} \sqrt{k_n}A(n/k_n) = \lambda,
\end{equation}
where $ \lambda $ is a finite constant.
Note conditions Eq. \ref{eq: second order cond} and Eq. \ref{cond lambda} imply that $ n / k_n \rightarrow \infty $ and that $ A $ is a regular-varying function with index $ \rho $.

Distributions that satisfy the second-order condition include the Cauchy, Student's $ t_\nu $, stable, Weibull and extreme value distributions (for more discussion on the second-order condition, see, for example, \citet{Drees:1998ve} and  \citet{Drees:2000ic}). 
In fact, any distribution with $ \bar{F}(x) = c_1 x^{-\alpha} + c_2 x^{-\alpha + \alpha \rho} (1 + o(1)) $ as $ x \rightarrow \infty $, where $ c_1 > 0 $, $ c_2 \ne 0 $, $ \alpha > 0 $ and $ \rho < 0 $, satisfies the second-order condition with the indicated values of $ \alpha $ and $ \rho $ (\citet{deHaan:2006bz}).

Pareto distributions with tail index $ \alpha > 0 $ ($ \bar{F}(x) = x^{-\alpha} $ for $ x \ge 1 $ and zero otherwise), however, do not satisfy the second-order condition, as the numerator on the left side of Eq. \ref{eq: second order cond} is zero when $ t $ is large enough. As will be seen later, the results can be readily extended to the case of Pareto distributions by replacing terms involving $ \rho $ with zero.

We now state the main result of this paper which establishes the functional convergence of the HEWE to a Gaussian process. 

\begin{theorem} \label{thm: convergence}
\normalfont
Assume the second-order condition Eq. \ref{eq: second order cond} holds and Eq. \ref{cond lambda} is satisfied for a given sequence $ k_n $ and $ \lambda $. Then\\
\textbf{(a)} there exist versions of $ H_{k, n}(\theta; \delta) $ and a standard Brownian Motion $ W $ defined on the same probability space such that as $ n \rightarrow \infty $,
\begin{multline}
H_{k, n}(\theta; \delta) 
 = \frac{ g_\delta(\theta)}{\alpha} + \frac{1}{\alpha} \frac{1}{\theta \sqrt{k_n}} \int_{\delta}^{\delta + \theta} \Big(1 - \frac{\delta}{x}\Big) dW(x)  
 +  A\Big(\frac{n}{ k_n }\Big) b_{\delta, \rho}(\theta)
 + o\Big(\frac{1} {\sqrt{k_n} }\Big), 
 \quad \text{a.s.}
\end{multline}
holds uniformly in $ (\theta, \delta) $ on compact subsets of $ (0, \infty) \times [0, \infty) $,
where
\begin{equation*}
g_\delta(\theta) = 
\begin{cases}
1, & \delta = 0, \\
1 - \frac{\delta}{\theta} \log\big(\frac{\theta}{\delta} + 1 \big), & \delta > 0, 
\end{cases}
\end{equation*}
\begin{equation*}
b_{\delta, \rho}(\theta) = 
\begin{cases}
\frac{ 1}{1 - \rho} \frac{1}{\theta^\rho}, & \delta = 0, \\
\frac{1 + (\theta/\delta) \rho - (\theta/\delta + 1)^\rho}{(\theta/\delta) (1-\rho)\rho} \frac{1}{(\delta + \theta)^\rho}, & \delta > 0.
\end{cases}
\end{equation*}
\textbf{(b)} For all $\delta \ge 0$,
\begin{equation*}
\sqrt{k_n} \Big(H_{k, n}(\cdot; \delta) - \alpha^{-1} g_\delta(\cdot)\Big) - \lambda b_{\delta, \rho}(\cdot) 
\stackrel{d}{\rightarrow} \alpha^{-1} G_\delta(\cdot)
\end{equation*}
in $ D(0, \infty) $, 
where
\begin{equation}
G_\delta(\theta) = \frac{1}{\theta}\int_{\delta}^{\delta + \theta} \Big(1 - \frac{\delta}{x} \Big) dW(x) 
\end{equation} 
is a Gaussian process with mean zero and covariance function
\begin{equation*}
\text{Cov}\big( G_\delta(\theta_1), G_\delta(\theta_2) \big) = 
\begin{cases}
\frac{1}{\theta_1 \theta_2} \bigg[
\theta_1 \wedge  \theta_2 
- 2 \delta \log \Big(1 + \frac{ \theta_1 \wedge  \theta_2}{\delta} \Big) + \frac{\delta (\theta_1 \wedge \theta_2) } {\delta + (\theta_1 \wedge \theta_2)} \bigg], & \delta > 0\\
\frac{1}{\theta_1 \vee \theta_2}, &\delta = 0.
\end{cases}
\end{equation*}
\end{theorem}

\begin{remark}
\normalfont
Theorem \ref{thm: convergence} states the weak convergence of $ H_{k,n}(\cdot; \delta) $ for all fixed $ \delta \ge 0 $. In fact, we have shown a stronger result (see Appendix) on the weak convergence of $H_{k,n}(\theta, \delta):=H_{k,n}(\theta; \delta)$ viewed now as a random field indexed by the pair $(\theta, \delta)$:
\begin{equation} \label{stronger result}
\sqrt{k_n} \bigg(H_{k, n}(\cdot, \cdot) - \alpha^{-1} \tilde{g}(\cdot, \cdot) \bigg) - \lambda  \tilde{b}_{\rho}(\cdot, \cdot) \stackrel{d}{\rightarrow} \alpha^{-1}\tilde{G}(\cdot, \cdot)
\end{equation}
in $ D((0, \infty)\times [0, \infty)) $, 
where $\tilde g(\theta,\delta)=g_\delta(\theta)$, $\tilde b_\rho(\theta,\delta)=b_{\delta,\rho}(\theta)$,
and $ \tilde{G}(\theta, \delta) = (1/\theta)\int_{\delta}^{\delta + \theta} (1 - \delta /x) dW(x)$ with mean zero and the following covariance function. If $ \delta_1 \vee \delta_2 > 0 $, 
\begin{align*}
& \text{Cov}\big( \tilde{G}(\theta_1, \delta_1), \tilde{G}(\theta_2, \delta_2) \big) \nonumber \\ 
& = \frac{1}{\theta_1 \theta_2} \bigg[
(\delta_1 + \theta_1) \wedge (\delta_2 + \theta_2) - (\delta_1 \vee \delta_2) \nonumber \\
& \qquad - (\delta_1 + \delta_2) \log \bigg(\frac{(\delta_1 + \theta_1) \wedge (\delta_2 + \theta_2)}{\delta_1 \vee \delta_2} \bigg) + \frac{\delta_1 \delta_2}{\delta_1 \vee \delta_2} - \frac{\delta_1 \delta_2} {(\delta_1 + \theta_1) \wedge (\delta_2 + \theta_2)} \bigg].
\end{align*}
If $ \delta_1 = \delta_2 = 0 $,
\begin{equation*}
\text{Cov}\big( \tilde{G}(\theta_1, 0), \tilde{G}(\theta_2, 0) \big) = \frac{1}{\theta_1 \vee \theta_2}.
\end{equation*}
\end{remark}

\begin{remark}
\normalfont
For fixed $ \theta $, the functions $ g_\delta $ and $b_{\delta, \rho} $ are continuous at $ \delta = 0 $. 
For iid Pareto variables $ X_1, X_2, \dots$ with tail index $ \alpha > 0 $, the result of Theorem \ref{thm: convergence} still holds with the bias term $ b_{\delta, \rho} $ replaced by zero.
\end{remark}

\begin{figure}[H] 
\centering
\begin{minipage}{.48\textwidth}
  \centering
  \includegraphics[width=0.9\linewidth]{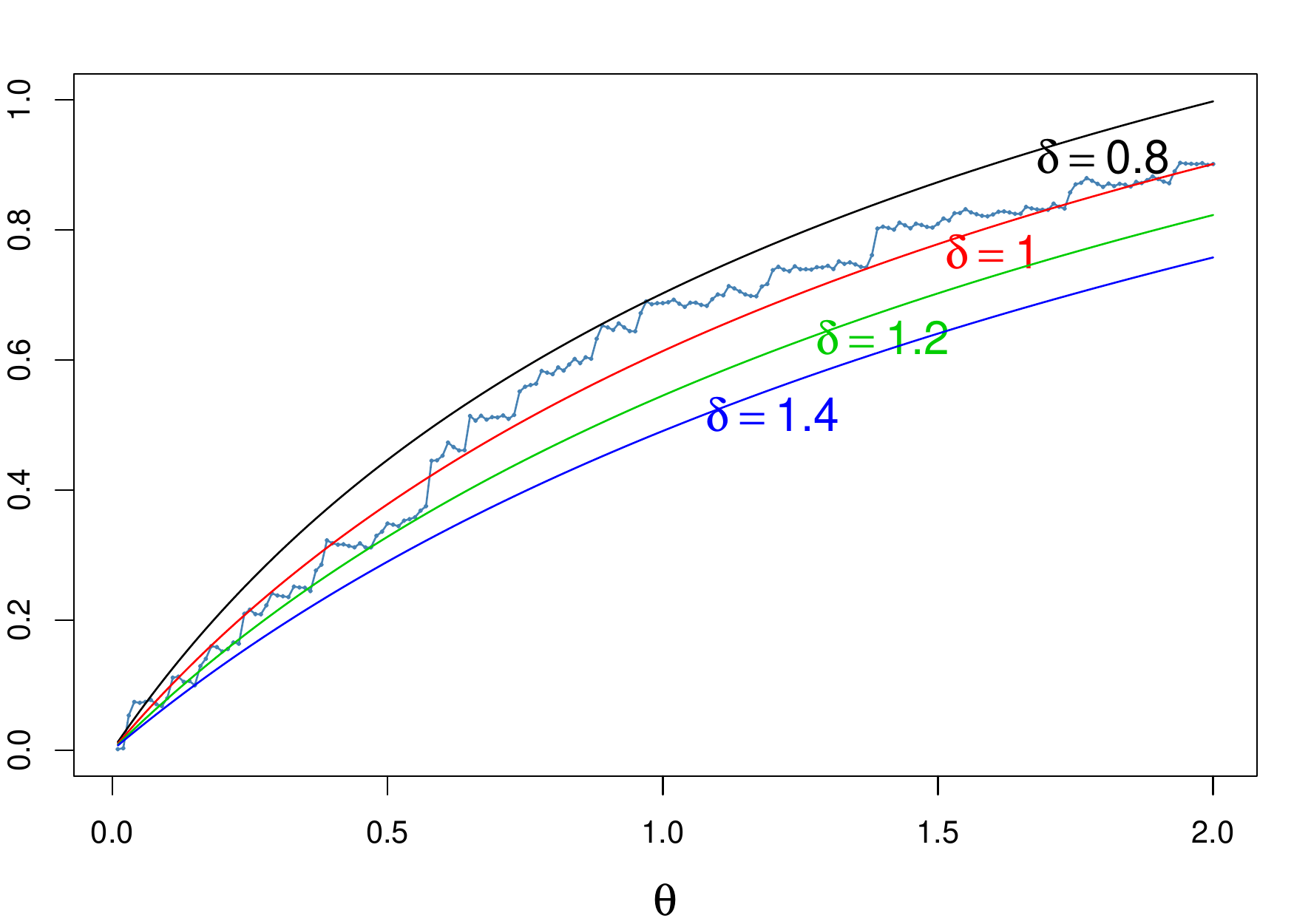} 
\end{minipage}
\hfill
\begin{minipage}{.48\textwidth}
  \centering
  \includegraphics[width=0.9\linewidth]{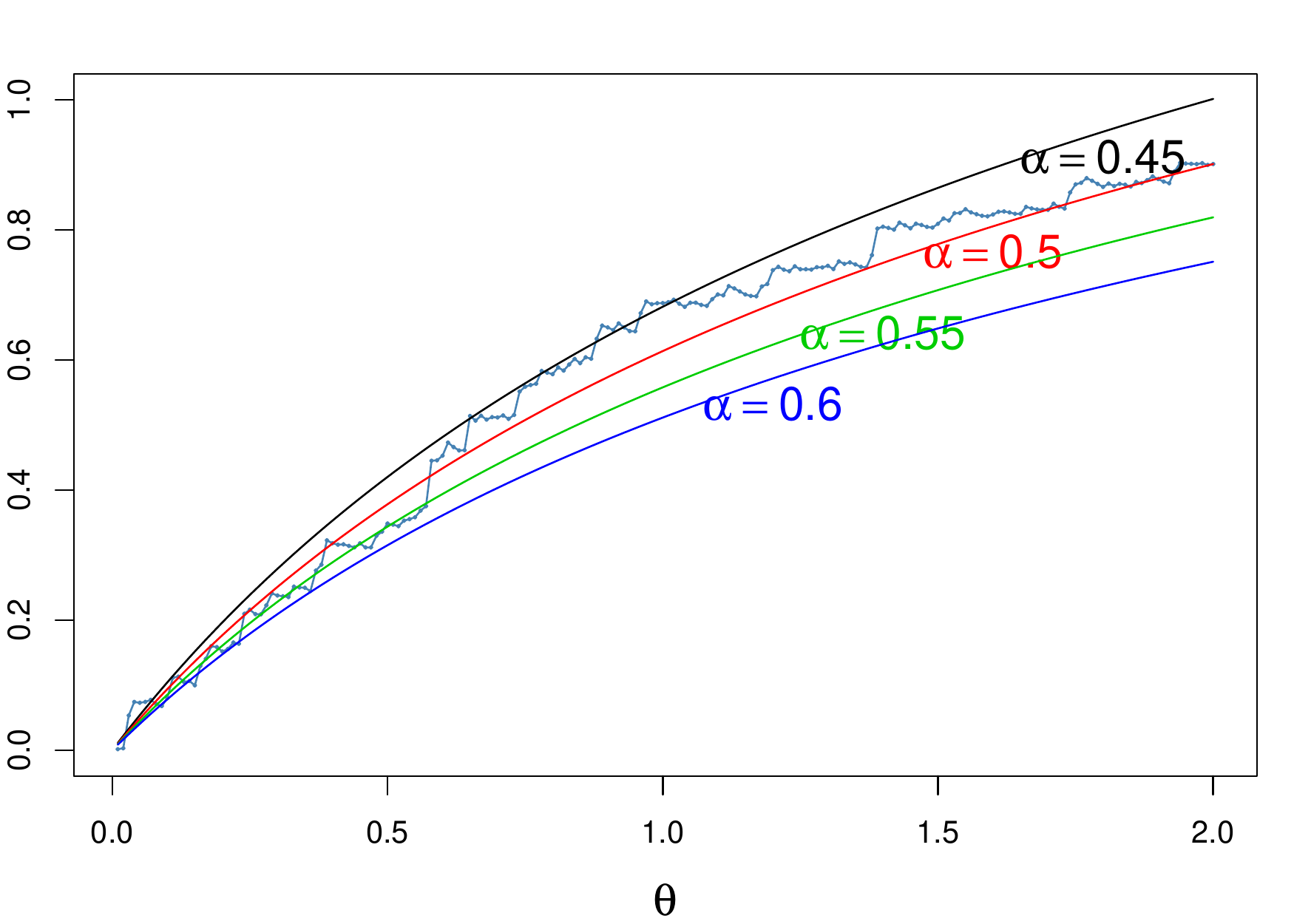} 
\end{minipage} 
\caption{Fitting mean curves with different values of parameters to the Hill plot for the Pareto sample as in Figure \ref{hill_pareto}. Left: fixing $ \alpha = 0.5 $. Right: fixing $ \delta = 1 $} \label{curves: fix A and D}
\end{figure}

Figure \ref{curves: fix A and D} shows the Hill estimates of the same sample from the Pareto distribution with $ \alpha = 0.5 $ as in Figure \ref{hill_pareto} overlaid with several mean curves. We chose $ k_n = 100 $ with the top $ 100 $ observations removed from the original sample. This implies $ \delta = 1 $.
In the left panel of Figure \ref{curves: fix A and D}, the Hill estimates are overlaid with the mean curves $ g_\delta(\theta) / \alpha $ of the Gaussian process with different values of $\delta $ while fixing the true value of $ \alpha = 0.5 $. 
The right panel of Figure \ref{curves: fix A and D} shows the mean curves with different values of $ \alpha $ while fixing the true value $ \delta = 1 $. In both plots, the Hill plot is closest to the mean curve corresponding to the true value of the parameter.

\begin{figure}[H]
\centering
\begin{minipage}{.48\textwidth}
  \includegraphics[width=0.9\textwidth]{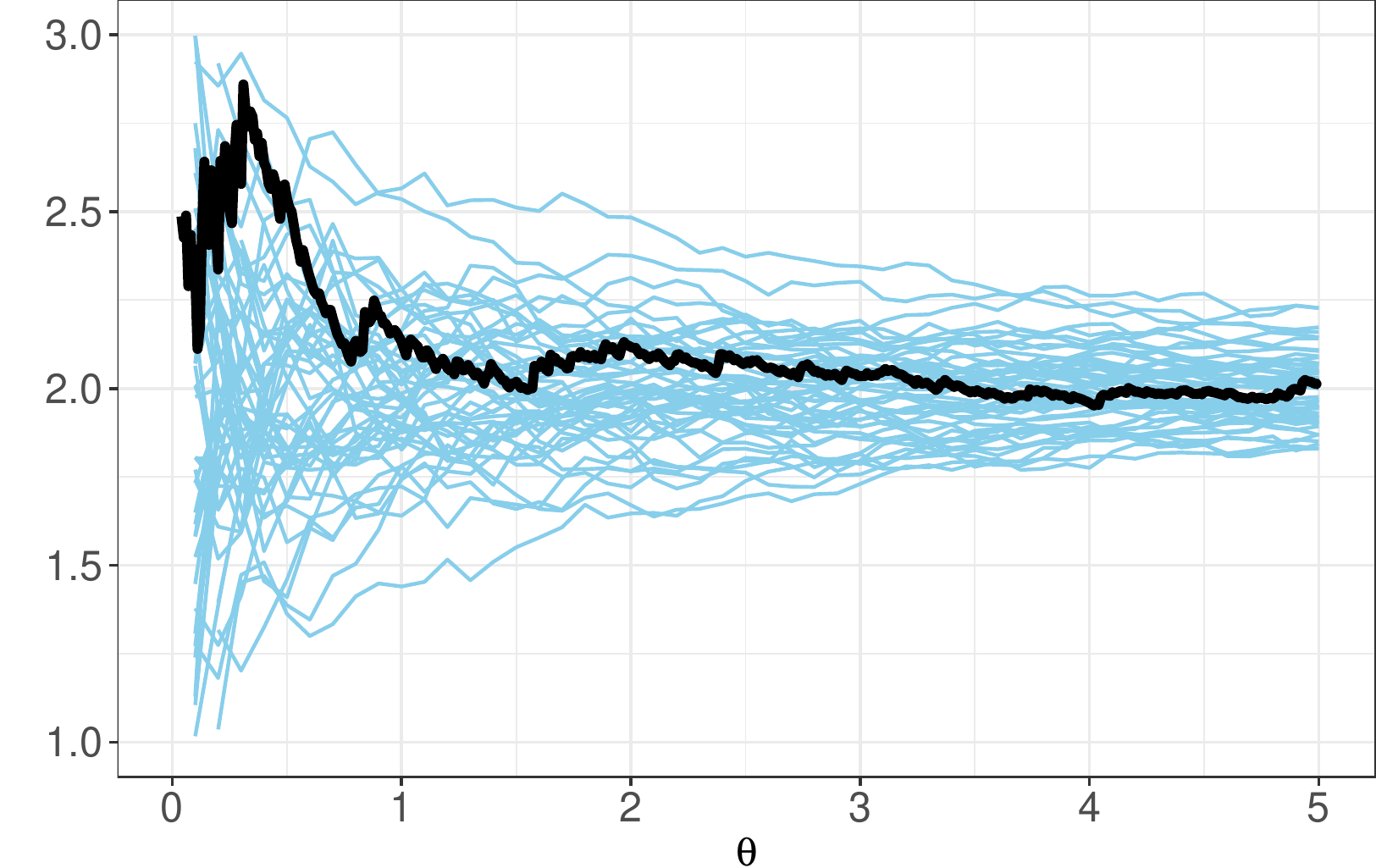}
\end{minipage}
\hfill
\begin{minipage}{.48\textwidth}
  \includegraphics[width=0.9\textwidth]{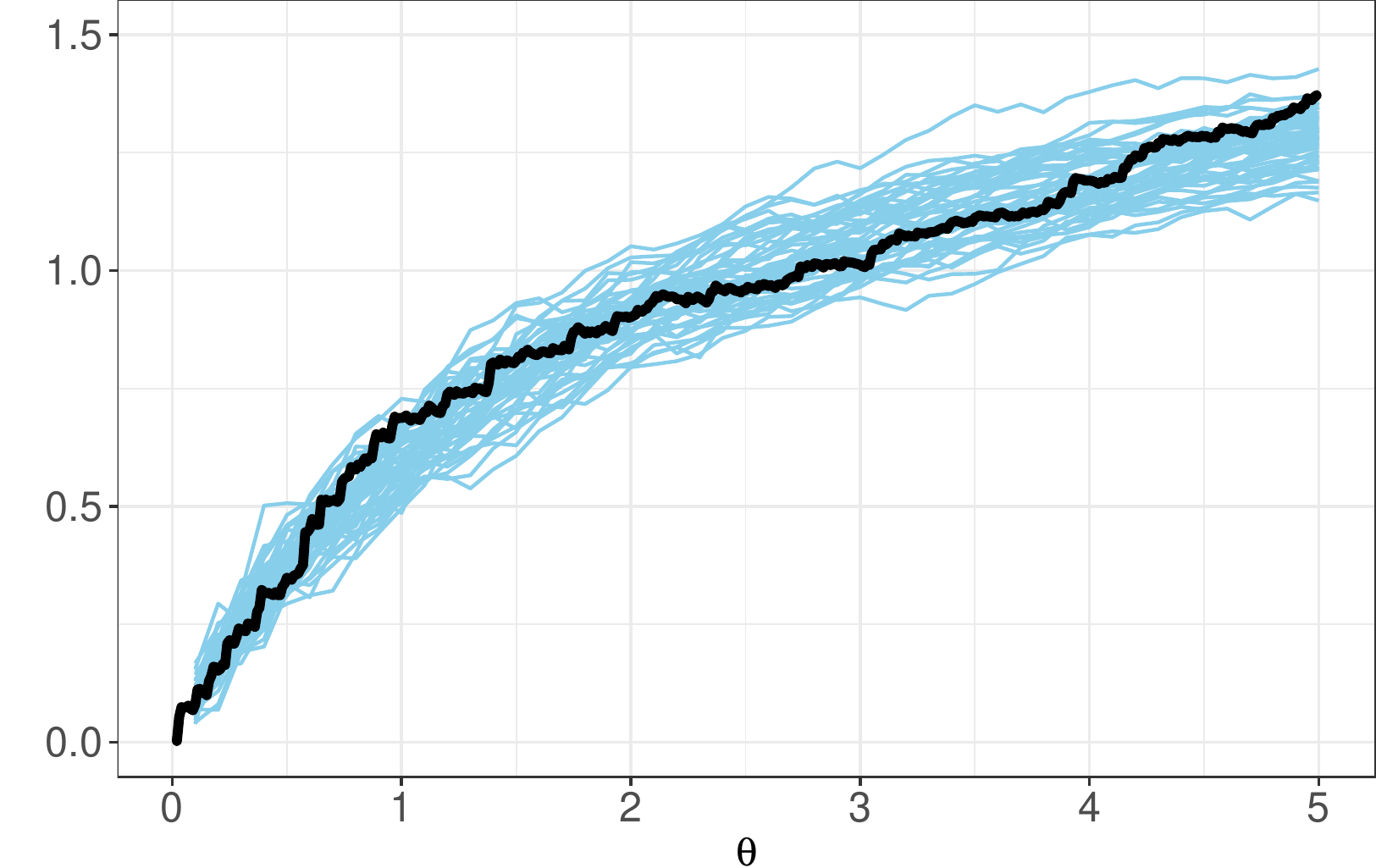}
\end{minipage}
\caption{Observed Hill plots for the Pareto sample (bold lines) and realizations from corresponding Gaussian processes (thin lines). Left: with the original sample. Right: top $ 100 $ extreme values removed} \label{gaussian_pareto}
\end{figure}

\begin{figure}[H]
\centering
\begin{minipage}{.48\textwidth}
  \includegraphics[width=0.9\textwidth]{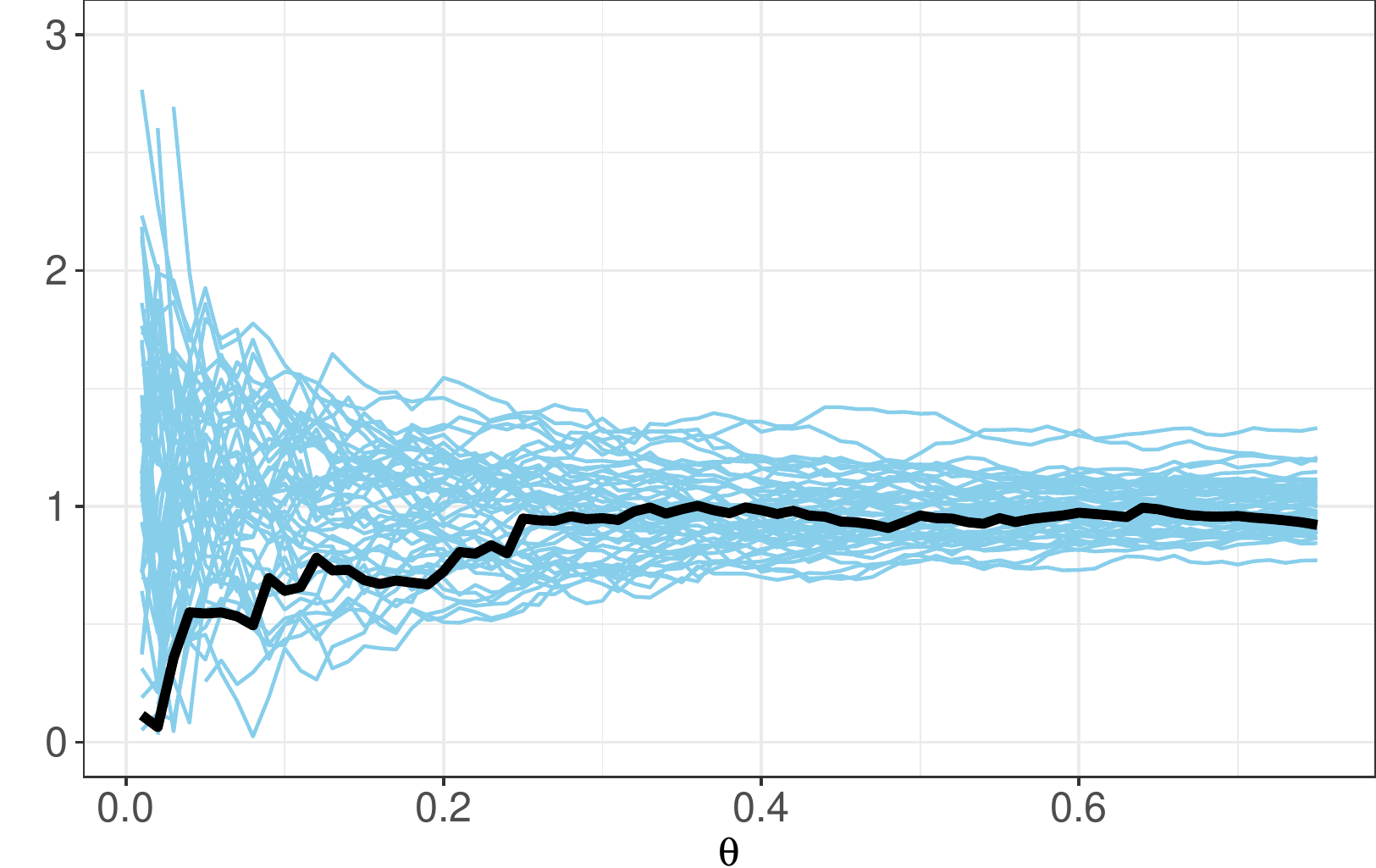}
\end{minipage}
\hfill
\begin{minipage}{.48\textwidth}
  \includegraphics[width=0.9\textwidth]{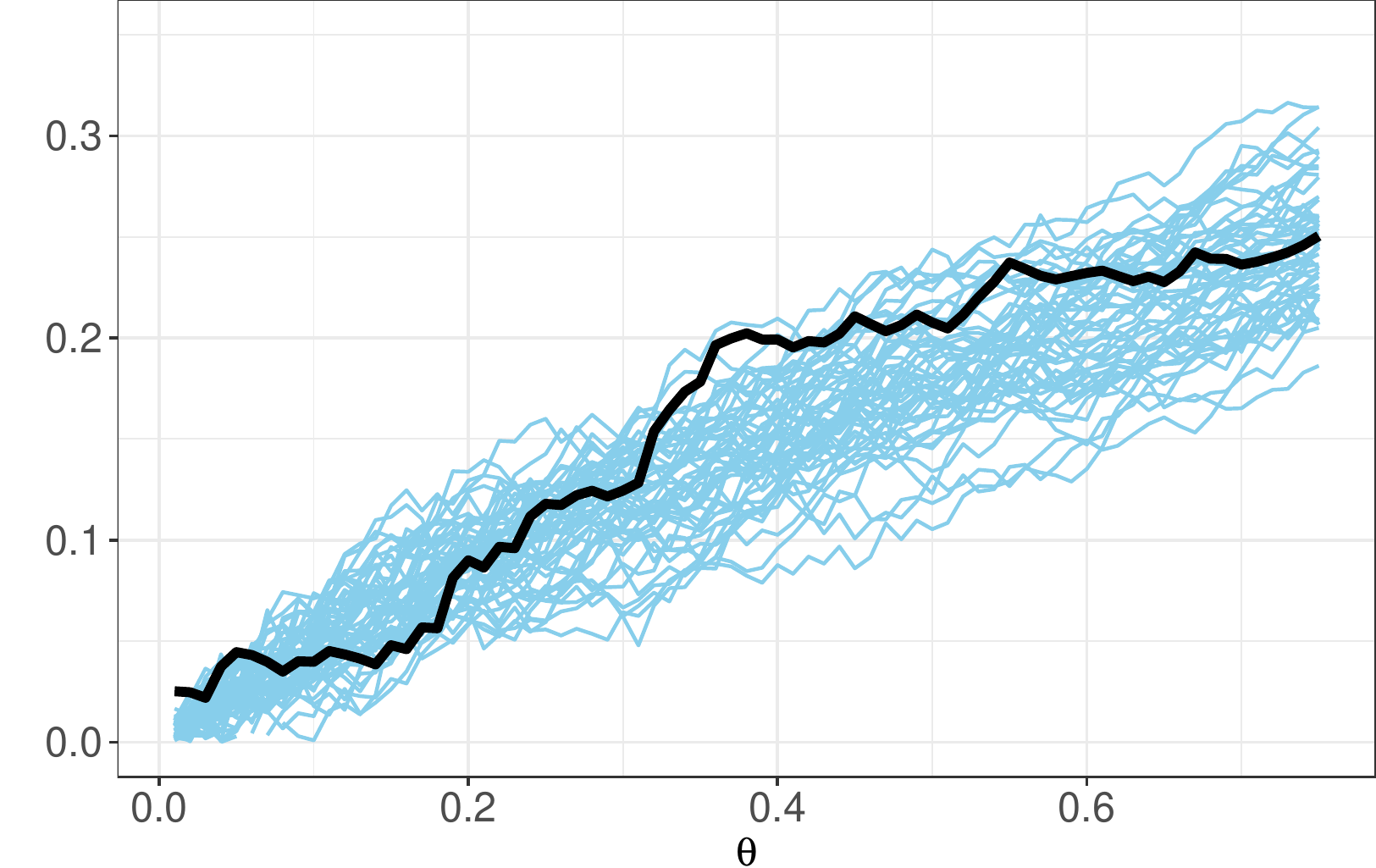}
\end{minipage}
\caption{Observed Hill plots for a Cauchy sample (bold lines) and realizations from corresponding Gaussian processes (thin lines). Left: with the original sample. Right: top $ 100 $ extreme values removed} \label{gaussian_cauchy}
\end{figure}

In order to demonstrate the variability generated by the limiting Gaussian process, we compare the Hill plots for samples from Pareto and Cauchy distributions with their Gaussian process approximations given by Theorem \ref{thm: convergence}.
Figure \ref{gaussian_pareto} presents the Hill plots for the same Pareto sample as in Figures \ref{hill_pareto} and \ref{curves: fix A and D}, without removal of extremes (left) and with the top $ 100 $ observations removed (right), along with $ 50 $ independent realizations from the corresponding Gaussian processes with bias $ b_{\delta, \rho} \equiv 0 $.
\noindent
Figure \ref{gaussian_cauchy} shows the Hill plots for a Cauchy sample ($n = 1000$, $ k_n = 100 $, $ \alpha = 1 $ and $ \rho = -2$), without removal of extremes and with the top $ 100 $ extremes removed, along with $ 50 $ independent realizations from the corresponding Gaussian processes with non-zero $ b_\rho $.

\section{Parameter Estimation} \label{section estimation}
Let $ X_1, X_2, \dots X_n $ be a sample from a distribution $ F $ satisfying the second-order regular variation condition Eq. \ref{eq: second order cond}, and let $ X_{(1)} \ge X_{(2)} \ge \dots \ge X_{(n)} $ denote the order statistics of $ \{X_i\} $. 
Suppose the $\lfloor \delta k_n \rfloor$ largest observations are unobserved in the data and the value of $ \delta $ is unknown.

In this section, we develop an approximate maximum likelihood estimation procedure for the unknown parameters $ \delta $, $ \alpha $ and $ \rho $ given the observed data. The procedure is based on the asymptotic distribution of $ H_{k,n}(\theta; \delta) $.
For typographical convenience we suppress the dependence of $\delta$ and $k_n$ and use the notation $H_n(\theta)$.

By Theorem \ref{thm: convergence}, the joint distribution of $(H_n(\theta_1), \dots, H_n(\theta_s))$ for fixed $ (\theta_1,\ldots, \theta_s)$ can be approximated, when $k_n$ is large, by a distribution with density function at ${\bf h}=(h_1,\ldots, h_s)$ given by 
\begin{equation} \label{eq: asym joint distr}
\frac{1}{\sqrt{(2\pi)^{s} |\mathbf{\Sigma}_{\alpha,\delta}|}} \exp \bigg[ -\frac{1}{2} \bigg(\mathbf{h} - \frac{\mathbf{g}_\delta}{\alpha} - \frac{\lambda \mathbf{b}_{\delta, \rho}}{ \sqrt{k_n}} \bigg)^\top \mathbf{\Sigma}_{\alpha, \delta}^{-1} \bigg(\mathbf{h} - \frac{\mathbf{g}_\delta}{\alpha}  - \frac{\lambda \mathbf{b}_{\delta, \rho}}{ \sqrt{k_n}} \bigg) \bigg],
\end{equation}
where
\begin{equation*}
\{\mathbf{g}_\delta\}_i = 
\begin{cases}
1, & \delta = 0, \\
1 - \frac{\delta}{\theta_i} \log\big(\frac{\theta_i}{\delta} + 1 \big), & \delta > 0, 
\end{cases}
\end{equation*}
\begin{equation*}
\{\mathbf{b}_{\delta, \rho}\}_i = 
\begin{cases}
\frac{ 1}{1 - \rho} \frac{1}{\theta_i^\rho}, & \delta = 0, \\
\frac{1 + (\theta_i/\delta) \rho - (\theta_i/\delta + 1)^\rho}{(\theta_i/\delta) (1-\rho)\rho} \frac{ 1}{(\delta + \theta_i)^\rho}, & \delta > 0, 
\end{cases}
\end{equation*}
and
\begin{equation*}
\mathbf{\Sigma}_{\alpha, \delta} (i,j) =
\begin{cases}
\frac{1}{\alpha^2 k_n}  \frac{1}{\theta_i \vee \theta_j}, & \delta = 0, \\
\frac{1}{\alpha^2 k_n} \frac{({\theta}_i \wedge {\theta}_j)^2}{\delta \theta_i \theta_j} \,v \big( \frac{\theta_i \wedge \theta_j}{\delta} \big), & \delta > 0, 
\end{cases}
\end{equation*}
with
\begin{equation*}
v(\theta) = \frac{1}{\theta} - \frac{2 \log(\theta + 1)}{\theta^2} + \frac{1}{\theta(\theta + 1)}.
\end{equation*}

To further simplify the calculation for the maximum likelihood estimator of $ \alpha $, $ \delta $ and $ \rho $, let
\begin{equation} \label{eq: H to T}
T_i = H_n(\theta_i) - \frac{\theta_{i-1}}{\theta_{i}} H_n(\theta_{i-1}),
\end{equation}
where $ \theta_0 = 0 $ is introduced for convenience.
Note that the $ T_i $ are asymptotically independent with the joint density function at $ \mathbf{t} = (t_1, \dots, t_s) $ being
\begin{equation} \label{eq: asym density T}
\frac{1}{\sqrt{(2\pi)^{s} |\tilde{\mathbf{\Sigma}}_{\alpha,\delta}|}} \exp \bigg[ -\frac{1}{2} \big(\mathbf{t} - \mathbf{m} \big)^\top \tilde{\mathbf{\Sigma}}_{\alpha, \delta}^{-1} \big(\mathbf{t} -  \mathbf{m} \big) \bigg],
\end{equation}
where 
\begin{equation*}
m_i =
\frac{1}{\alpha} \Big(\{\mathbf{g}_\delta\}_i - \frac{\theta_{i-1}}{\theta_{i}} \{\mathbf{g}_\delta\}_{i-1} \Big) + \frac{\lambda}{\sqrt{k_n}} \Big(\{\mathbf{b}_{\delta, \rho}\}_i - \frac{\theta_{i-1}}{\theta_{i}} \{\mathbf{b}_{\delta, \rho}\}_{i-1} \Big)
\end{equation*}
and $ \tilde{\mathbf{\Sigma}}_{\alpha, \delta} $ is a diagonal matrix, in which 
\begin{equation*}
\tilde{\mathbf{\Sigma}}_{\alpha, \delta} (i, i) = 
\begin{cases}
\frac{1}{\alpha^2 k_n} \big( \frac{1}{\theta_i} - \frac{\theta_{i-1}}{\theta_{i}^2} \big), 
& \delta = 0, \\ 
\frac{1}{\alpha^2 k_n \delta} \big( v \big(\frac{\theta_i}{\delta} \big) - \big( \frac{\theta_{i-1}}{\theta_i} \big)^2 v \big(\frac{\theta_{i-1}}{\delta} \big) \big),  
&  \delta > 0.
\end{cases}
\end{equation*}
The log-likelihood corresponding to the density Eq. \ref{eq: asym density T} is
\begin{equation} \label{eq: loglik T_new}
C + s \log(\alpha) + \frac{1}{2} \sum_{i=1}^s \log(w_i) - \frac{1}{2} \alpha^2 k_n \sum_{i=1}^s w_i (t_i - m_i)^2,
\end{equation}
where $ C $ is a constant independent of $ \alpha $, $ \delta $ and $ \rho $. 
For $ \delta > 0 $, 
\begin{equation*}
w_i = \delta \bigg/ \bigg(v \Big(\frac{\theta_i}{\delta} \Big) - \Big( \frac{\theta_{i-1}}{\theta_i} \Big)^2 v \Big(\frac{\theta_{i-1}}{\delta} \Big) \bigg).
\end{equation*}
For $ \delta = 0 $,
\begin{equation*}
w_i = 1 \bigg/ \bigg( \frac{1}{\theta_i} - \frac{\theta_{i-1}}{\theta_i^2} \bigg).
\end{equation*}

For fixed $\alpha$ and $\delta$, the only part of the log-likelihood Eq. \ref{eq: loglik T_new} that needs to be optimized is the weighted sum of squares
\begin{equation} \label{eq: weighted sum}
\sum_{i=1}^s w_i(t_i-m_i)^2,                
\end{equation}
and it is minimized over the values of $\rho$ and $\lambda$. 
Note the value of $ \lambda $ depends on the choice of $ k_n $ through Eq. \ref{cond lambda}. When $ k_n $ is fixed, $ \lambda $ is viewed as an independent nuisance parameter and appears in $ m_i $ via
\begin{equation*}
\frac{\lambda}{\sqrt{k_n}}\left( \{ {\bf b}_{\delta,\rho}\}_i -
    \frac{\theta_i}{\theta_{i-1}} \{ {\bf
      b}_{\delta,\rho}\}_{i-1}\right), 
\end{equation*}
which we denote by $\lambda \{ {\bf f}_{\delta,\rho}\}_i / \sqrt{k_n}$, where 
\begin{equation*}
\{\mathbf{f}_{\delta, \rho} \}_i = 
\begin{cases}
\frac{1}{1 - \rho} \frac{1}{\theta_i^\rho} - \frac{\theta_{i-1}}{\theta_{i} } \frac{1}{1 - \rho} \frac{1}{\theta_{i-1}^\rho}, & \delta = 0, \\
\frac{1 + (\theta_i/\delta) \rho - (\theta_i/\delta + 1)^\rho}{(\theta_i/\delta) (1-\rho)\rho} \frac{1}{(\delta + \theta_i)^\rho} 
- \frac{\theta_{i-1}}{\theta_{i} } 
\frac{1 + (\theta_{i-1}/\delta) \rho - (\theta_{i-1}/\delta + 1)^\rho}{(\theta_{i-1}/\delta) (1-\rho)\rho} \frac{1}{(\delta + \theta_{i-1})^\rho} , & \delta > 0. 
\end{cases}
\end{equation*}
Minimizing Eq. \ref{eq: weighted sum} over $\lambda$ and $\rho$ results in 
\begin{equation*}
\hat{\rho}_{\alpha, \delta} = \arg \min_{\rho \le 0} 
\sum_{i=1}^s w_i \bigg(t_i - \frac{1}{\alpha} \Big(\{\mathbf{g}_\delta\}_i - \frac{\theta_{i-1}}{\theta_{i}} \{\mathbf{g}_\delta\}_{i-1}\Big) - \frac{\hat{\lambda}_{\alpha, \delta, \rho}}{\sqrt{k_n}} \{\mathbf{f}_{\delta, \rho} \}_i  \bigg)^2,
\end{equation*}
where
\begin{align} \label{eq: solve lambda}
\hat{\lambda}_{\alpha, \delta, \rho} 
= \sqrt{k_n} \, \frac{\sum_{i = 1}^s w_i 
\big(t_i - (\{\mathbf{g}_\delta\}_i - \frac{\theta_{i-1}}{\theta_{i}} \{\mathbf{g}_\delta\}_{i-1})/\alpha \big) 
\{\mathbf{f}_{\delta, \rho} \}_i } 
{\sum_{i = 1}^s w_i \{\mathbf{f}_{\delta, \rho} \}_i^2}.
\end{align}
Note that this estimation approach, in which $ \lambda $ is viewed as a nuisance parameter, adjusts for the choice of $ k_n $ automatically. If a different $ k_n $ is selected, the estimate of $ \lambda $ will adapt to reflect this change.

Once we have found the optimal values of $\rho$ and $\lambda$, we optimize the resulting expression in Eq. \ref{eq: loglik T_new} by examining its values on a selected grid of $ (\alpha, \delta) $. Alternatively, an iterative procedure can be used, where in each step one of $ \alpha, \delta, \rho $ is updated given values of the other two parameters until convergence of the log-likelihood function. Details on the implementation of the optimization algorithm are described in Section 4.1.

\section{Simulation Studies} \label{sect: sim}
In this section we test our procedure on simulated data. 
In each of the following simulations, we generate $ 200 $ independent samples of size $ n $ from a regular-varying distribution function with tail index $ \alpha $. Given a $ k_n $, we remove the largest $ \lfloor \delta k_n \rfloor$ observations from each of the original samples and apply the proposed method to the samples after the removal. 

For comparison, we also apply the method in \citet{Beirlant:2016gs} to the same samples. In \citet{Beirlant:2016gs}, $ \alpha $ and the threshold $ T $ over which the observations are discarded are estimated with the MLE based on the truncated Pareto distribution. The odds ratio of the truncated observations under the un-truncated Pareto distribution is estimated by solving an equation involving the estimates of $ \alpha $ and $ T $. Finally, the number of truncated observations is calculated given the odds ratio and the observed sample size.
  
For each combination of distribution and parameters, we start from $ \theta_1 = 5 / k_n $ and let $ \theta_i = \theta_{i-1} + 1/k_n $ for $ 1 < i \le s $. 
We consider a sequence of different endpoints $ \theta_s k_n $ to examine the influence of the range of order statistics included in the estimation.
For each value of  $ \theta_s $, we solve for the estimates of $ \alpha $ and $ \delta $ based on the asymptotic density of $ ( H_n(\theta_1), \dots, H_n(\theta_s) ) $ following the procedure described in Section \ref{section estimation}.

Simulations from both Pareto and non-Pareto distributions show that the proposed method provides reliable estimates of the tail index and performs particularly well in estimating the number of missing extremes. The advantages of the proposed method become more apparent in dealing with non-Pareto samples.

\subsection{Pareto Distribution}

First we examine Pareto samples with $n = 500$ and $ \alpha = 0.5 $. Let $ k_n = 50 $ and $ \delta = 1 $ so that $ \delta k_n = 50 $ top extreme observations are removed from the original data.

We apply the estimation procedure introduced in Section \ref{section estimation} to the Pareto data. A series of different values of $\theta_s$ are selected and for each fixed $\theta_s$ the estimation is based on the largest $\theta_s k_n$ values in the data. 
First we calculate the Hill estimates $(H_n(\theta_1), \dots, H_n(\theta_s))$, whose joint distribution is given by Eq. \ref{eq: asym joint distr}. 
To simplify the maximum likelihood estimation, we further transform the calculated $ \{H_n(\theta_i) \}$ to the series $\{T_i \}$ via Eq. \ref{eq: H to T}, which  has joint distribution Eq. \ref{eq: asym density T}. 
The unknown parameters in the log-likelihood are $\alpha, \delta, \rho$ and a nuisance parameter $\lambda$. The parameters are estimated following a two-step procedure; 
for each pair of fixed $\alpha$ and $\delta$, the optimization of the log-likelihood can be further reduced to the optimization of the weighted sum of squares Eq. \ref{eq: weighted sum}. 
For each value of $\rho$, the solution of the optimal $\lambda$ has an explicit form Eq. \ref{eq: solve lambda} involving $\rho$, so that the weighted sum of squares becomes a function of $\rho$ only and can be optimized readily with existing optimization algorithms for continuous functions. 
As the first step of the optimization, we find the optimal $\rho$ with the function optimize() in R 3.4.0.
In the second step, we search for the optimal values of $\alpha$ and $\delta$ on a selected grid of values. While the precision of the estimation depends on the fineness of the selected grid, upon experimenting with different sizes of the grid, we observe the optimization is generally robust and did not appear to be trapped in local maxima. For demonstrative purpose, in all examples of Section \ref{sect: sim},  the fineness of the grid of $\alpha$ is on the scale of $0.01$ and the fineness of the grid of $\delta$ is on the scale of  $0.001$.

\begin{figure}[H]
  \centering
  \begin{minipage}[b]{0.45\textwidth}
  \centering
    \includegraphics[width=0.95\textwidth]{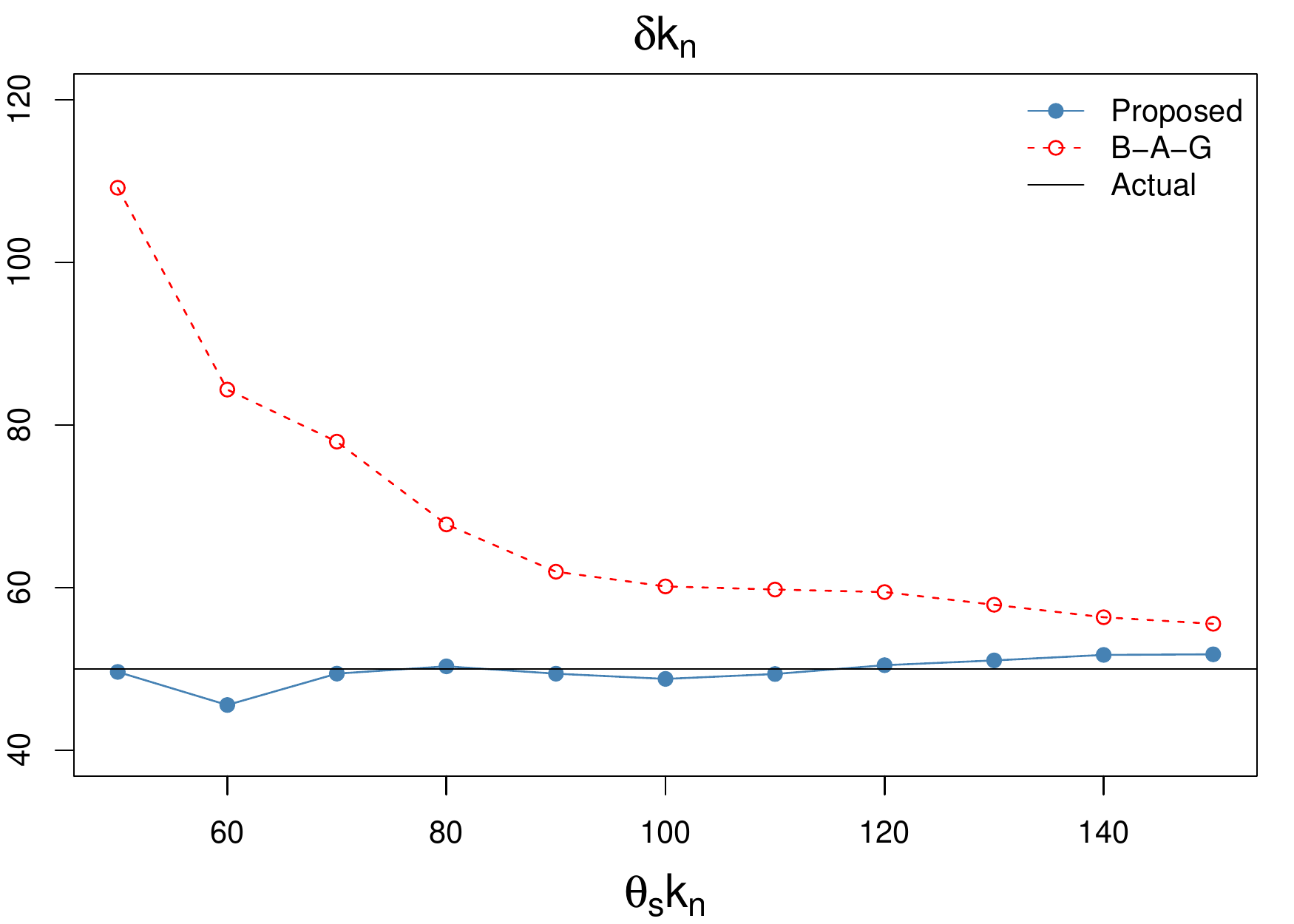}
  \end{minipage}
  \hfill
  \begin{minipage}[b]{0.45\textwidth}
  \centering
    \includegraphics[width=0.95\textwidth]{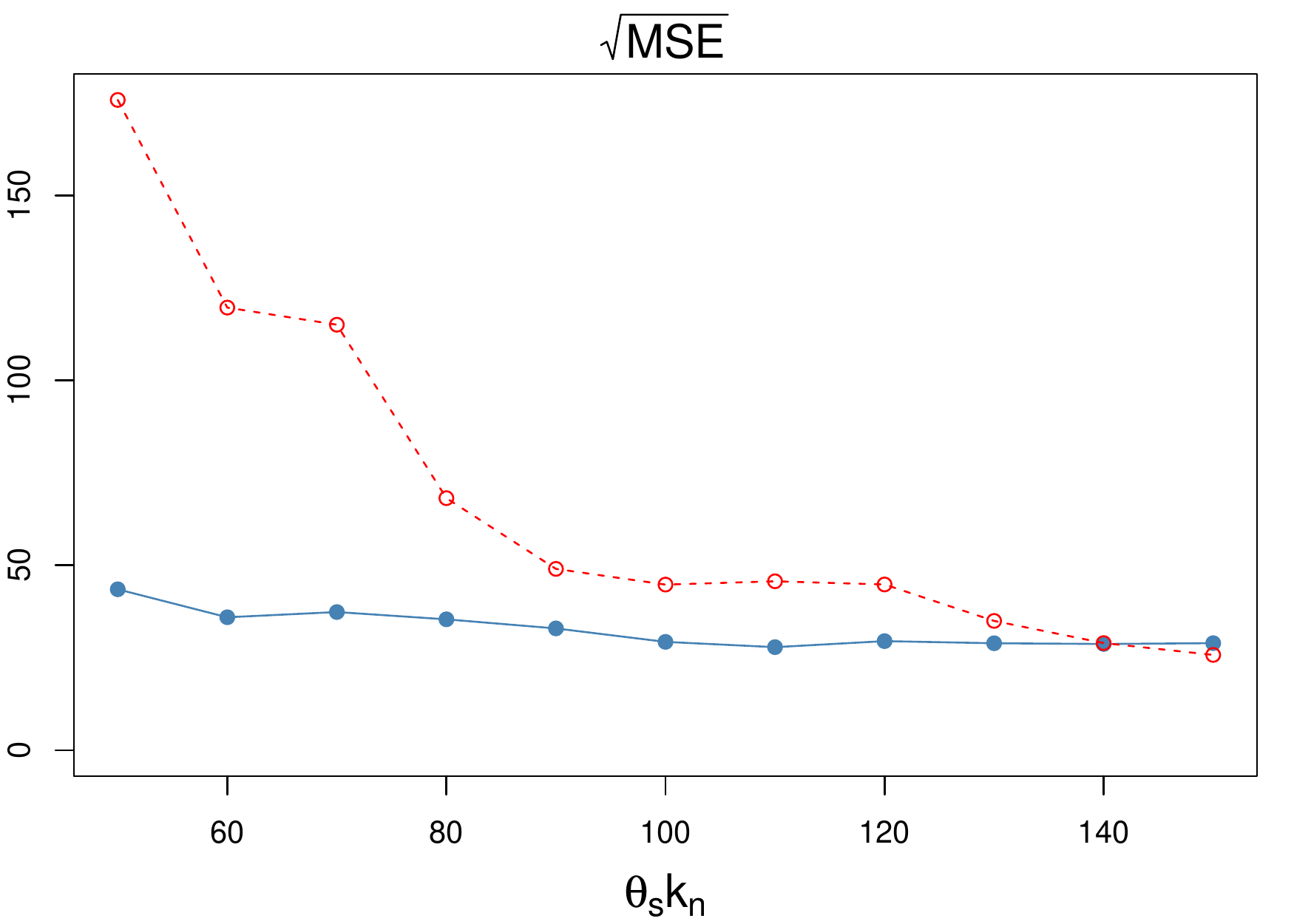}
  \end{minipage}
  \caption{Estimated number of missing extremes and $ \sqrt{\text{MSE}} $ for Pareto samples. $n = 500$, $ \alpha = 0.5 $, $ k_n = 50 $, $ \delta = 1 $} \label{figure: tr pareto delta 1}
\end{figure}

\begin{figure}[H]
  \centering
  \begin{minipage}[b]{0.45\textwidth}
  \centering
    \includegraphics[width=0.95\textwidth]{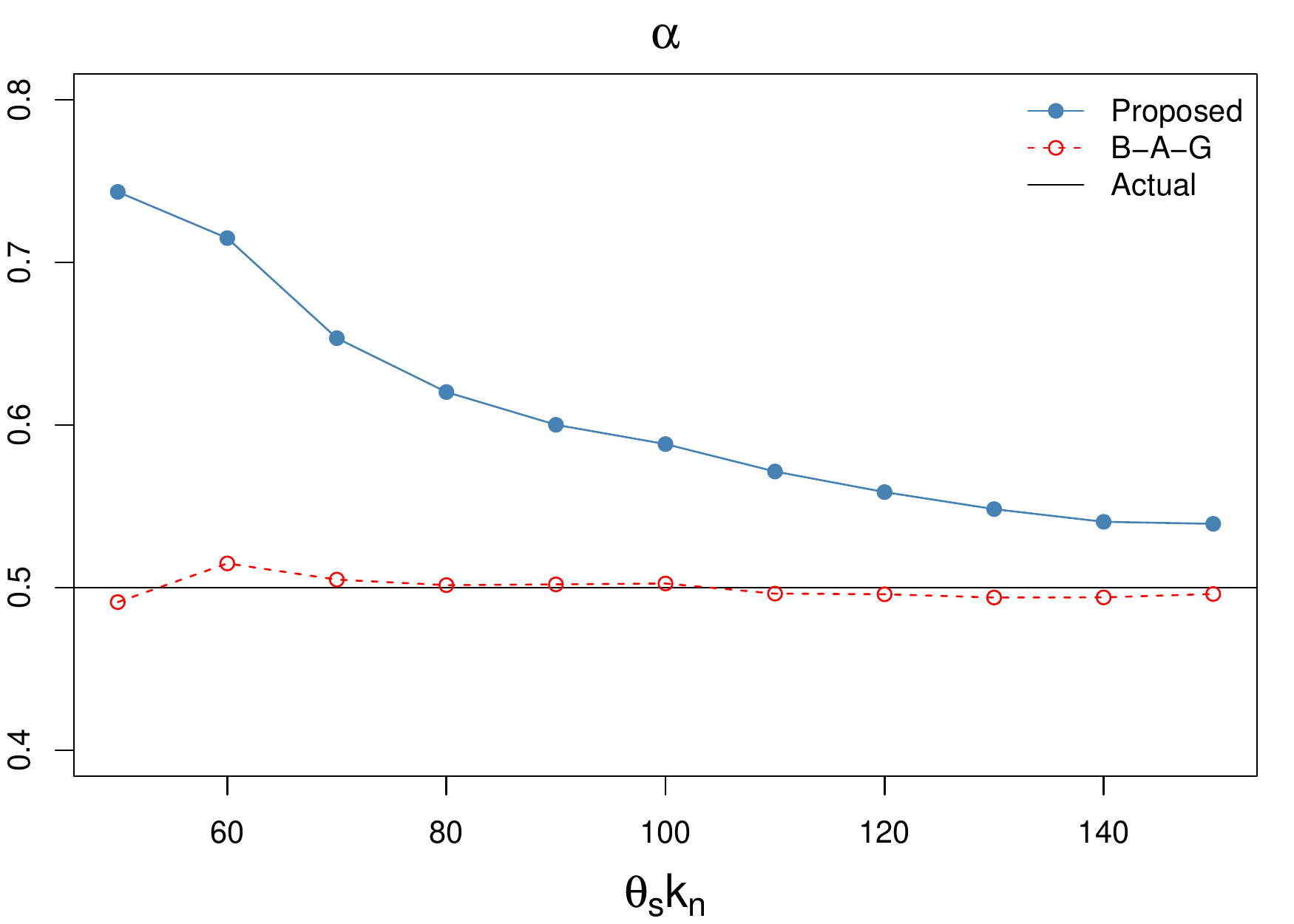}
  \end{minipage}
  \hfill
  \begin{minipage}[b]{0.45\textwidth}
  \centering
    \includegraphics[width=0.95\textwidth]{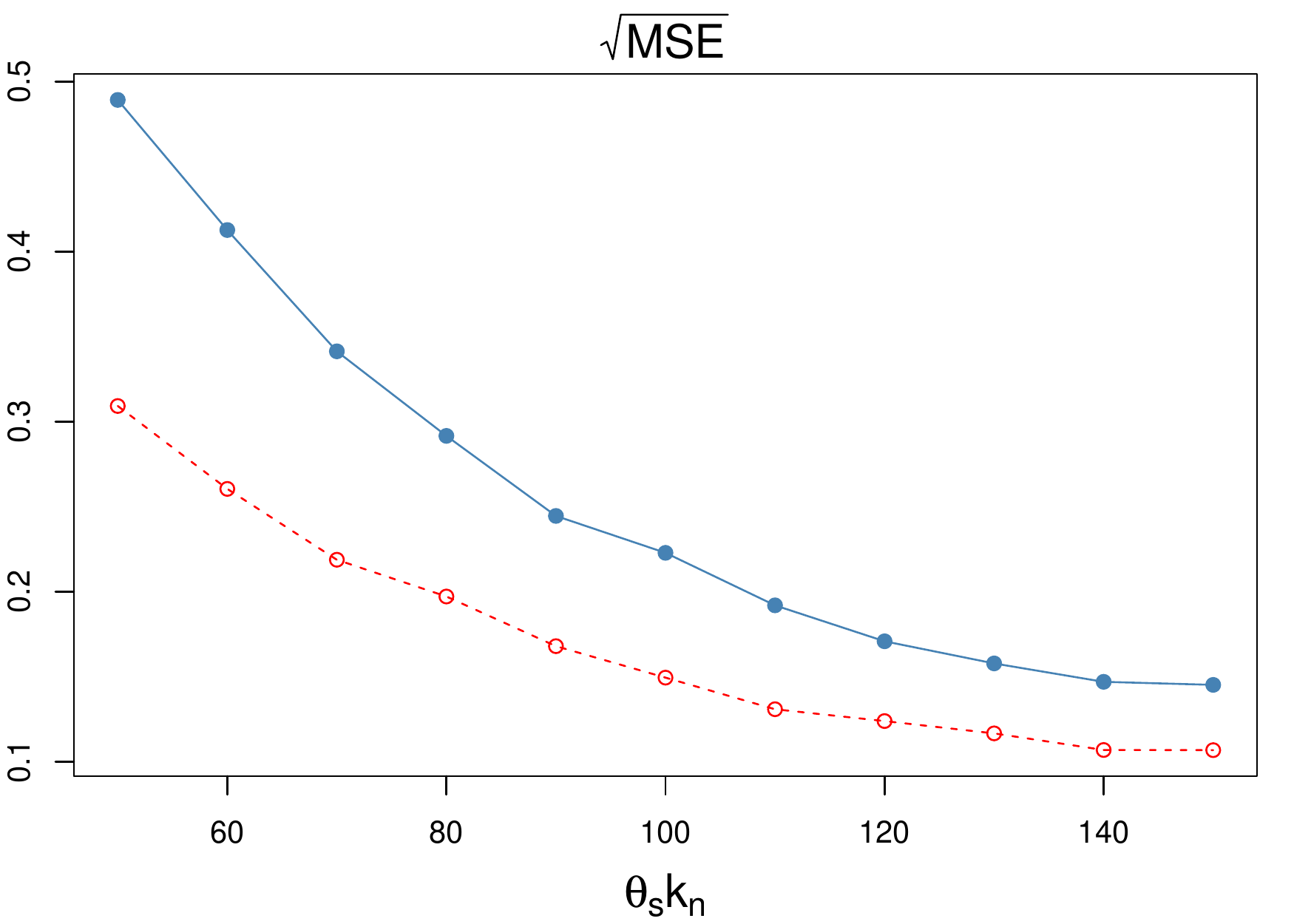}
  \end{minipage}
  \caption{Estimated tail index and $ \sqrt{\text{MSE}} $ for Pareto samples. $n = 500$, $ \alpha = 0.5 $, $ k_n = 50 $, $ \delta = 1 $} \label{figure: alpha pareto delta 1}
\end{figure}

Figures \ref{figure: tr pareto delta 1} and \ref{figure: alpha pareto delta 1} show the averaged estimates of $ \alpha $ and $ \delta k_n $ as well as the estimated mean squared errors (MSE) with different $ \theta_s k_n $. Estimates by the proposed method are plotted in solid lines while those by the method in \citet{Beirlant:2016gs} are in dashed lines.
The proposed method overestimates the tail index $ \alpha $, especially when the number of upper order statistics included in the estimation is small. This is not unexpected, as the method does not assume the data are from a Pareto distribution and thus does not benefit from the extra information that the bias term in the likelihood should be zero.
However, the proposed method estimates the number of missing extreme values accurately, and the estimation is robust to different numbers of upper order statistics included. 

\begin{figure}[H]
  \centering
  \begin{minipage}[b]{0.45\textwidth}
  \centering
    \includegraphics[width=0.95\textwidth]{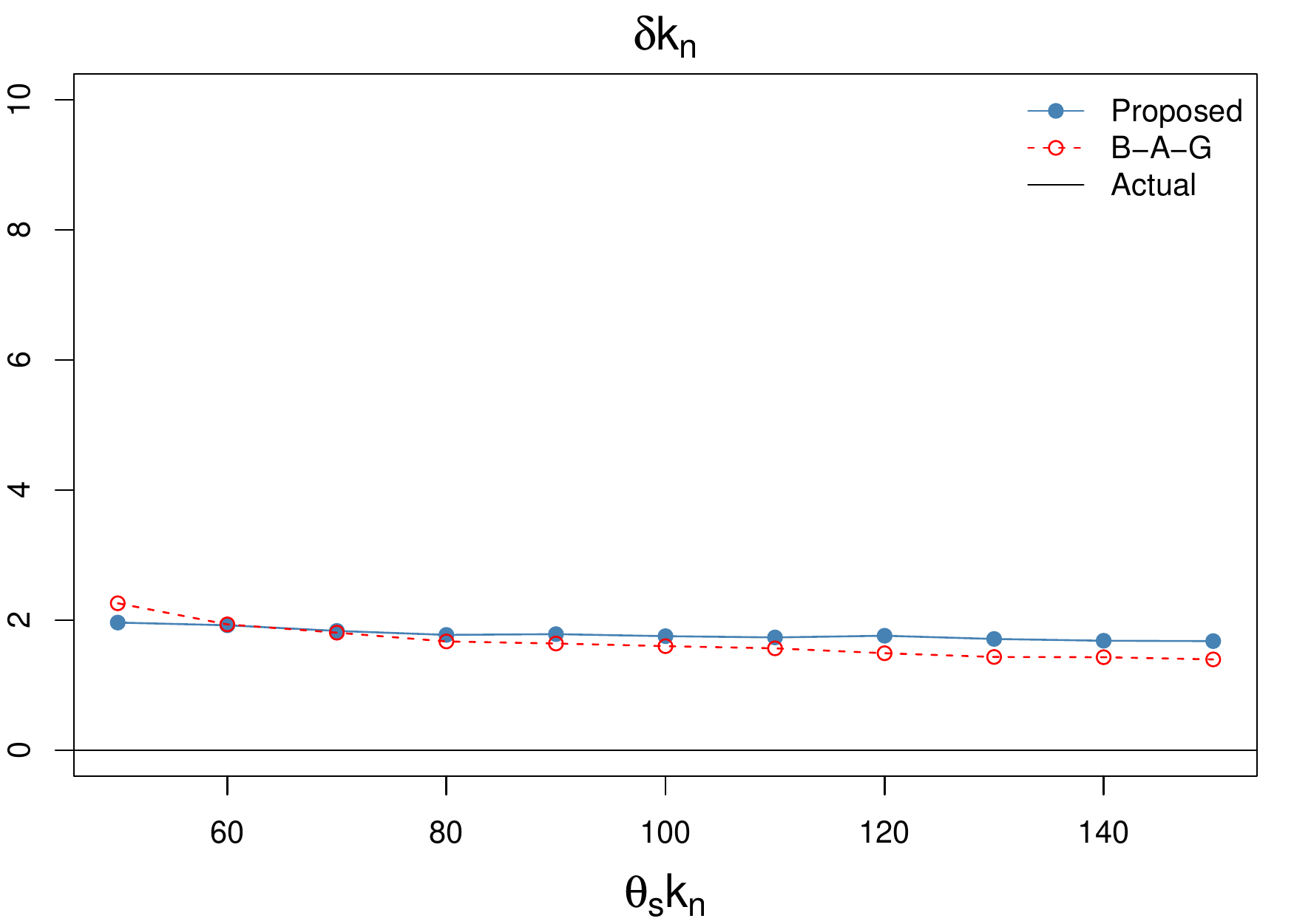}
  \end{minipage}
  \hfill
  \begin{minipage}[b]{0.45\textwidth}
  \centering
    \includegraphics[width=0.95\textwidth]{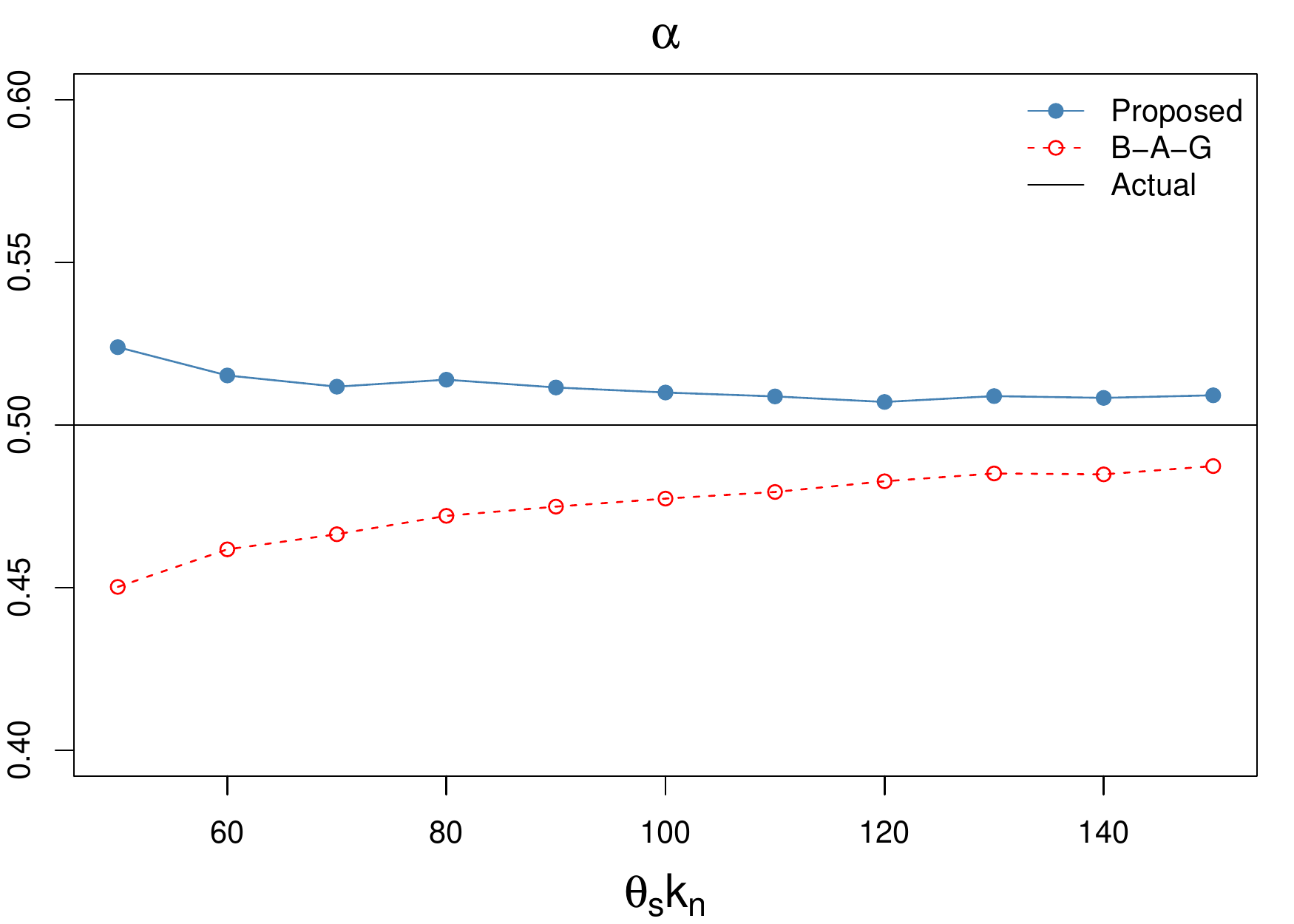}
  \end{minipage}
  \caption{Averages of estimated number of missing extremes and tail index for $ 200 $ Pareto samples. $n = 500$, $ \alpha = 0.5 $, $ k_n = 50 $, $ \delta = 0 $} \label{figure: tr pareto 0}
\end{figure}

We also examine the efficacy of the estimation procedure for $ 200 $ independent Pareto samples without any extreme values missing ($\delta=0$). 
Figure \ref{figure: tr pareto 0} shows that both methods give accurate estimates of the tail index and are able to estimate the number of missing extremes to be close to zero. 


\subsection{Non-Pareto Examples}

Next we examine the scenarios when the data are not from a Pareto distribution. Observations used here are generated from Cauchy and Student's $ t $-distributions. The following results show that the proposed method continues to perform well in estimating the number of missing extremes, even for distributions whose tail indices are more challenging to estimate when the top extremes are unobserved.

\subsubsection{Cauchy Distribution}

\begin{figure}[H]
  \centering
  \begin{minipage}[b]{0.45\textwidth}
  \centering
    \includegraphics[width=0.95\textwidth]{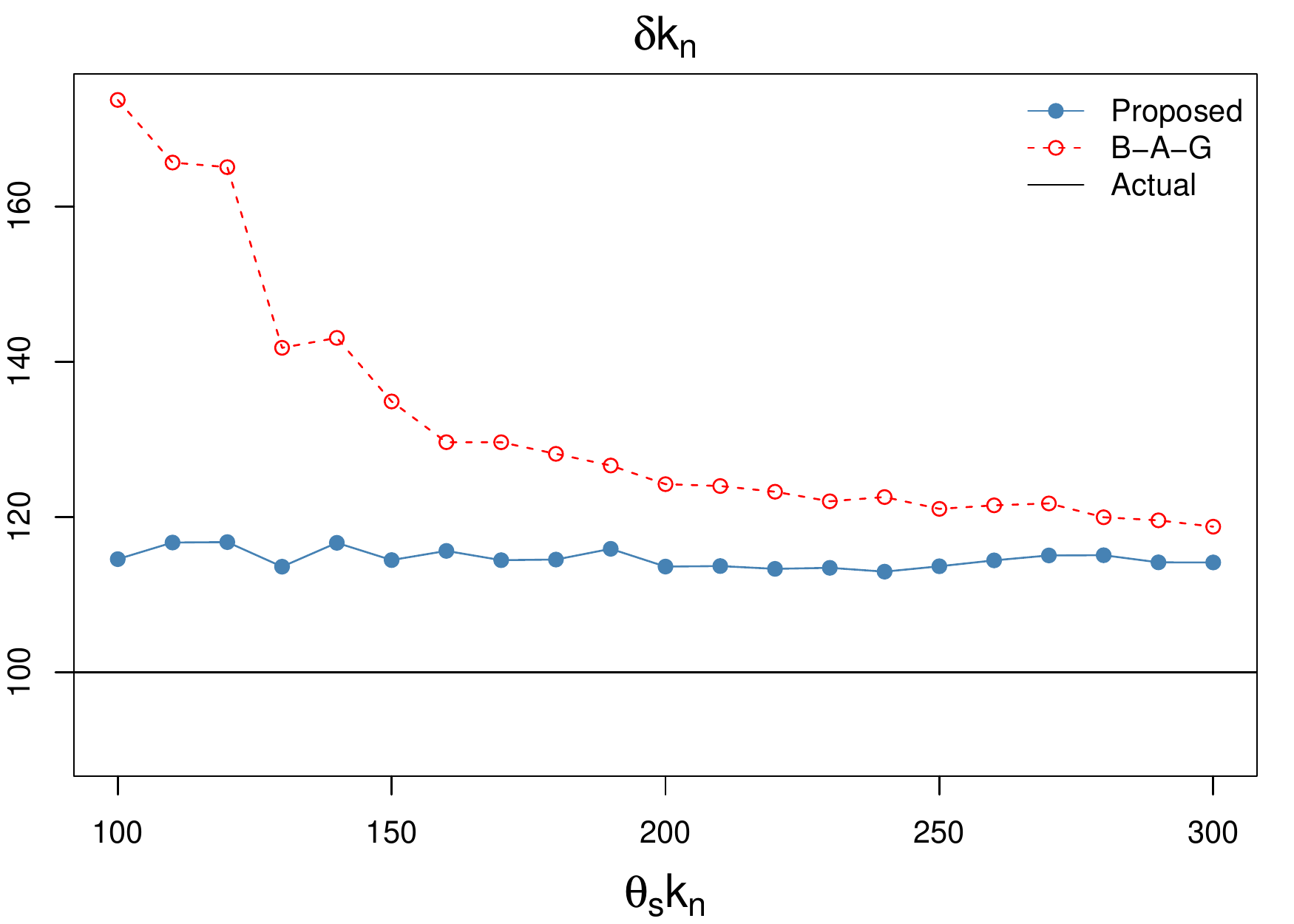}
  \end{minipage}
  \hfill
  \begin{minipage}[b]{0.45\textwidth}
  \centering
    \includegraphics[width=0.95\textwidth]{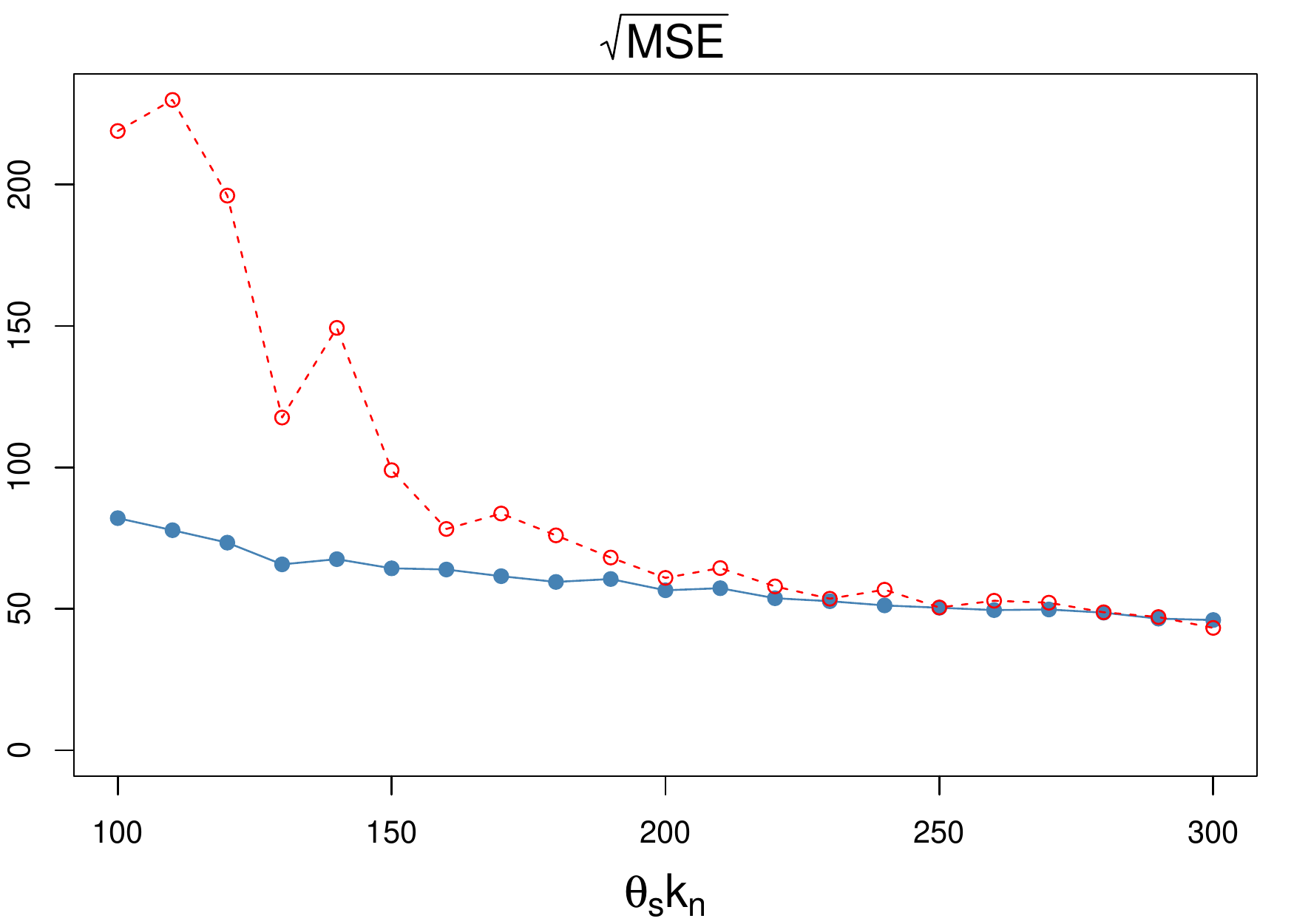}
  \end{minipage}
  \caption{Averages of estimated number of missing extremes and $ \sqrt{\text{MSE}} $ for $ 200 $ Cauchy samples. $n = 2000$, $ \alpha = 1 $, $ k_n = 100 $, $ \delta = 1 $} \label{figure: tr cauchy 1}
\end{figure}

\begin{figure}[H]
  \centering
  \begin{minipage}[b]{0.45\textwidth}
  \centering
    \includegraphics[width=0.95\textwidth]{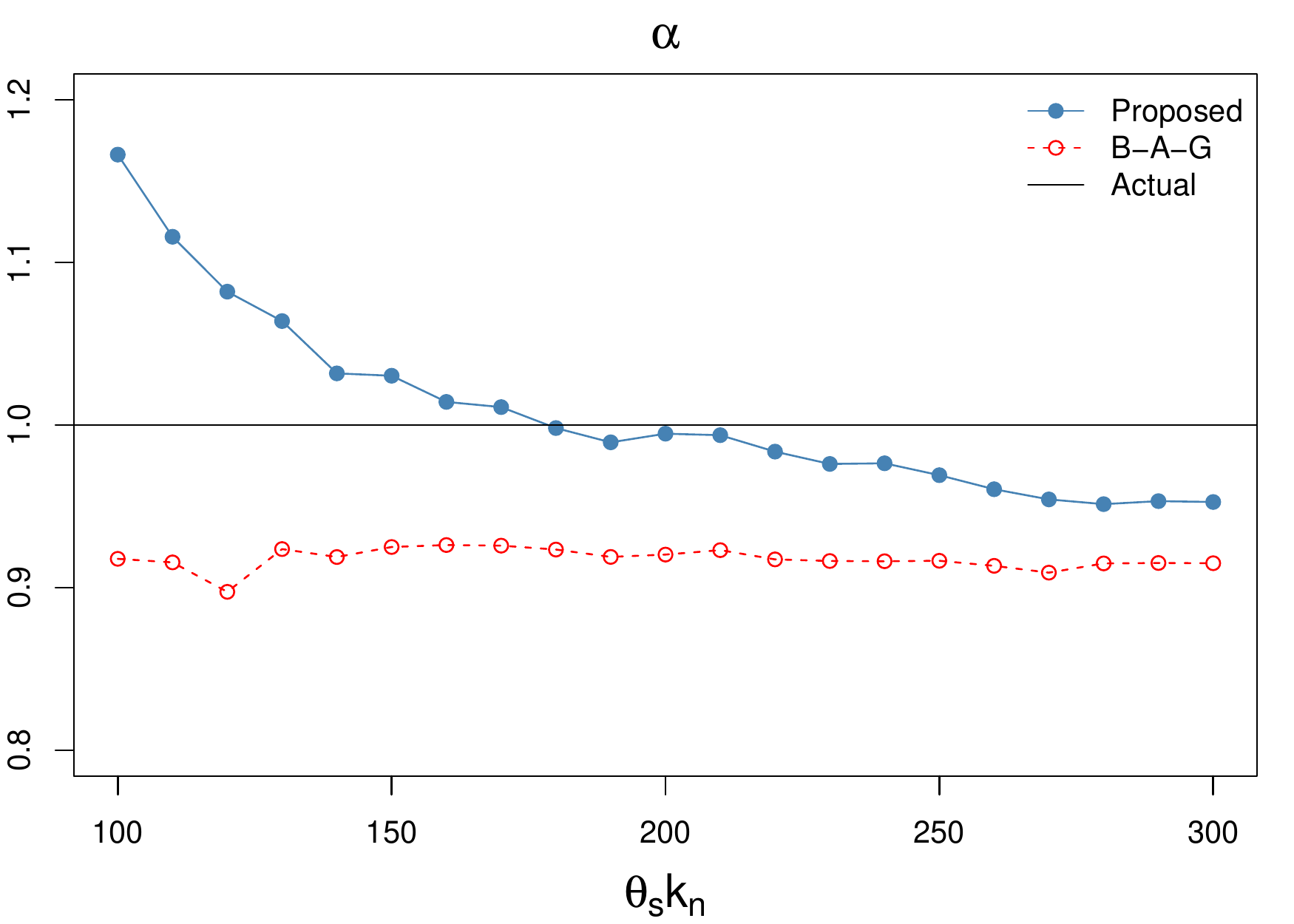}
  \end{minipage}
  \hfill
  \begin{minipage}[b]{0.45\textwidth}
  \centering
    \includegraphics[width=0.95\textwidth]{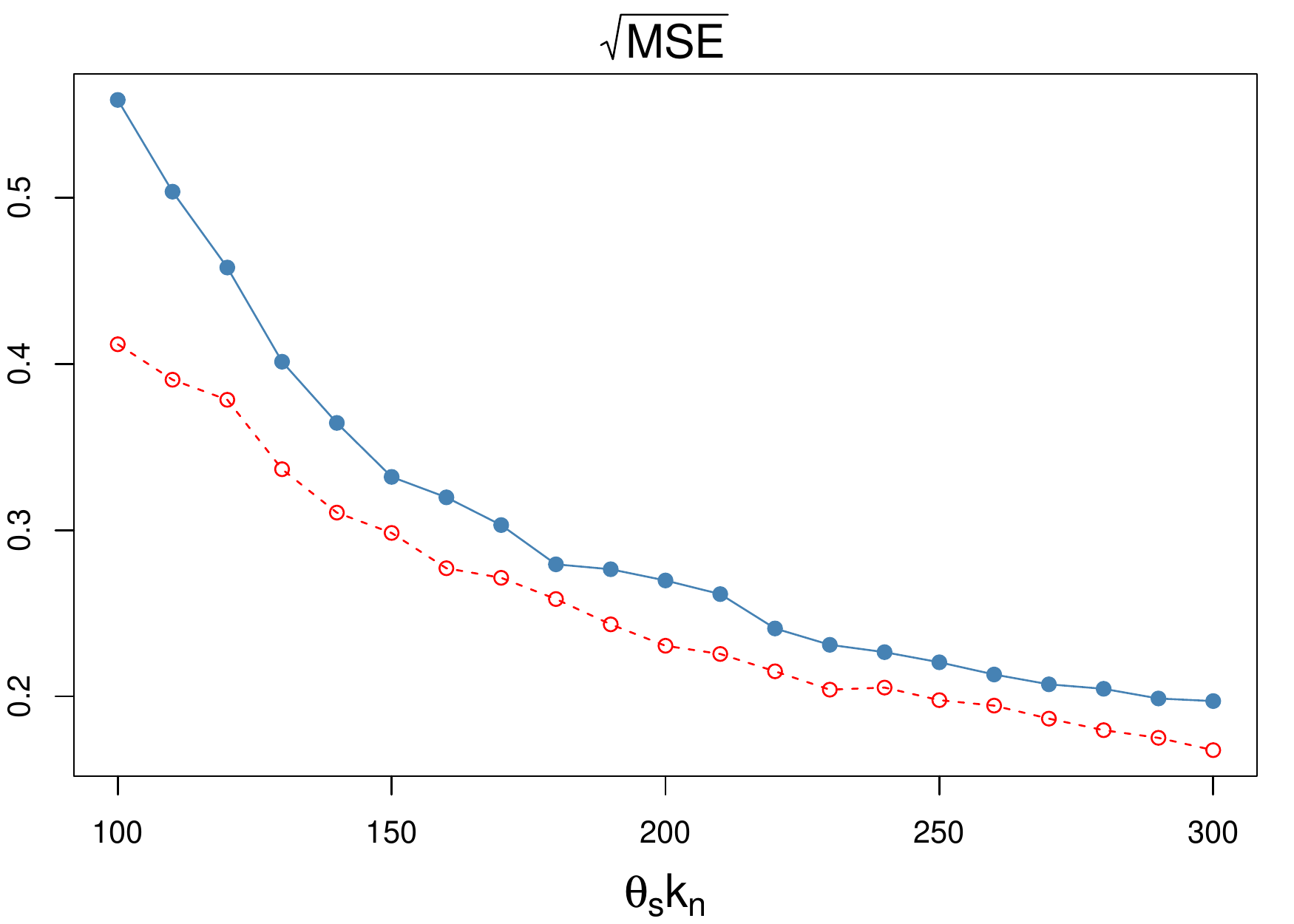}
  \end{minipage}
  \caption{Averages of estimated tail index and $ \sqrt{\text{MSE}} $ for $ 200 $ Cauchy samples. $n = 2000$, $ \alpha = 1 $, $ k_n = 100 $, $ \delta = 1 $} \label{figure: alpha cauchy 1}
\end{figure}

Figures \ref{figure: tr cauchy 1} and \ref{figure: alpha cauchy 1} show averaged estimates for $ 200 $ independent Cauchy samples with the largest $ 100 $ observations removed from each sample.

\begin{figure}[H]
  \centering
  \begin{minipage}[b]{0.45\textwidth}
  \centering
    \includegraphics[width=0.92\textwidth]{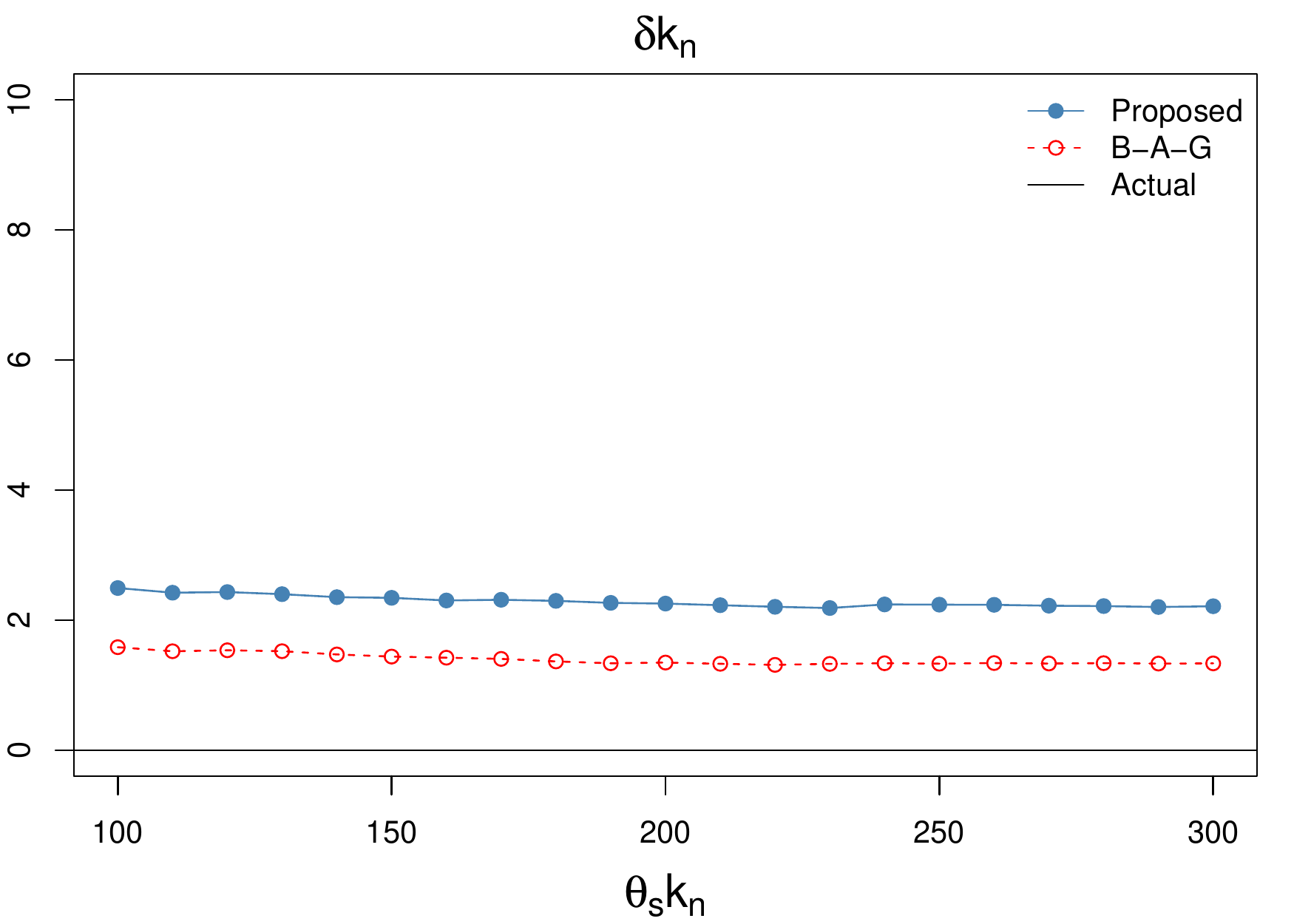}
  \end{minipage}
  \hfill
  \begin{minipage}[b]{0.45\textwidth}
  \centering
    \includegraphics[width=0.92\textwidth]{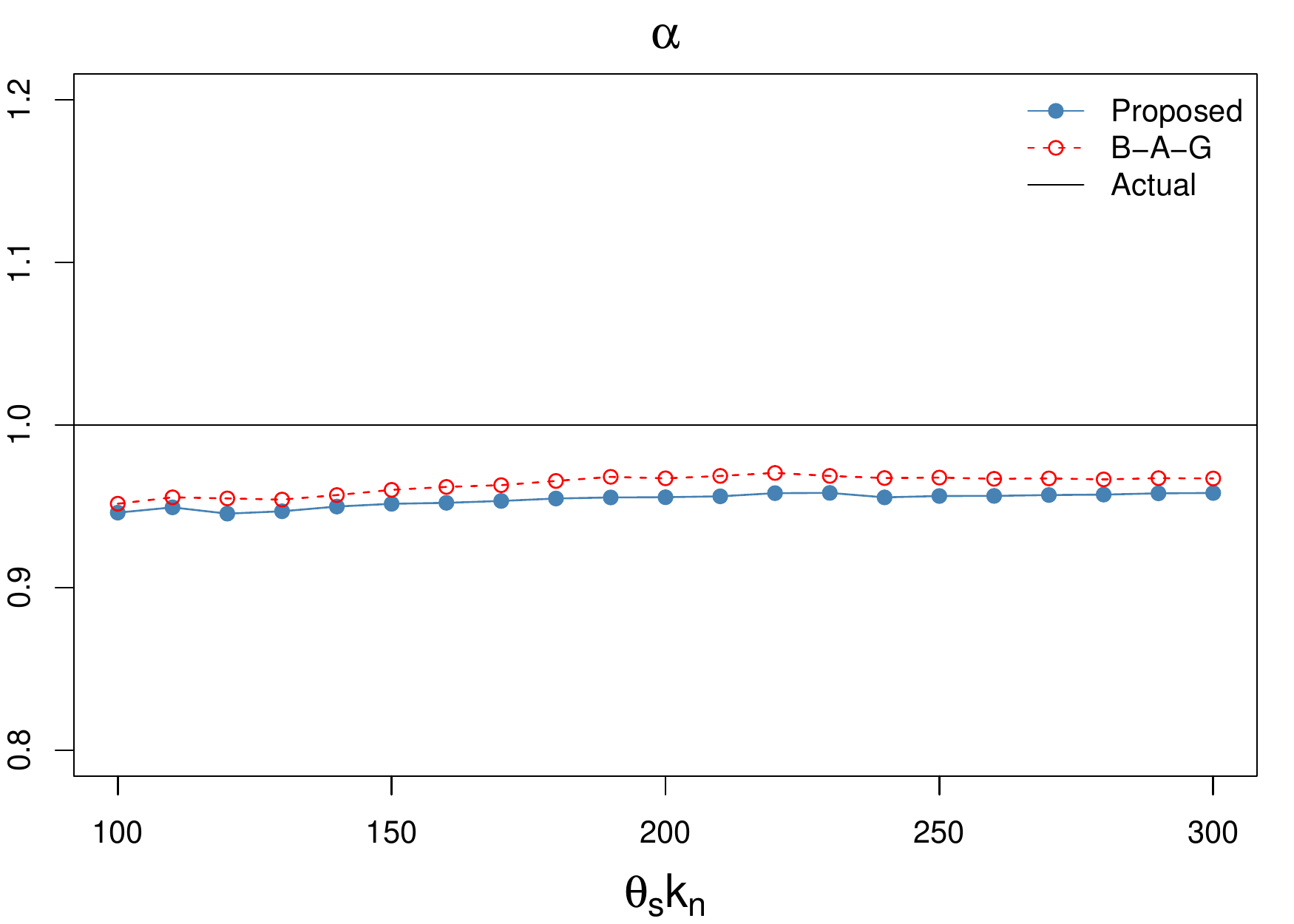}
  \end{minipage}
  \caption{Averages of estimated number of missing extremes and tail index for $ 200 $ Cauchy samples. $n = 2000$, $ \alpha = 1 $, $ k_n = 100 $, $ \delta = 0 $} \label{figure: cauchy 0}
\end{figure}

Figure \ref{figure: cauchy 0} shows the estimates for $ 200 $ independent Cauchy samples without any extremes missing. Both methods produce accurate results for the zero number of missing extremes and the tail index. 

\subsubsection{Student's $ t_{2.5} $  Distribution}
\begin{figure}[H]
  \centering
  \begin{minipage}[b]{0.45\textwidth}
  \centering
    \includegraphics[width=0.92\textwidth]{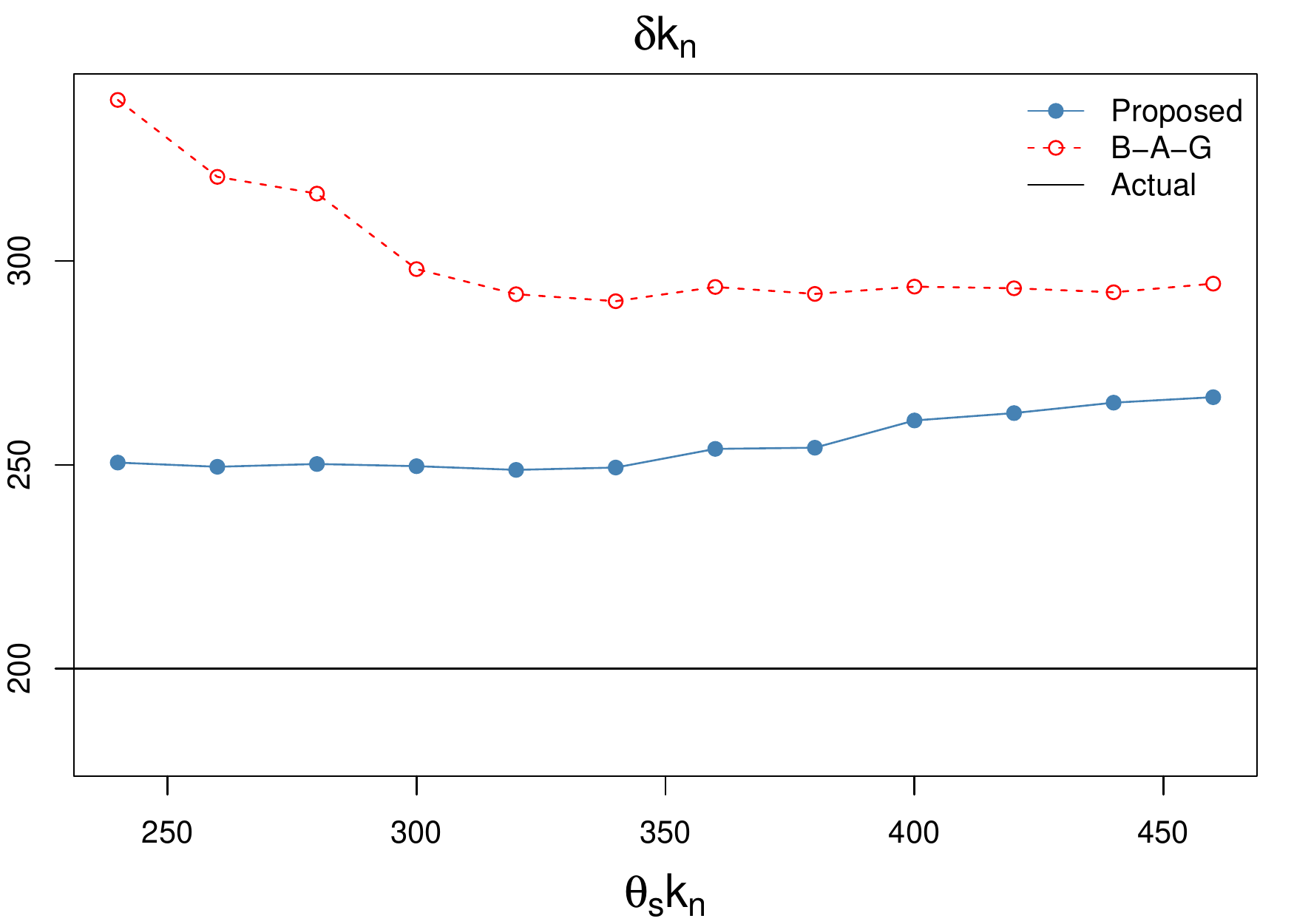}
  \end{minipage}
  \hfill
  \begin{minipage}[b]{0.45\textwidth}
  \centering
    \includegraphics[width=0.92\textwidth]{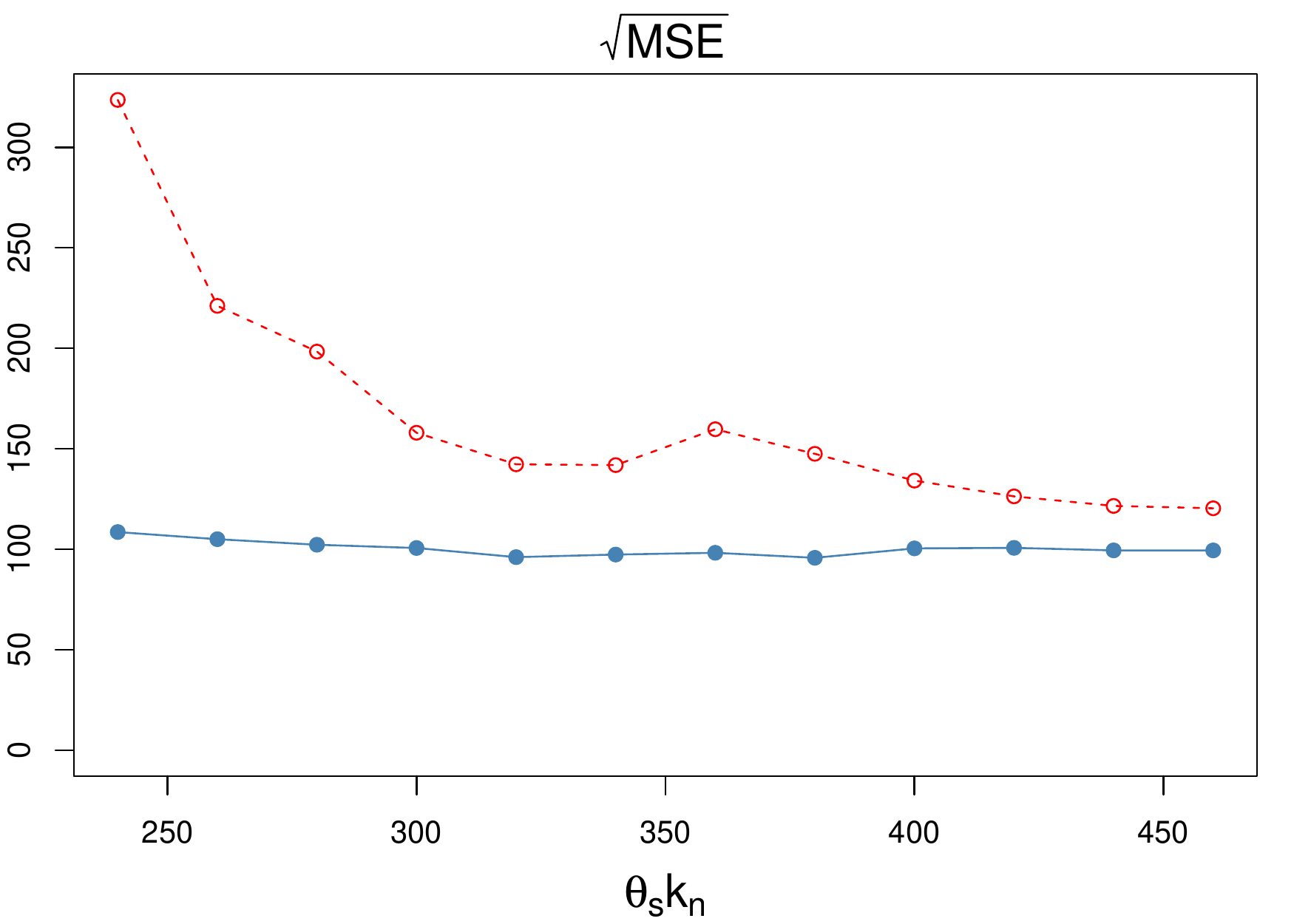}
  \end{minipage}
  \caption{Averages of estimated number of missing extremes and $ \sqrt{\text{MSE}} $ for $ 200 $ Student's $ t_{2.5} $ samples. $n = 10,000$, $ \alpha = 2.5 $, $ k_n = 200 $, $ \delta = 1 $} \label{figure: tr t 1}
\end{figure}

\begin{figure}[H]
  \centering
  \begin{minipage}[b]{0.45\textwidth}
  \centering
    \includegraphics[width=0.95\textwidth]{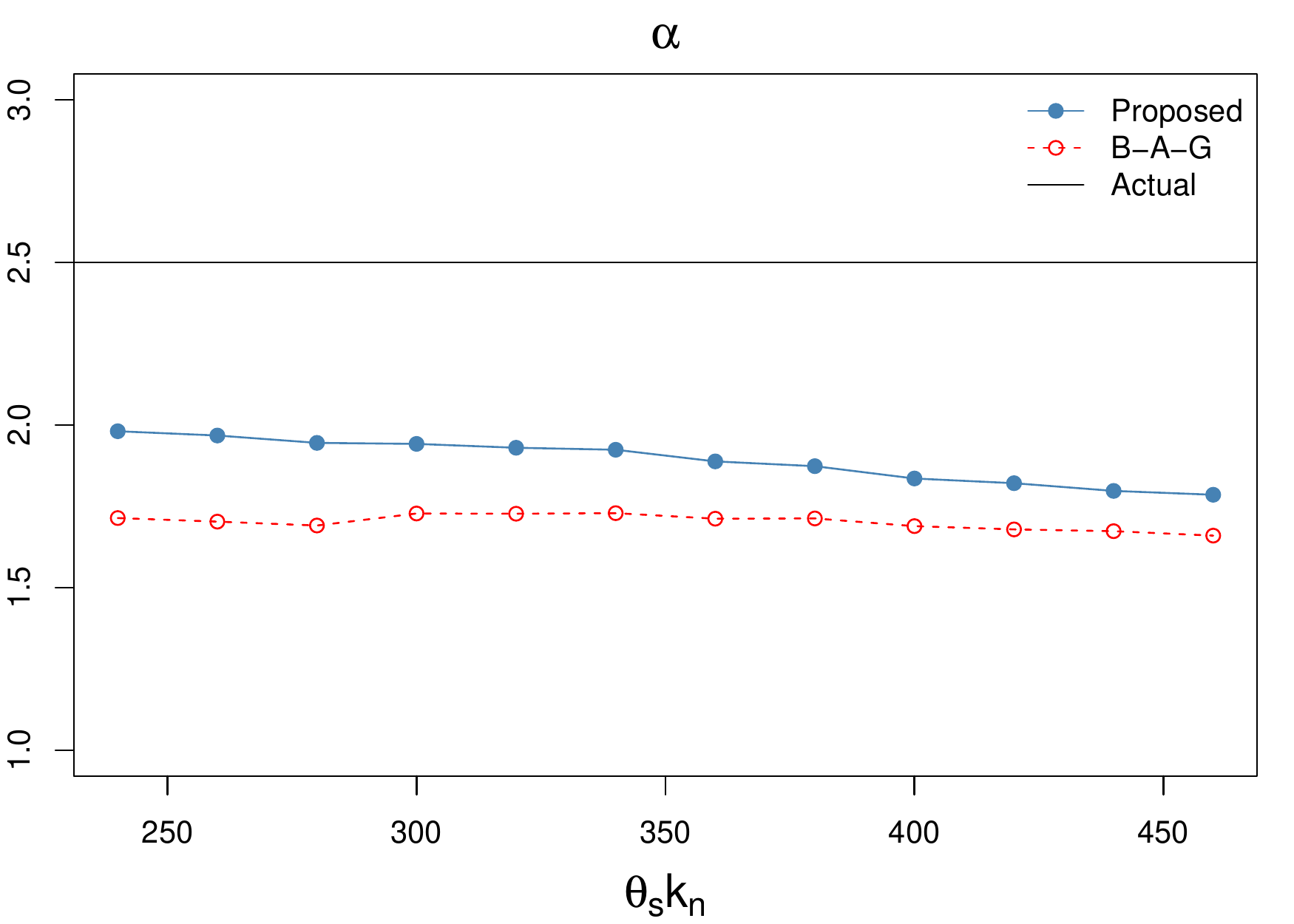}
  \end{minipage}
  \hfill
  \begin{minipage}[b]{0.45\textwidth}
  \centering
    \includegraphics[width=0.95\textwidth]{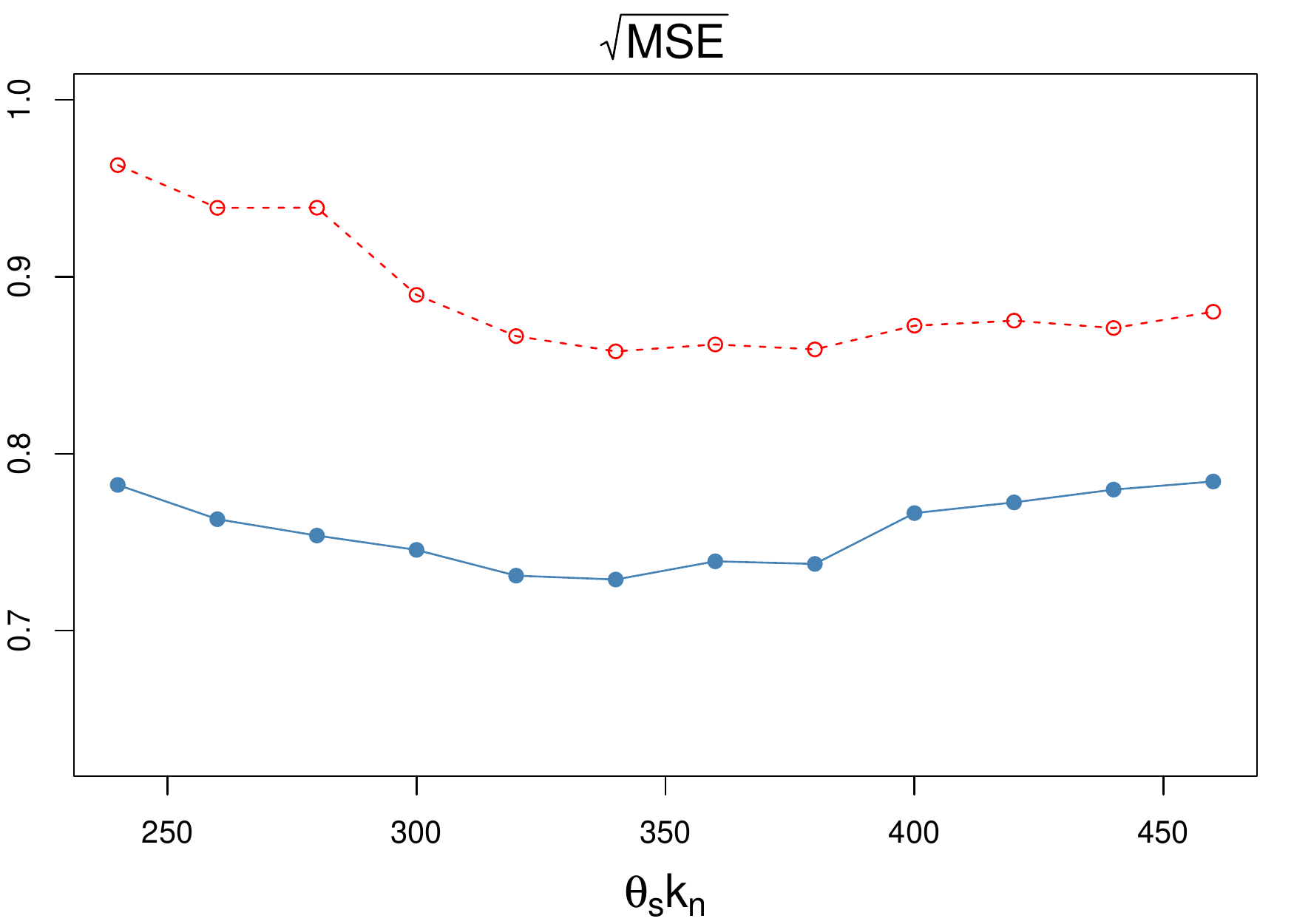}
  \end{minipage}
  \caption{Averages of estimated tail index and $ \sqrt{\text{MSE}} $ for $ 200 $ Student's $ t_{2.5} $ samples. $n = 10000$, $ \alpha = 2.5 $, $ k_n = 200 $, $ \delta = 1 $} \label{figure: alpha t 1}
\end{figure}

Figures \ref{figure: tr t 1} and \ref{figure: alpha t 1} show the estimates for $ 200 $ independent samples from the Student's $ t $-distribution with degrees of freedom $ df = 2.5 $. The tail index $ \alpha = df $. In each sample there are $ n = 10000 $ observations originally. Let $ k_n = 200 $ and $ \delta = 1 $ so that the largest $ 200 $ observations have been removed from each of the original samples.

\subsection{Robustness to Model Parameters}

To examine the robustness of the proposed method to different model parameters, we applied the proposed estimation procedure to data generated from Pareto and Cauchy distributions with different parameter values and compared the accuracy of the estimation for these different settings. 

\begin{figure}[H]
  \centering
  \begin{minipage}[b]{0.45\textwidth}
  \centering
    \includegraphics[width=0.95\textwidth]{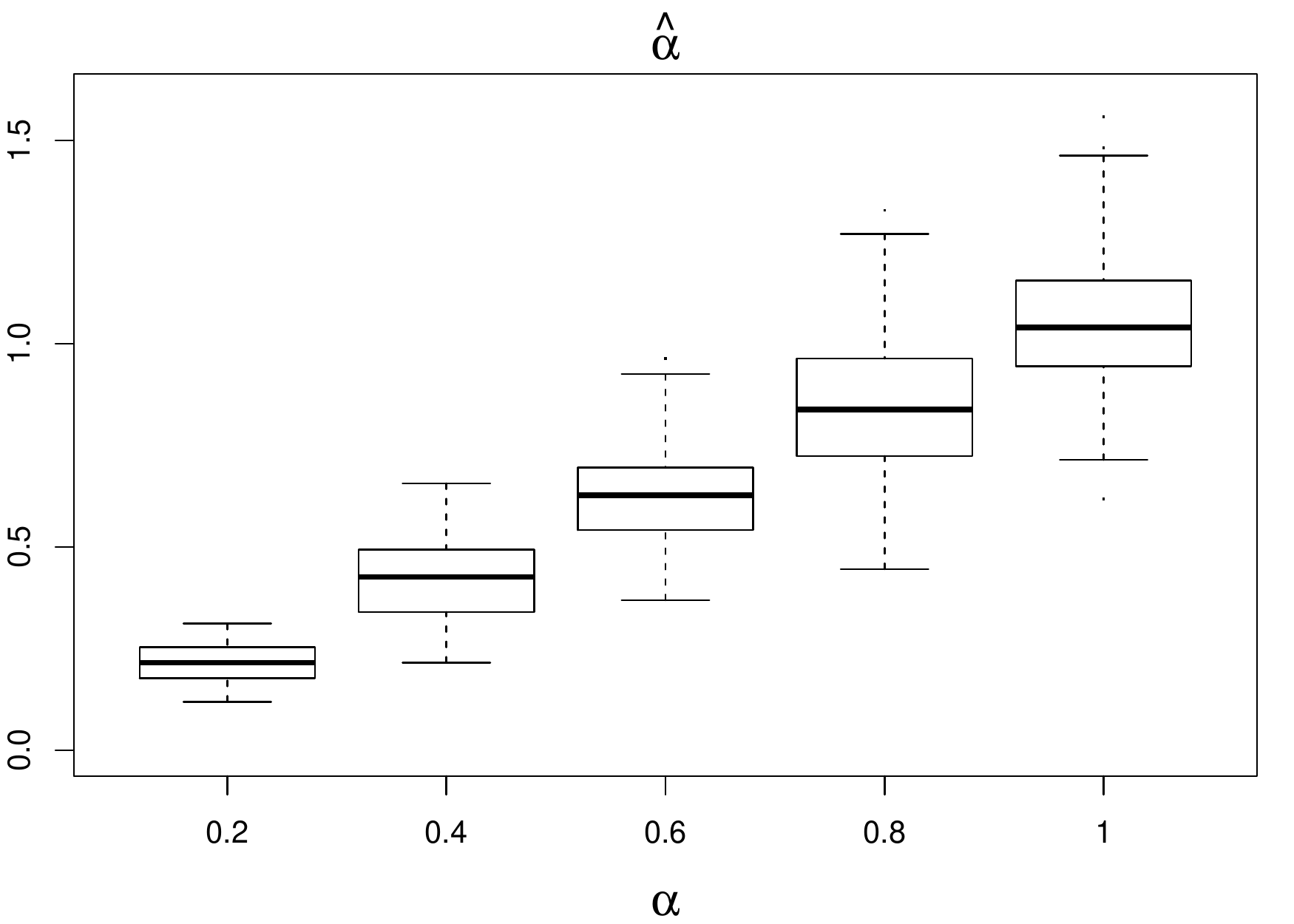}
  \end{minipage}
  \hfill
  \begin{minipage}[b]{0.45\textwidth}
  \centering
    \includegraphics[width=0.95\textwidth]{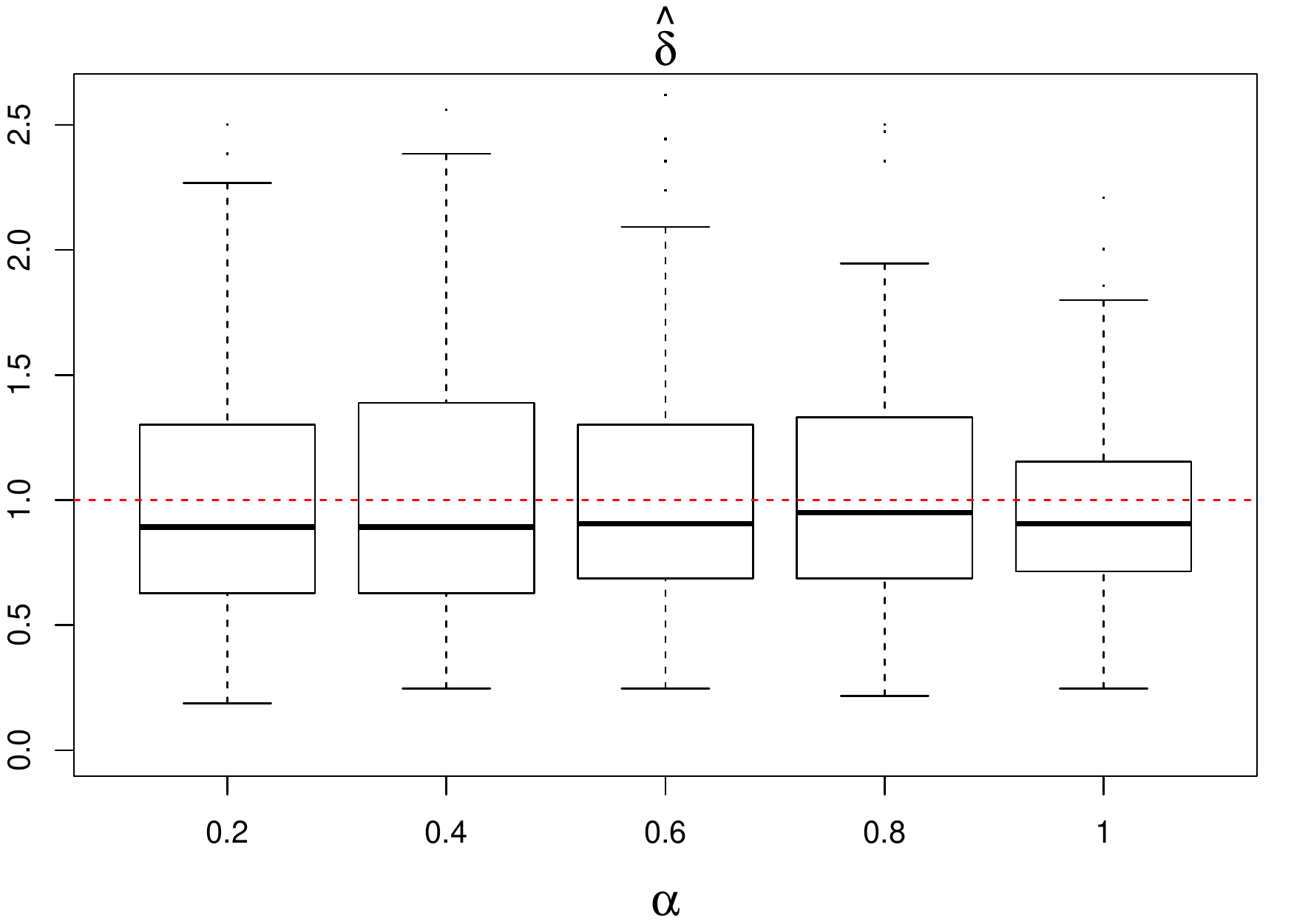}
  \end{minipage}
  \caption{Boxplots of  $\hat{\alpha}$ (left) and $\hat{\delta}$ (right) estimated using the proposed method with $100$ independent samples of size $n = 500$ from the Pareto distribution with $\alpha = 0.2, 0.4, 0.6, 0.8, 1$, respectively. $k_n = 50$ and $\delta = 1$ for all Pareto samples in removing the top extremes. } \label{figure: pareto_change_alpha}
\end{figure}

Figure \ref{figure: pareto_change_alpha} shows the estimation of the tail index $\alpha$ and the parameter $\delta$ for removing top extremes with data generated from the Pareto distribution with different values of $\alpha$ for $\delta = 1$ and $n = 500$ fixed. Each boxplot summarizes the major quantiles (1\%, 25\%, 50\%, 75\%, 99\%) of the estimation results from $100$ independent samples under the designated parameter setting. The $x$-axis indicates the values of $\alpha$ in the model from which the data are generated. In obtaining the estimation results, we use a fixed range $\theta_s k_n = 180$ as the top extremes included in the estimation procedure. 
Results did not appreciably change for different choices of $\theta_s k_n$ that are in a reasonable range. 
In practice, it is suggested to perform the estimation procedure over a series of values of $\theta_s k_n$ to determine a reasonable value for the estimation.

\begin{figure}[H]
  \centering
  \begin{minipage}[b]{0.45\textwidth}
  \centering
    \includegraphics[width=0.95\textwidth]{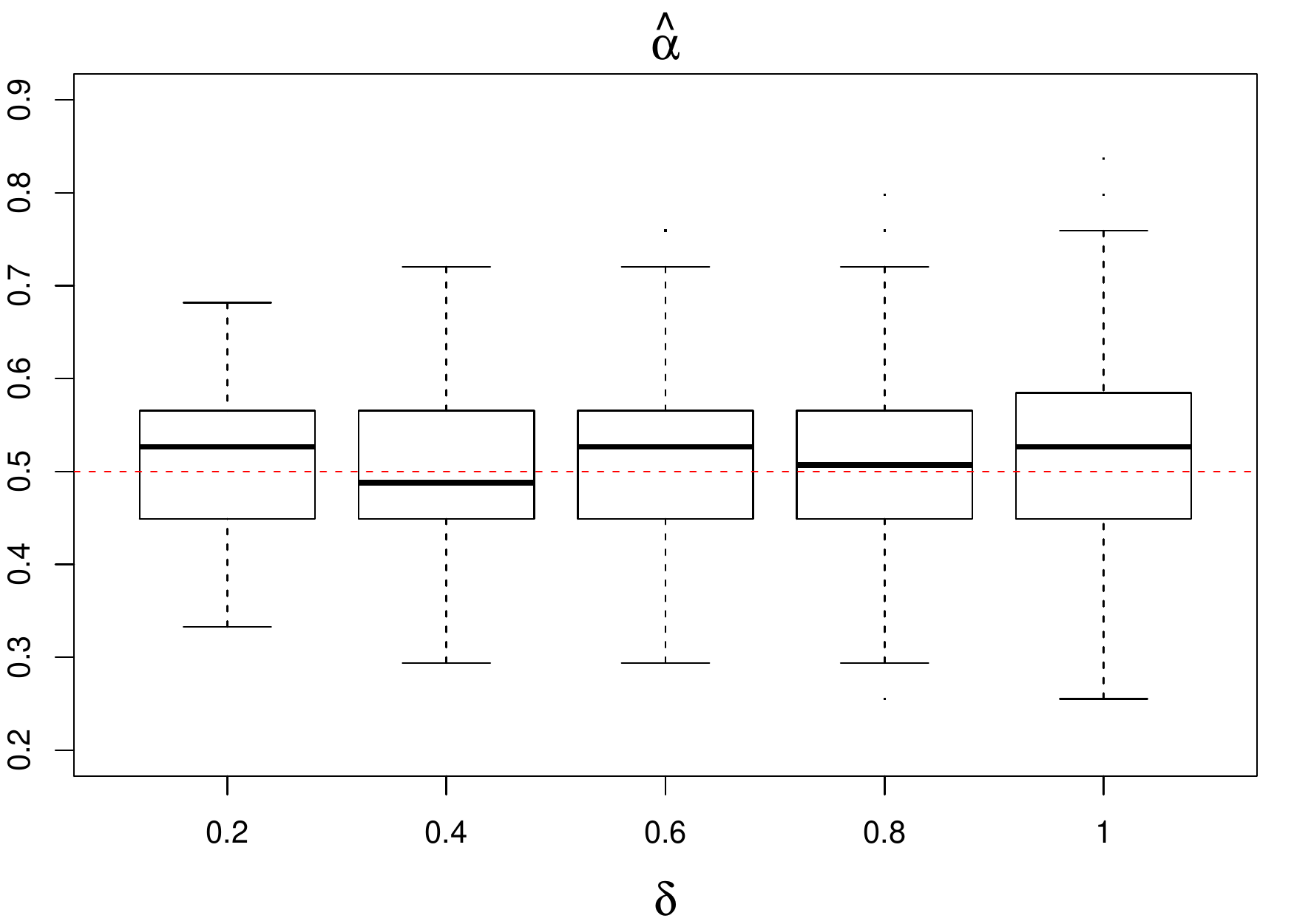}
  \end{minipage}
  \hfill
  \begin{minipage}[b]{0.45\textwidth}
  \centering
    \includegraphics[width=0.95\textwidth]{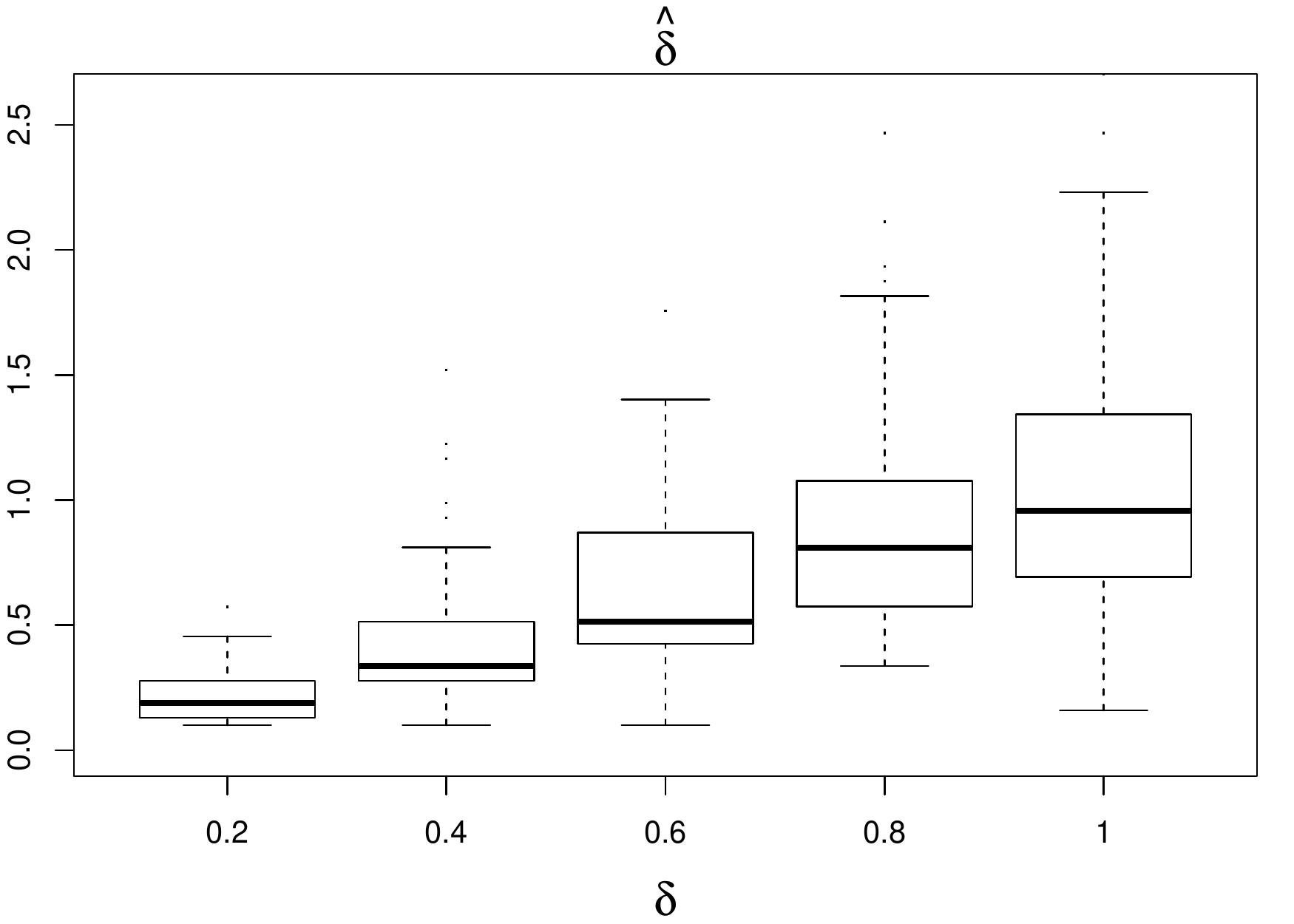}
  \end{minipage}
  \caption{Boxplots of  $\hat{\alpha}$ (left) and $\hat{\delta}$ (right) estimated using the proposed method with $100$ independent samples of size $n = 500$ from the Pareto distribution with $k_n = 50$, $\alpha = 0.5$ and $\delta = 0.2, 0.4, 0.6, 0.8, 1$, respectively, in removing the top extremes. } \label{figure: pareto_change_delta}
\end{figure}

Similarly, Figure \ref{figure: pareto_change_delta} shows the estimation of the tail index $\alpha$ and the parameter $\delta$ with data generated from the Pareto distribution with $n = 500$, $\alpha = 0.5$ fixed and different values of $\delta$. The $x$-axis of each plot indicates the values of $\delta$ in removing top extremes from the Pareto samples.

\begin{figure}[H]
  \centering
  \begin{minipage}[b]{0.45\textwidth}
  \centering
    \includegraphics[width=0.95\textwidth]{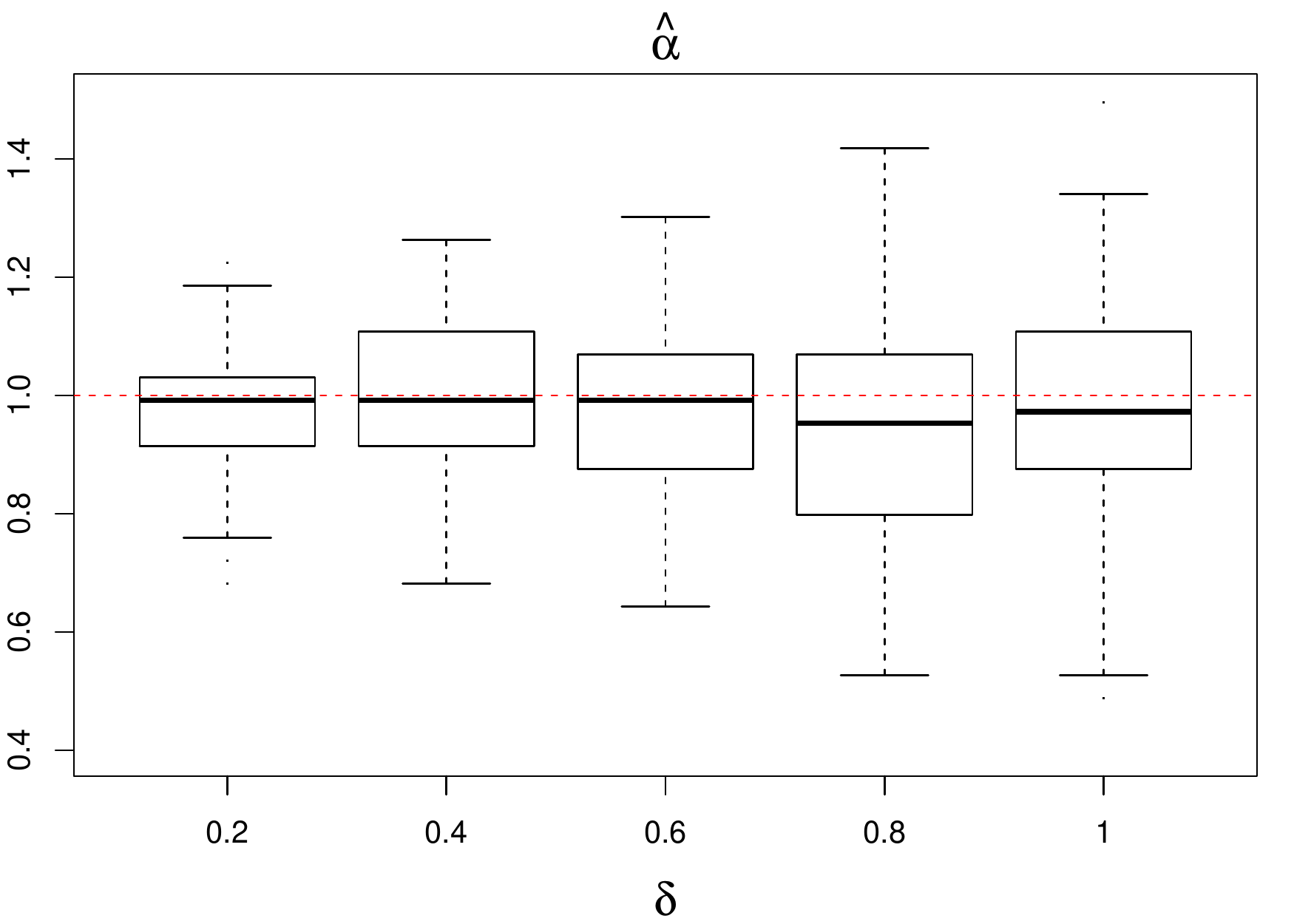}
  \end{minipage}
  \hfill
  \begin{minipage}[b]{0.45\textwidth}
  \centering
    \includegraphics[width=0.95\textwidth]{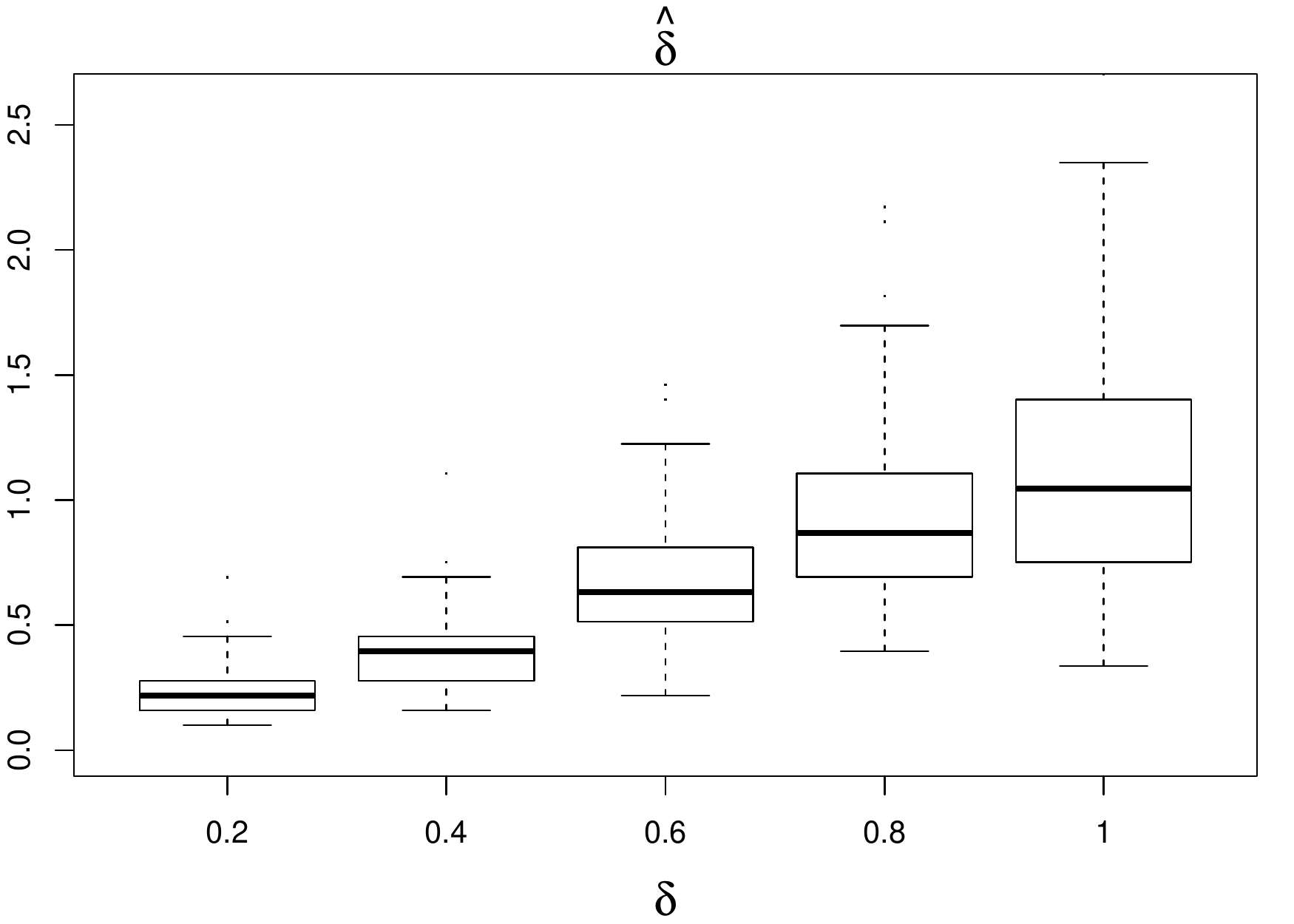}
  \end{minipage}
  \caption{Boxplots of  $\hat{\alpha}$ (left) and $\hat{\delta}$ (right) estimated using the proposed method with $100$ independent samples of size $n = 2000$ from the Cauchy distribution with $k_n = 100$, $\alpha = 1$ and $\delta = 0.2, 0.4, 0.6, 0.8, 1$, respectively, in removing the top extremes. } \label{figure: cauchy_change_delta}
\end{figure}

As an example of non-Pareto distributions, Figure \ref{figure: cauchy_change_delta} summarizes the estimation of the tail index $\alpha$ and the parameter $\delta$ for removing top extremes with data generated from the Cauchy distribution with $n = 2000$, $\alpha = 1$ and different values of $\delta$. The selected range of top extremes to include in the estimation is $\theta_s k_n = 320$.

In both Pareto and non-Pareto cases the proposed estimation procedure produces results that are reasonably robust to changing model parameters. 

\begin{figure}[H]
  \centering
  \begin{minipage}[b]{0.45\textwidth}
  \centering
    \includegraphics[width=0.95\textwidth]{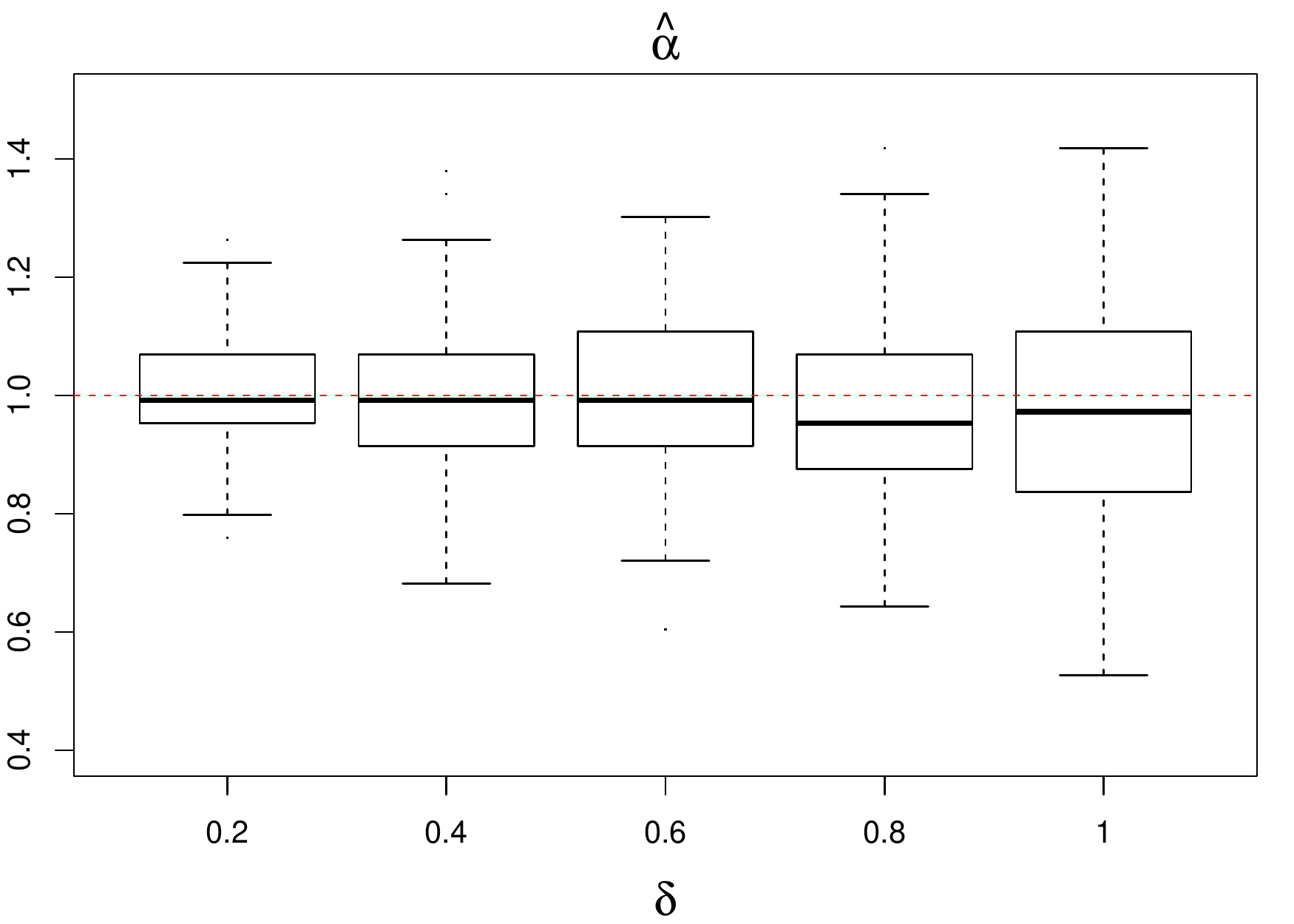}
  \end{minipage}
  \hfill
  \begin{minipage}[b]{0.45\textwidth}
  \centering
    \includegraphics[width=0.95\textwidth]{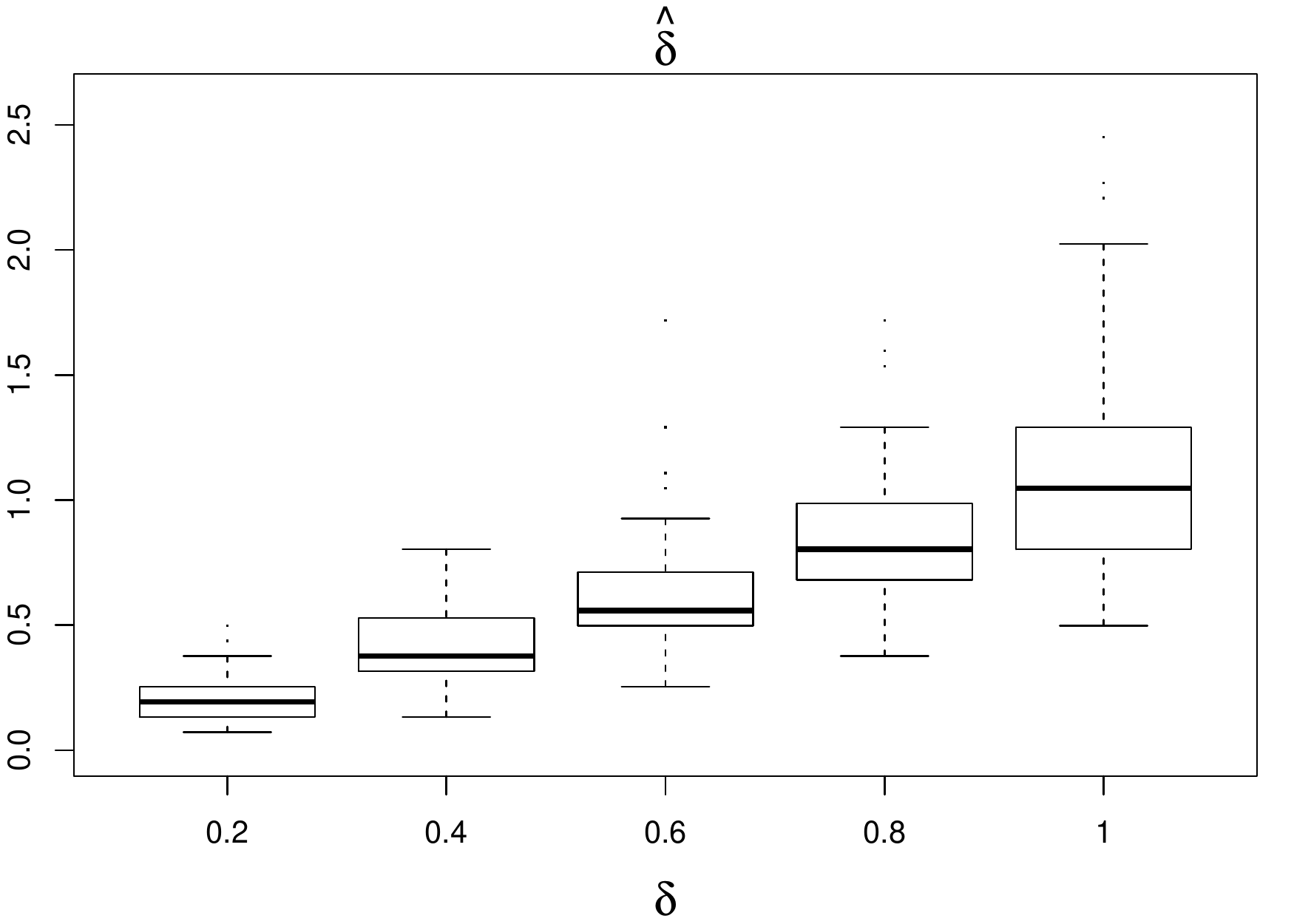}
  \end{minipage}
  \caption{Boxplots of  $\hat{\alpha}$ (left) and $\hat{\delta}$ (right) estimated using the proposed method with $100$ independent samples of size $n = 4000$ from the Cauchy distribution with $k_n = 200$, $\alpha = 1$ and $\delta = 0.2, 0.4, 0.6, 0.8, 1$, respectively, in removing the top extremes. } \label{figure: cauchy_change_delta_2}
\end{figure}

In addition, Figure \ref{figure: cauchy_change_delta_2} shows the estimation results with independent samples of size $n = 4000$ from the Cauchy distribution with different values of $\delta$. The selected range of top extremes for the estimation $\theta_s k_n = 440$.
The comparison of Figure \ref{figure: cauchy_change_delta} and \ref{figure: cauchy_change_delta_2} indicates the proposed method produces  robust results despite different sample sizes.

\subsection{Light-tailed Example}

Simulations in the above sections focused on heavy-tailed samples. One might ask if the Hill curve of a light-tailed sample would exhibit similar patterns as the Hill plot of a heavy-tailed sample with missing extremes and whether the proposed method is capable of identifying the different cases. 

\begin{figure}[H]
  \centering
  \begin{minipage}[b]{0.45\textwidth}
  \centering
    \includegraphics[width=0.92\textwidth]{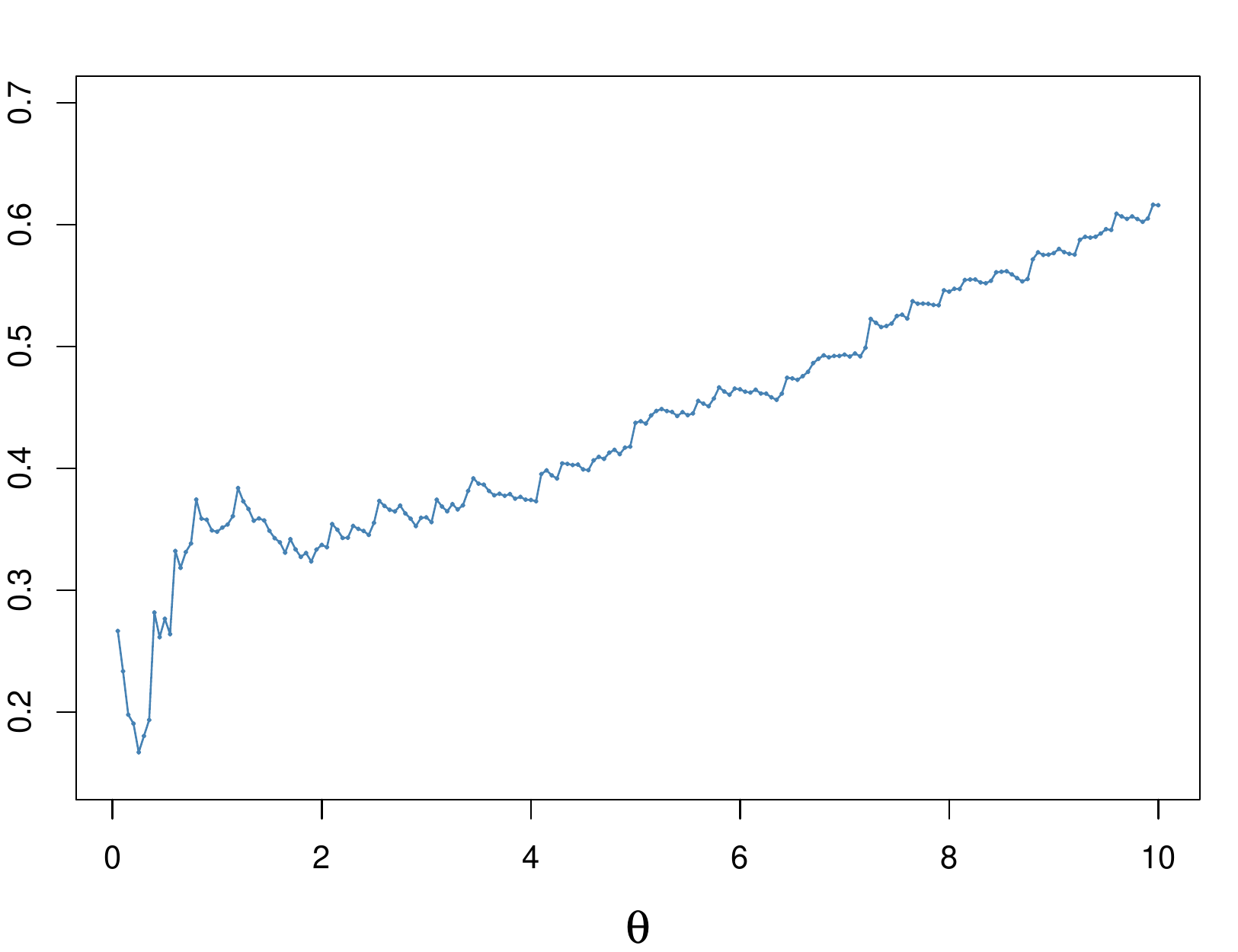}
  \end{minipage}
  \hfill
  \begin{minipage}[b]{0.45\textwidth}
  \centering
    \includegraphics[width=0.92\textwidth]{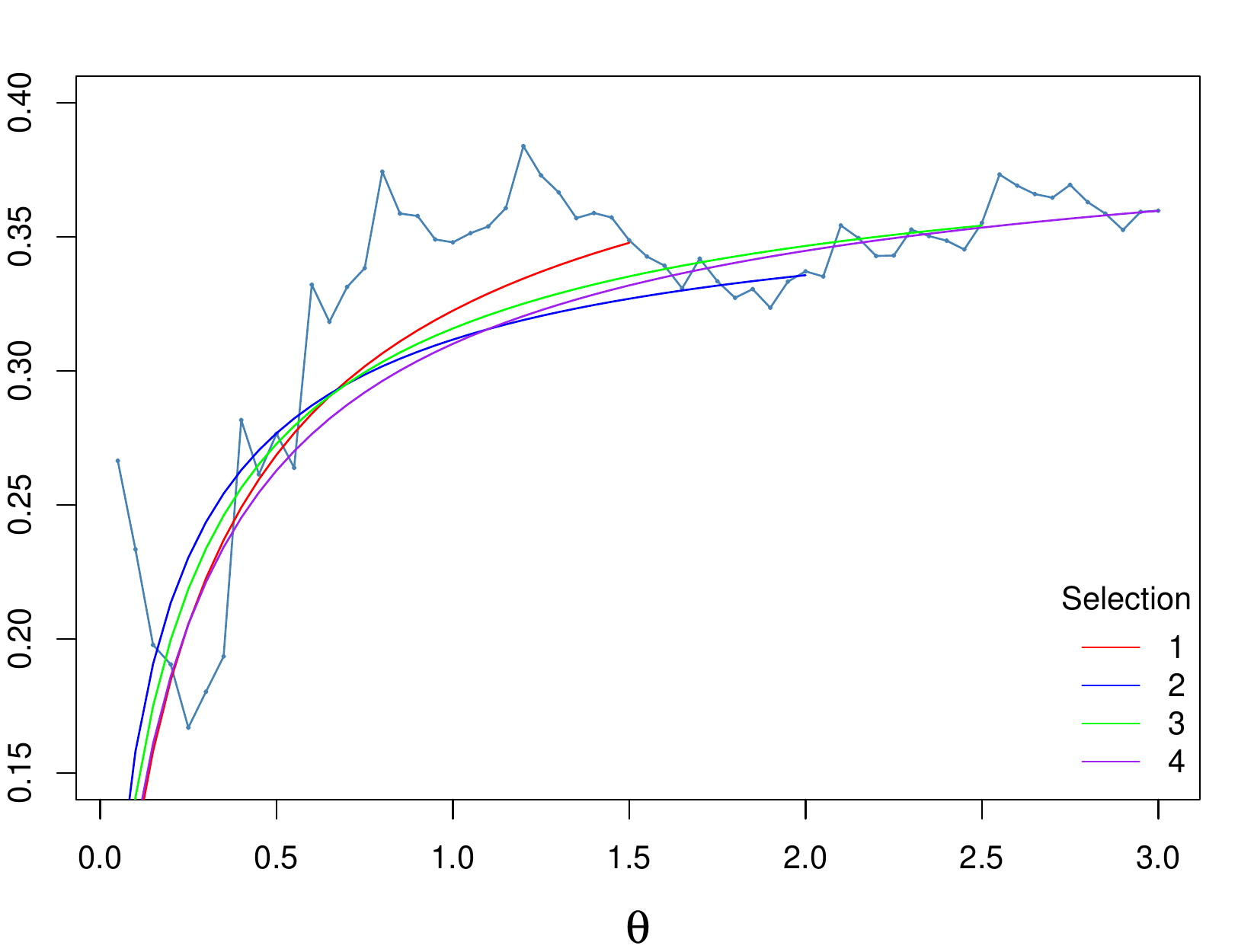}
  \end{minipage}
  \caption{Left: Hill plot for a standard exponential sample of size $n = 500$, $ \gamma = 0 $. Right: Hill plot overlaid with estimated mean curves of Gaussian processes using $ k_n = 20 $. Selection 1: $  \hat{\gamma}= 0.33 $, $ \hat{\delta} k_n = 2.26 $. 2: $ \hat{\gamma} = 0.30 $, $ \hat{\delta} k_n = 1.14 $. 3: $ \hat{\gamma} = 0.31 $, $ \hat{\delta} k_n = 1.58 $. 4: $ \hat{\gamma} =  0.33$, $ \hat{\delta} k_n = 1.92 $} \label{figure: exp 0 hill}
\end{figure}

Here we demonstrate that the proposed method can indeed differentiate between the light- and heavy-tailed cases with an example of light-tailed data without any missing values. The left panel of Figure \ref{figure: exp 0 hill} is the Hill plot based on a sample of $ 500 $ from the standard exponential distribution. Although the curve is generally increasing, it is not as smooth as in the case of heavy-tailed data with missing extremes. In the right panel of Figure \ref{figure: exp 0 hill}, the Hill plot is overlaid with mean curves of Gaussian processes estimated using different parts of the observed Hill curve based on the method in Section \ref*{section estimation}. The estimates of missing extremes range from $ 0 $ to $ 3.6 $, which reflect the truth that there are no extreme values missing from the data. The true value of $ \gamma = 0 $ and the proposed method is also able to estimate $ \gamma $ with relatively small values. 

In summary, we have applied our estimation procedure to both Pareto and non-Pareto heavy tailed distributions. We have considered
both the standard scenario when all the extremes are present, and the scenario when some of the extremes are missing. Our method
is competitive in all cases, and it appears to work better in the non-Pareto cases due to its self-adjusting mechanism of reducing the bias. 
The simulation results show that our method is able to simultaneously estimate the tail index and the number of the missing extremes 
with a reasonable accuracy.

\section{Applications} \label{sect: data}

We now apply the proposed method to real data. In practice, the number of missing extreme values and the reason for their absence are usually unknown. The consistency of an estimation procedure can be tested by artificially removing a number of additional extremes from the observed data. Consistency requires that, in a certain range, such additional removal should not have a major effect
on the estimated tail index. Further, the estimated number of the originally missing upper order statistics should stay, approximately,
the same after accounting for the artificially removed observations. Here we examine a massive Google+ social network dataset and a moderate-sized earthquake fatality dataset, and in both cases the proposed procedure provides reasonable results.

\subsection{Google+}
We first apply our method to the data from the Google+ social network introduced in Section \ref{sec: intro}. The data contain one of the largest weakly connected components of a snapshot of the network taken on October 19, 2012. A weakly connected component of the network is created by treating the network as undirected and finding all nodes that can be reached from a randomly selected initial node.
There are 76,438,791 nodes and 1,442,504,499 edges in this component. 
The quantities of interest are the in- and out-degrees of nodes in the network, which often exhibit heavy-tailed properties (see, for example, Chapter 8 of \citet{Newman:2010ur}). 

We use, for estimation purposes, the largest $ 5000 $ values of the in-degree observations as the data set. We choose $k_n=200$. Next, we repeat the estimation procedure after artificially removing $ 400 $ largest of the $ 5000 $ values. 
In the estimation, we start from $ \theta_1 = 1 / k_n $ and let $ \theta_i = \theta_{i-1} + 1/k_n $ for $ 1 < i \le s $. 
As in the simulation studies, we consider a sequence of different endpoints $ \theta_s k_n $ and obtain estimates corresponding to different values of $ \theta_s k_n $.
For comparison, we also apply the estimation procedure of \citet{Beirlant:2016gs} to the dataset.

\begin{figure}[H]
\centering
\begin{minipage}{.45\textwidth}
  \centering
  \includegraphics[width=0.95\textwidth]{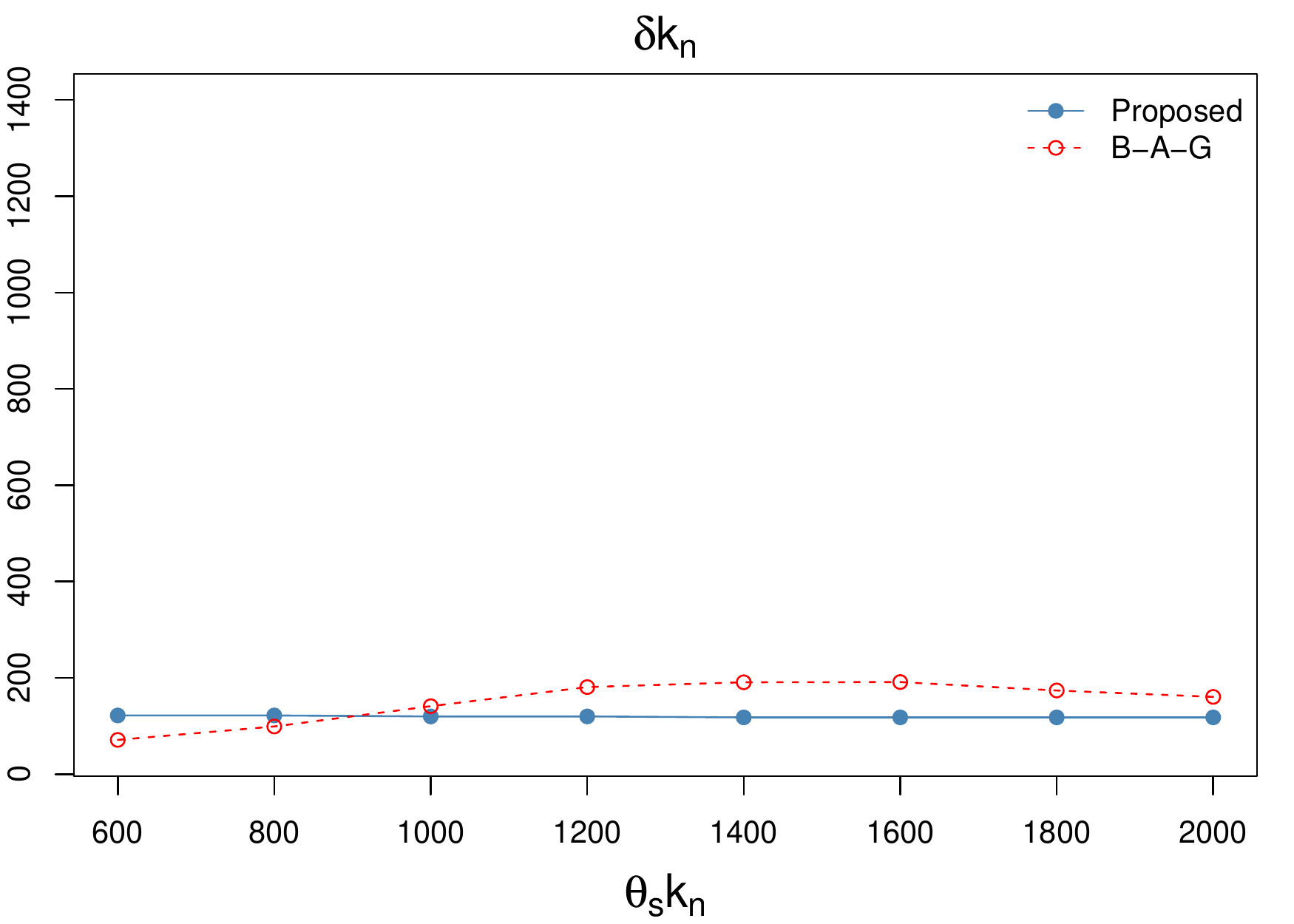}
\end{minipage}%
\hfill
\begin{minipage}{.45\textwidth}
  \centering
  \includegraphics[width=0.95\textwidth]{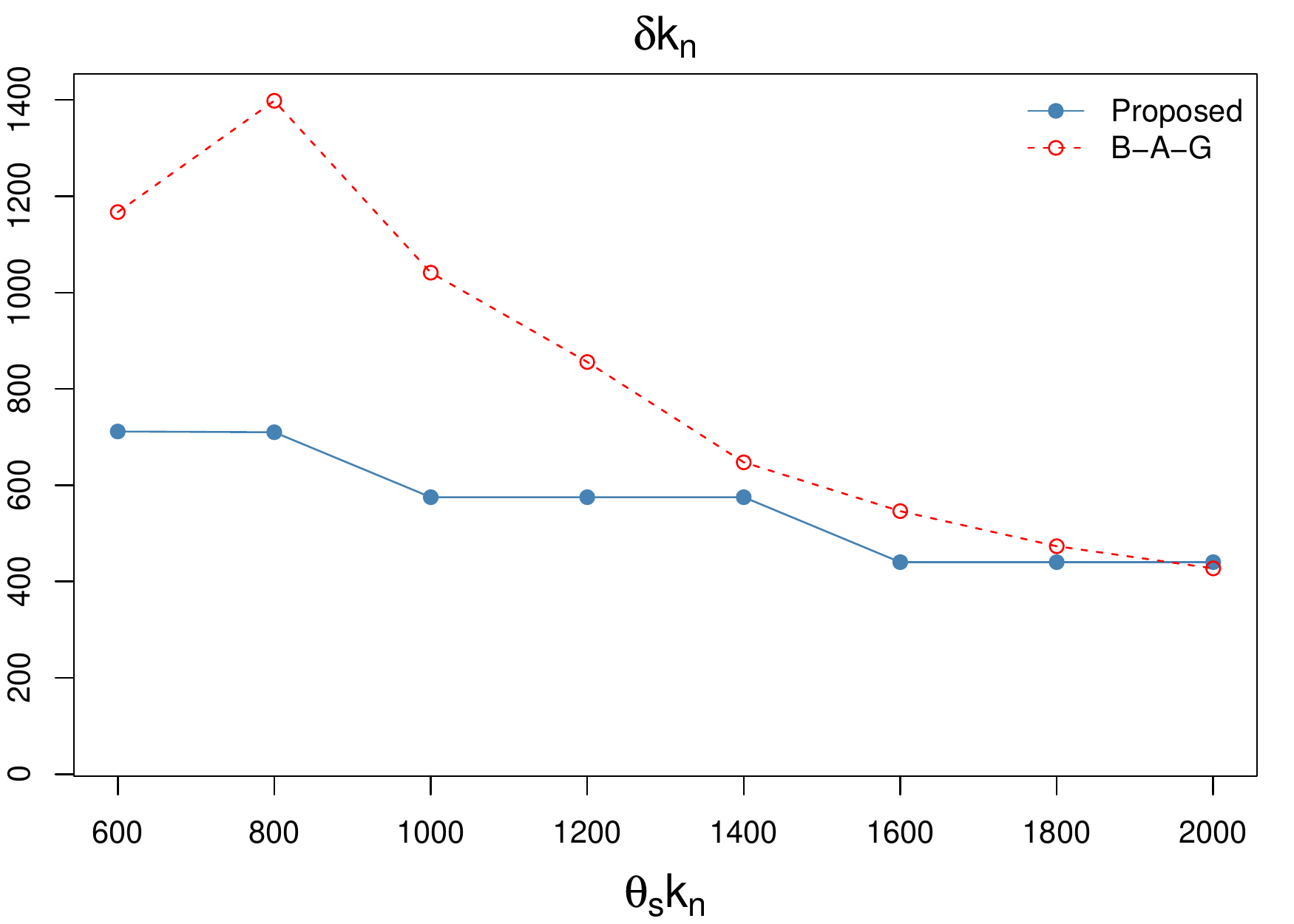}
\end{minipage}
\caption{Estimated number of missing extremes. Left: with the original $ 5000 $ observations. Right: top $ 400 $ values removed} \label{GPlus_truncated}
\end{figure}

\begin{figure}[H]
\centering
\begin{minipage}{.45\textwidth}
  \centering
  \includegraphics[width=0.95\textwidth]{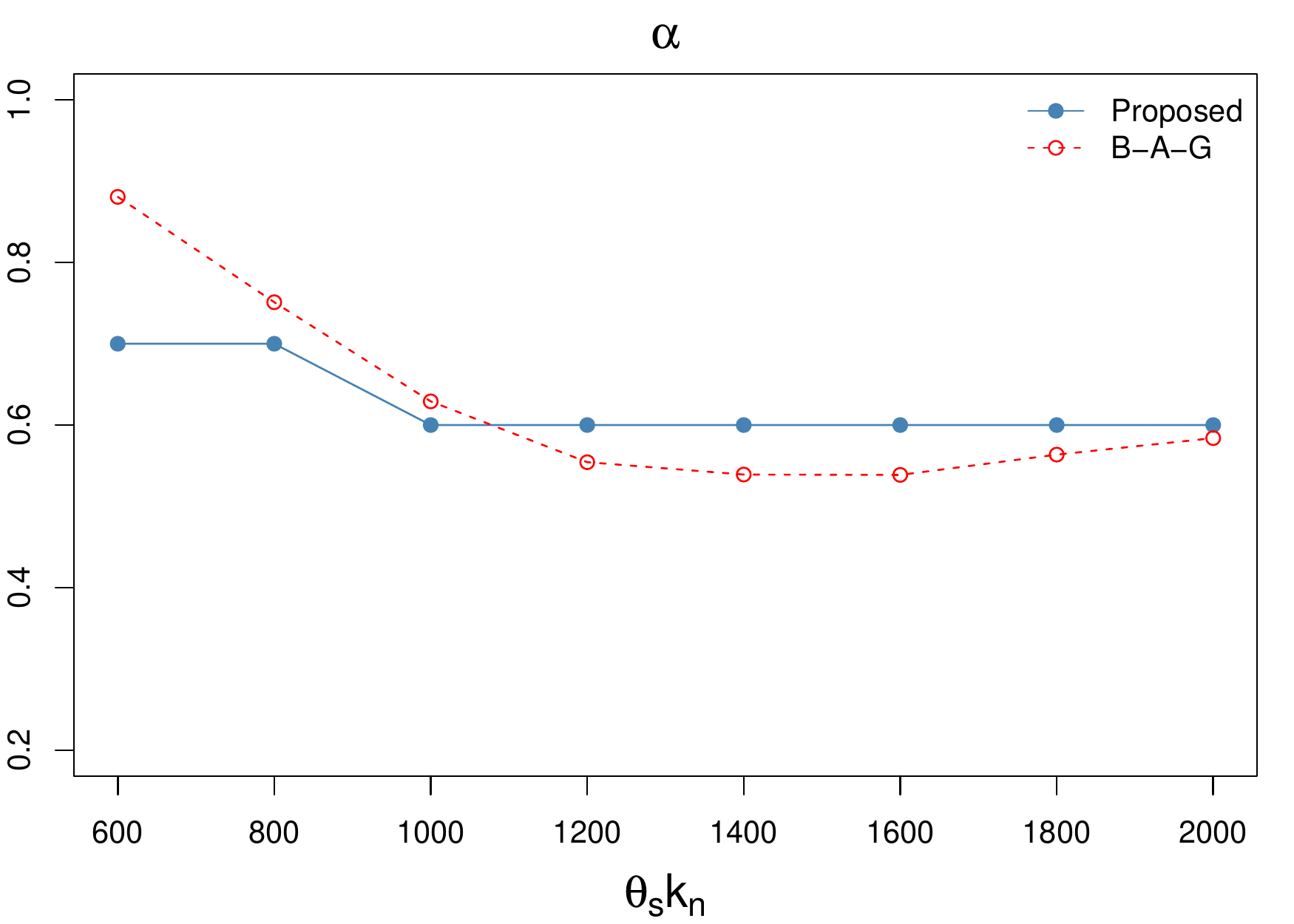}
\end{minipage}%
\hfill
\begin{minipage}{.45\textwidth}
  \centering
  \includegraphics[width=0.95\textwidth]{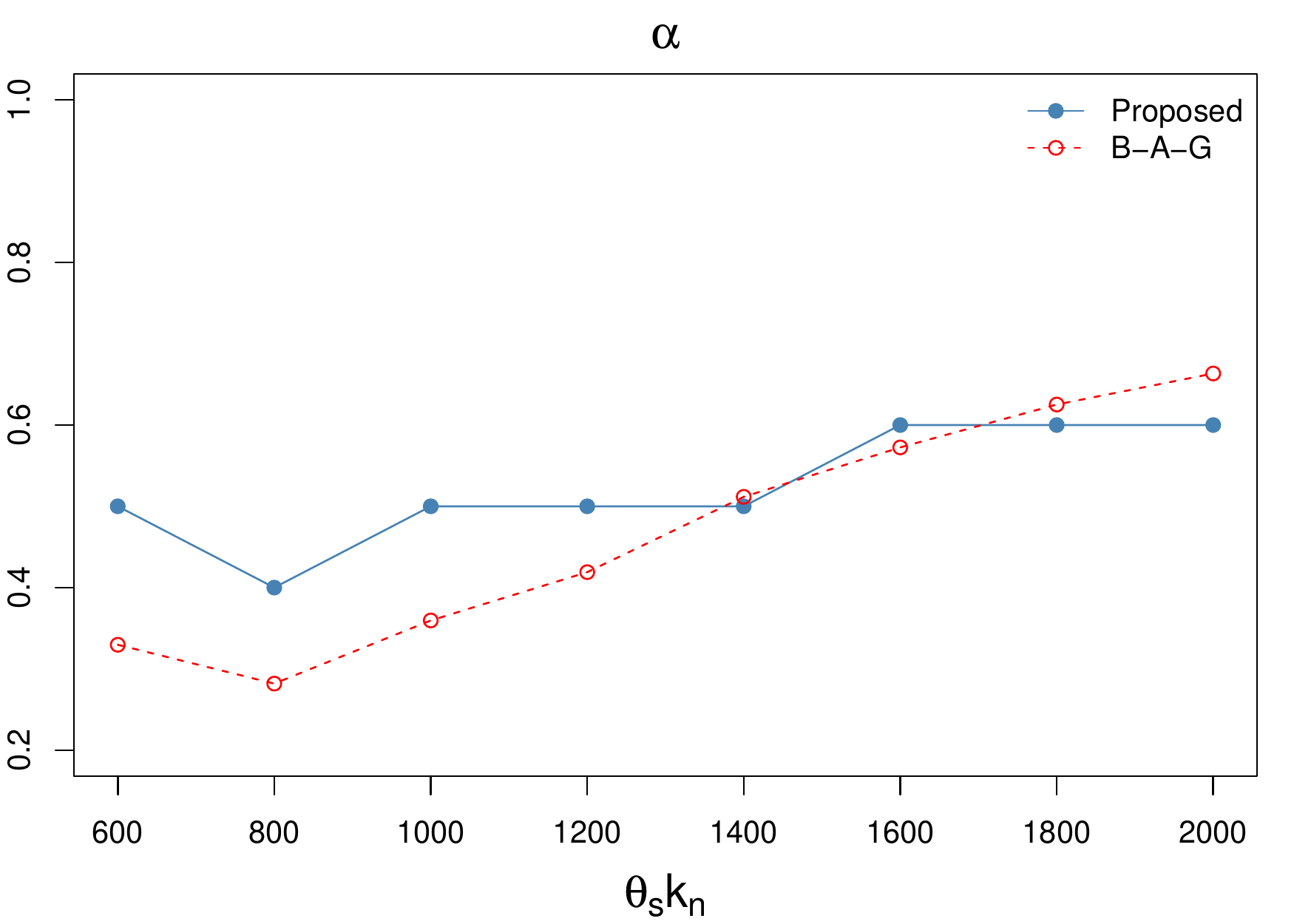}
\end{minipage}
\caption{Estimated tail index. Left: with the original $ 5000 $ observations. Right: top $ 400 $ values removed} \label{GPlus_alpha}
\end{figure}

Figures \ref{GPlus_truncated} and \ref{GPlus_alpha} show, respectively, the estimates of the number of missing extremes and the tail index of the in-degree, before and after the artificial removal. 
It can be seen by comparing the plots on the left and right panels of Figure \ref{GPlus_truncated} that the estimates given the proposed method, which are around $ 150 $ before and $ 550 $ after the artificial removal, reflect reasonably well the additional removal of $ 400 $ top values.
The tail index is mostly estimated to be in the range of $ 0.5 - 0.6 $ and the estimates are reasonably consistent before and after the artificial removal (Figure \ref{GPlus_alpha}).

\subsection{Earthquakes}

While power-law distributions are widely used to model natural disasters such as earthquakes, forest fires and floods, some studies (\citet{Burroughs:2001kc, 2001PApGe.158..741B, Burroughs:2002ip}, \citet{Clark:2013te}, \citet{Beirlant:2016gs, 2016arXiv160602090B}) have observed evidence of truncation in the data available for such events. Causes for the truncation are complex. Possible explanations include physical limitations on the magnitude of the events (\citet{Clark:2013te}), spatial and temporal sampling limitations and changes in the mechanisms of the events (\citet{Burroughs:2001kc, 2001PApGe.158..741B, Burroughs:2002ip}).
In addition, improved detection and rescue techniques might have led to reduction in disaster-related fatalities occurred in recent years.

We apply our method to the dataset of earthquake fatalities (\url{http://earthquake.usgs.gov/earthquakes/world/world_deaths.php}) published by the U.S. Geological Survey, which was also used for demonstration in \citet{Beirlant:2016gs}. 
The dataset is of moderate sample size. It contains information of $ 125 $ earthquakes causing 1,000 or more deaths from 1900 to 2014. 
In the estimation procedure we choose $k_n=10$. Initially the procedure is applied to the original data set. Then we repeat the procedure after artificially removing the $ 10 $ largest of the $ 125 $ values. 
In the estimation, we start from $ \theta_1 = 1 / k_n $ and let $ \theta_i = \theta_{i-1} + 1/k_n $ for $ 1 < i \le s $. 
We consider a sequence of different endpoints $ \theta_s k_n $ and estimate the number of missing extremes and the tail index with different values of $ \theta_s k_n $.
Since the top $ k $ order statistics in the data after removing the top $ 10 $ extreme values are the top $ k  + 10 $ in the original data without the $ 10 $ largest observations, in comparing results before and after the removal, the range of $ \theta_s k_n $ for the data after the removal is shifted to the left by $ 10 $.

\begin{figure}[H]
\centering
\begin{minipage}{.45\textwidth}
  \centering
  \includegraphics[width=0.95\textwidth]{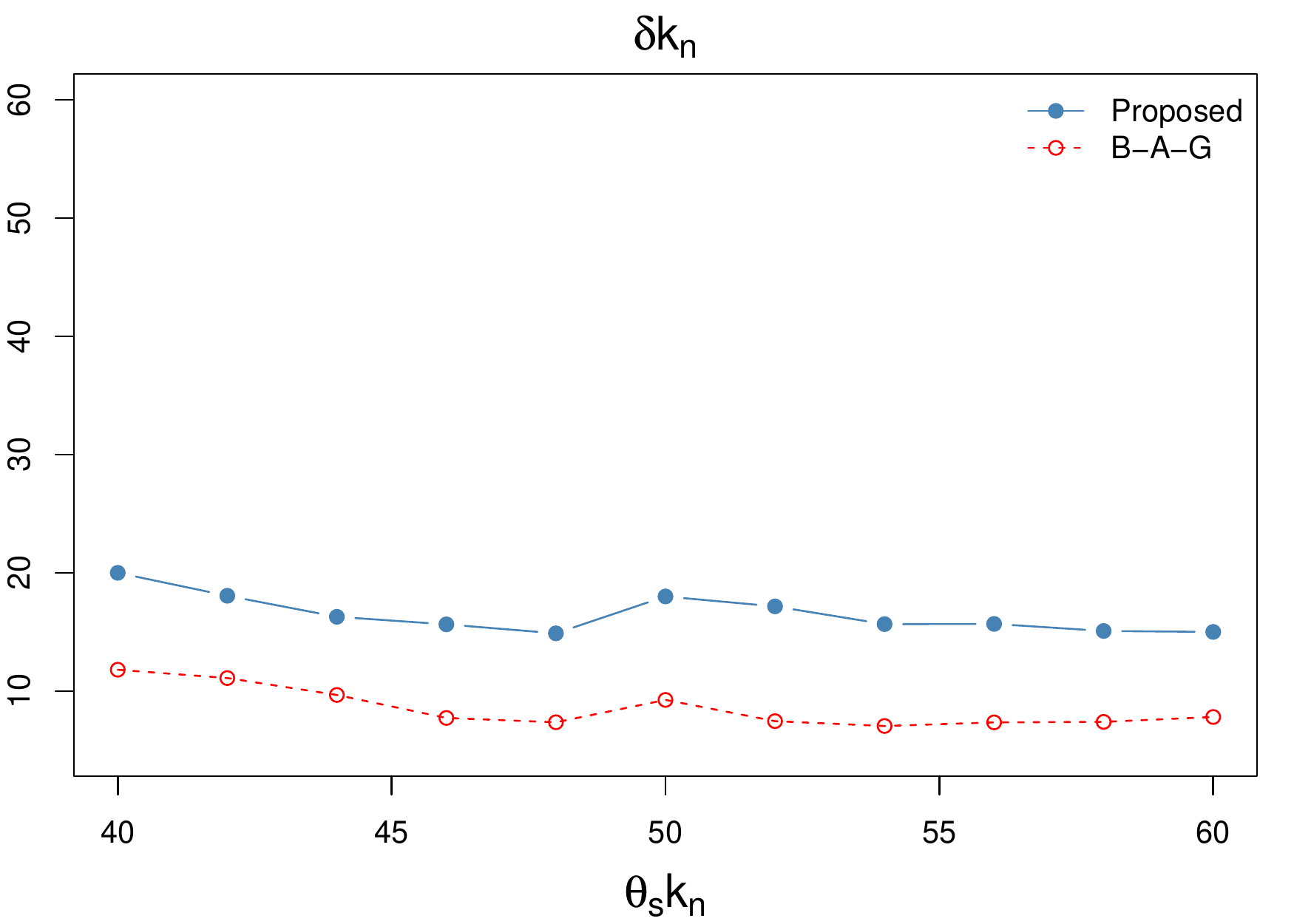}
\end{minipage}
\hfill
\begin{minipage}{.45\textwidth}
  \centering
  \includegraphics[width=0.95\textwidth]{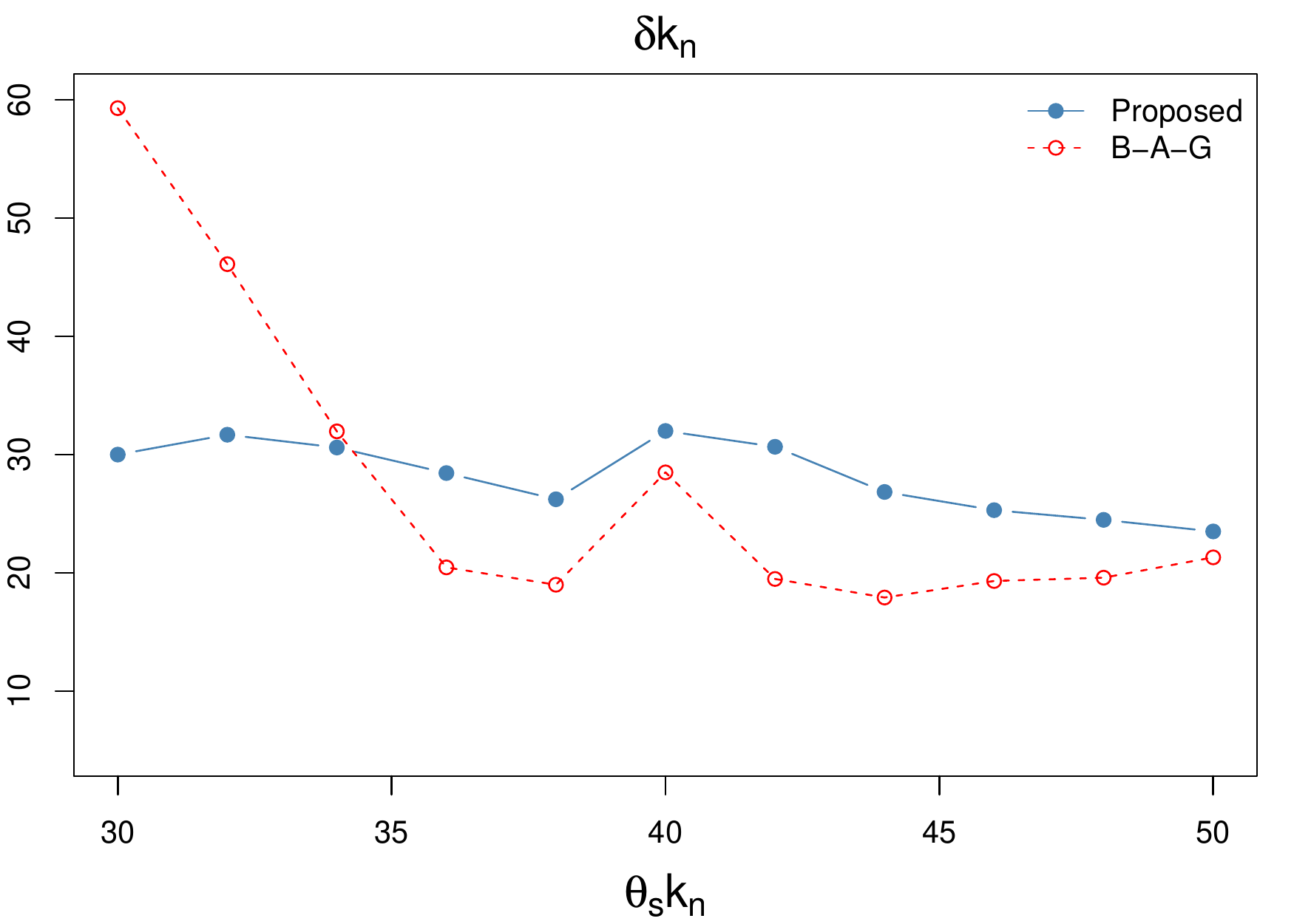}
\end{minipage}
\caption{Estimated number of missing extremes. Left: with the original $ 125 $ observations. Right: with top $ 10 $ values removed} \label{earthquake_truncated}
\end{figure}

\begin{figure}[H]
\centering
\begin{minipage}{.45\textwidth}
  \centering
  \includegraphics[width=0.95\textwidth]{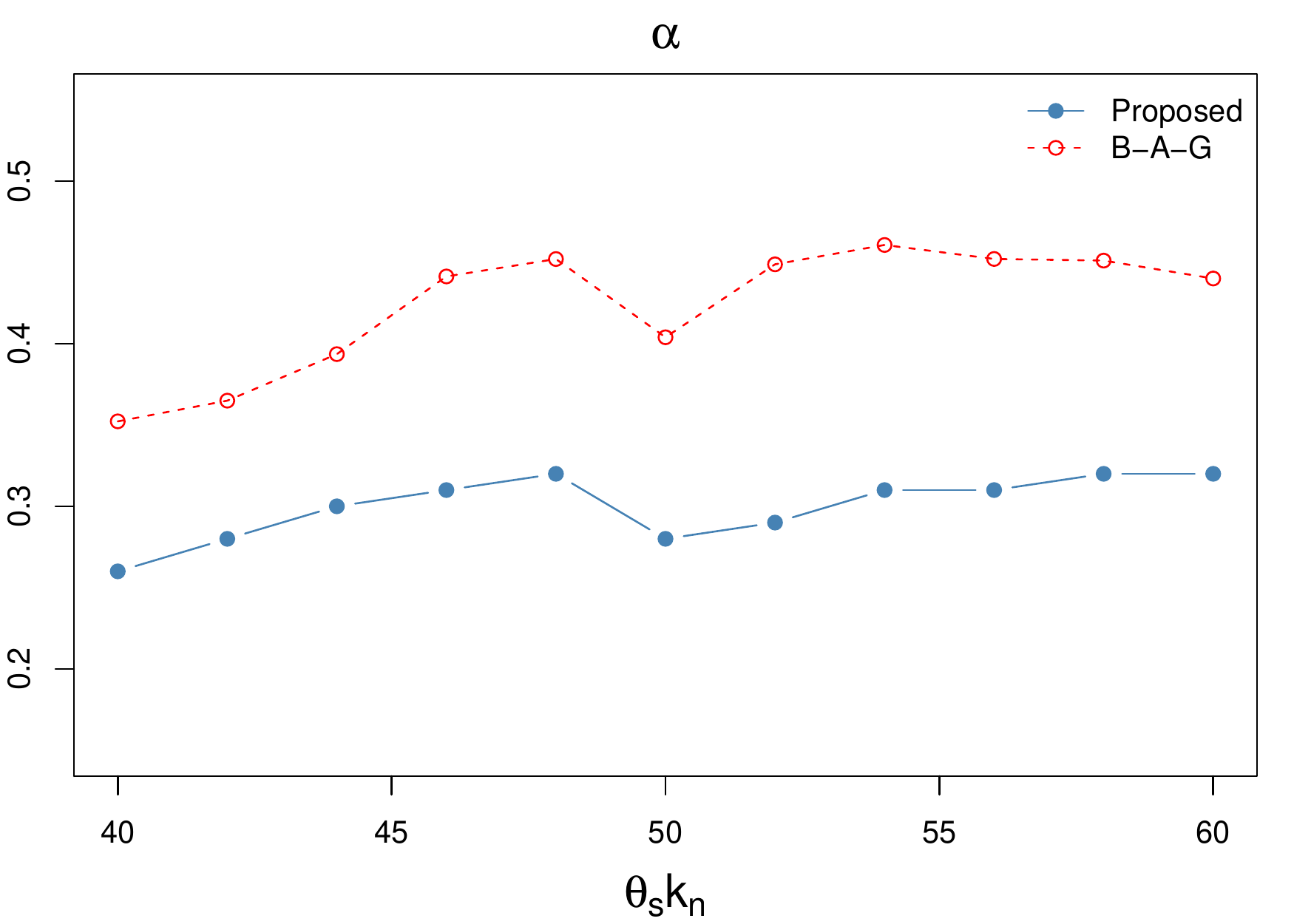}
\end{minipage}%
\hfill
\begin{minipage}{.45\textwidth}
  \centering
  \includegraphics[width=0.95\textwidth]{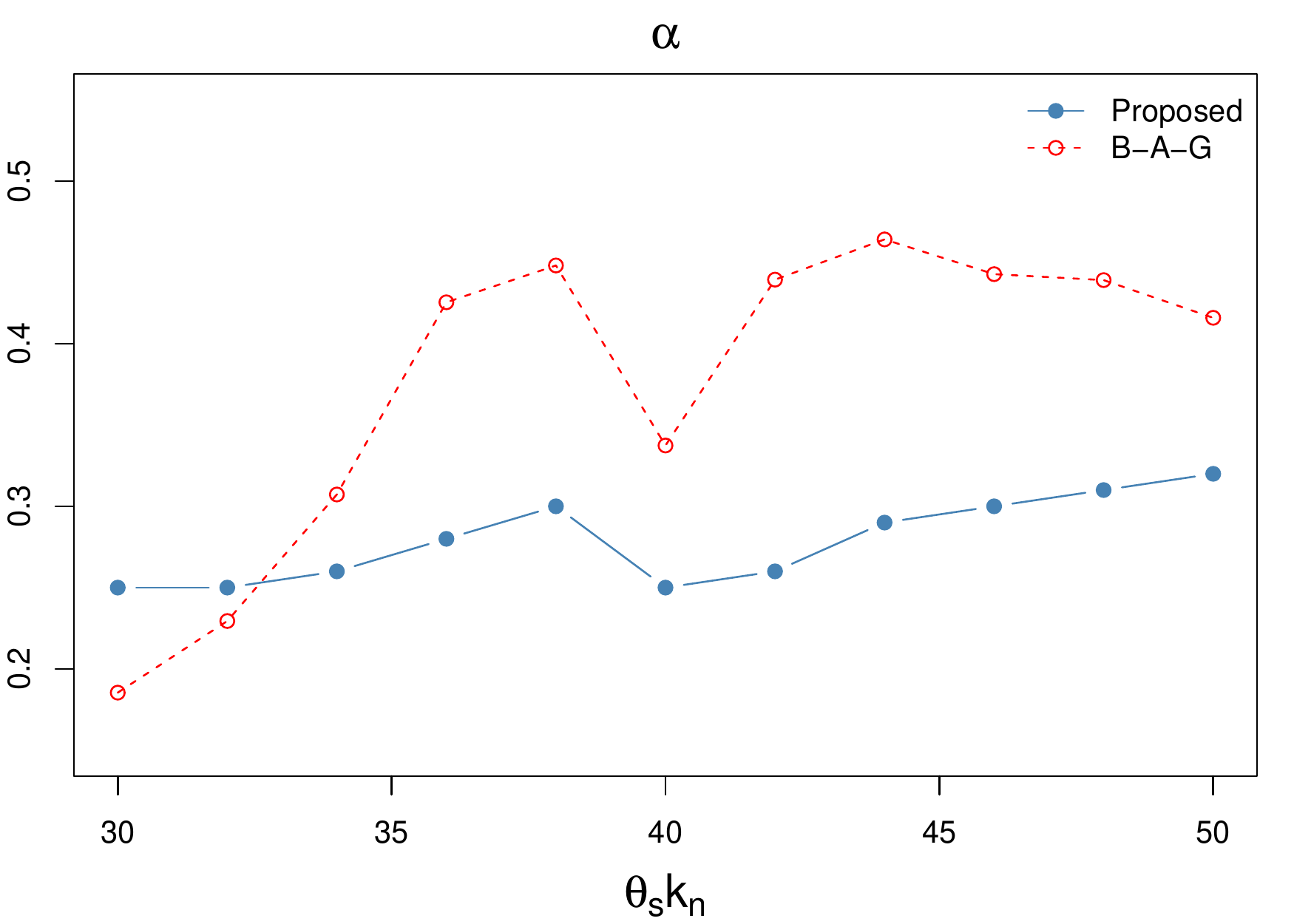}
\end{minipage}
\caption{Estimated tail index. Left: with the original $ 125 $ observations. Right: with top $ 10 $ values removed} \label{earthquake_alpha}
\end{figure}

Figures \ref{earthquake_truncated} and \ref{earthquake_alpha} show the estimates of the number of missing extremes and the tail index of the fatalities. The number of missing extremes is estimated to be around $ 15 - 20 $ for the original data.
After removing the top $ 10 $ earthquakes with the most fatalities, the estimates are now around $ 25 - 30 $, which reflect reasonably well the additional removal (see the left and right panels of Figure \ref{earthquake_truncated}). The estimates of the tail index are reasonably consistent and remain to be in the range of $ 0.25 - 0.3 $ after the additional removal (Figure \ref{earthquake_alpha}).

\newpage

\section*{Appendix}

\subsection{Proof of Theorem 2.1}

Define 
\begin{equation}
S_j = \sum_{i=1}^{j} E_i
\end{equation}
where $ E_i $ are iid standard exponential random variables. Then by Corollary 1.6.9 of \citet{Reiss:1989bw} (see also \citet{Kaufmann:1998fn}),
\begin{equation*}
(X_{(1)},\dots, X_{(n)}) \stackrel{d}{=} \bigg(U\Big(\frac{S_{n+1}}{S_{1}}\Big), \dots, U\Big(\frac{S_{n+1}}{S_{n}}\Big) \bigg).
\end{equation*}
With the second-order condition Eq. \ref{eq: second order cond},  we have from (4.1) - (4.4) of \citet{Drees:2000ic} that for all $ j = 1, \dots,  \lfloor \theta k_n \rfloor $,
\begin{multline*}
\bigg | \log \frac{U(S_{n+1} / S_{\lfloor \delta k_n \rfloor + j})}{U(S_{n+1} / S_{\lfloor \delta k_n \rfloor + \lfloor \theta k_n \rfloor + 1})}
- \frac{1}{\alpha} \log \Big(\frac{S_{\lfloor \delta k_n \rfloor + \lfloor \theta k_n \rfloor + 1}}{S_{\lfloor \delta k_n \rfloor + j}} \Big)\\
- A\Big(\frac{n}{\delta k_n  + \theta k_n  }\Big)  \frac{(1+ \theta / \delta)^\rho - 1}{\rho}
 \bigg|
=o \Big(A\Big(\frac{n}{\delta k_n  + \theta k_n  }\Big)  \Big) \quad \text{a.s.}
\end{multline*}
uniformly in $ (\theta, \delta) \in [m, M]\times [0, M]$ for any $ 0 < m < M$. It follows that
\begin{align*}
H_{k,n}(\theta; \delta) & = \frac{1}{\lfloor \theta k_n \rfloor} \sum_{j = 1}^{\lfloor \theta k_n \rfloor} 
\log U(S_{n+1}/S_{\lfloor \delta k_n \rfloor + j}) - \log U(S_{n+1}/S_{\lfloor \delta k_n \rfloor + \lfloor \theta k_n \rfloor + 1}) \\
& = \frac{1}{\alpha}  \frac{1}{\lfloor \theta k_n \rfloor} \sum_{j = 1}^{\lfloor \theta k_n \rfloor}  
\log \Big(\frac{S_{\lfloor \delta k_n \rfloor + \lfloor \theta k_n \rfloor + 1}}{S_{\lfloor \delta k_n \rfloor + j}} \Big) \\
& \qquad+ \frac{1}{\lfloor \theta k_n \rfloor}  A\Big(\frac{n}{ \delta k_n +  \theta k_n }\Big) 
\rho^{-1} \sum_{j=1}^{\lfloor \theta k_n \rfloor} \Big[\Big(\frac{\theta k_n + \delta k_n +1}{\delta k_n + j}\Big)^\rho- 1 \Big] 
+ o \Big(A\Big(\frac{n}{ k_n  }\Big)  \Big) \quad \text{a.s.}.
\end{align*}
Since the $ E_j^* :=  j \log (S_{j+1}/S_{j} )$ are iid standard exponential random variables (\citet{Reiss:1989bw}), observe that
\begin{align*}
& \frac{1}{\lfloor \theta k_n \rfloor}  \sum_{j=1}^{\lfloor \theta k_n \rfloor} \log \Big(\frac{S_{\lfloor \delta k_n \rfloor + \lfloor \theta k_n \rfloor + 1}}{S_{\lfloor \delta k_n \rfloor + j}} \Big) 
= \frac{1}{\lfloor \theta k_n \rfloor}  \sum_{j=\lfloor \delta k_n \rfloor + 1}^{\lfloor \delta k_n \rfloor + \lfloor \theta k_n \rfloor} \Big(1 - \frac{\lfloor \delta k_n \rfloor}{j}\Big) E_j^* \\
& =  \frac{1}{\lfloor \theta k_n \rfloor} \sum_{j=\lfloor \delta k_n \rfloor + 1}^{\lfloor \delta k_n \rfloor + \lfloor \theta k_n \rfloor}
\Big(1 - \frac{\lfloor \delta k_n \rfloor}{j}\Big) 
(E_j^* - 1) 
+ 1 - \frac{\lfloor \delta k_n \rfloor}{ \lfloor \theta k_n \rfloor}  \sum_{j=\lfloor \delta k_n \rfloor + 1}^{\lfloor \delta k_n \rfloor + \lfloor \theta k_n \rfloor} \frac{1}{j}\\
& = \RN{1} + 1 - \frac{\delta}{\theta} \log\Big(\frac{\theta}{\delta} +1 \Big) + o(1/\sqrt{k_n})
\end{align*}
uniformly in $(\theta, \delta)\in [m, M] \times [0, M]$,
where $ \RN{1} =  \sum_{j=\lfloor \delta k_n \rfloor + 1}^{\lfloor \delta k_n \rfloor + \lfloor \theta k_n \rfloor} 
(1 - \lfloor \delta k_n \rfloor / j) (E_j^* - 1) /  \lfloor \theta k_n \rfloor $.
Using the Koml\'{o}s - Major - Tusn\'{a}dy approximation (\citet{Komlos:1975ch, Komlos:1976kn}), there exists a standard Brownian motion $ \tilde{W} $ such that
\begin{equation*}
\RN{1}  =  \frac{1}{\lfloor \theta k_n \rfloor}   \int_{\lfloor \delta k_n \rfloor + 1}^{\lfloor \delta k_n \rfloor + \lfloor \theta k_n \rfloor} \Big(1 - \frac{\lfloor \delta k_n \rfloor}{\lfloor y \rfloor} \Big) d\tilde{W}(y) + o(1/\sqrt{k_n}), \quad \text{a.s.}
\end{equation*}
Consider the time change $ x = y/k_n  $ and let $ W(x) = \tilde{W}(x k_n)/\sqrt{k_n} $, then
\begin{equation*}
\RN{1} = \frac{1}{\theta \sqrt{k_n}} \int_{\delta}^{\delta + \theta} \Big( 1-\frac{\delta}{x} \Big) dW(x) + o(1/\sqrt{k_n}), \quad \text{a.s.}
\end{equation*}
Summarizing, 
\begin{equation*}
H_{k,n}(\theta; \delta) = \frac{1}{\alpha} \Big(1 - \frac{\delta}{\theta} \log  \Big(\frac{\theta}{\delta} +1\Big)\Big) + \frac{1}{\alpha} \frac{1}{\theta \sqrt{k_n}} \int_{\delta}^{\delta + \theta} \Big(1 - \frac{\delta}{x} \Big) dW(x) + o(1/\sqrt{k_n}) + \RN{2}, \quad \text{a.s.},
\end{equation*}
where
\begin{equation*}
\RN{2} = \frac{1}{\lfloor \theta k_n \rfloor}  A\Big(\frac{n}{ \delta k_n +  \theta k_n }\Big) \frac{1}{\rho}
\sum_{j=1}^{\lfloor \theta k_n \rfloor} \Big[\Big(\frac{\theta k_n + \delta k_n +1}{\delta k_n + j}\Big)^\rho- 1 \Big].
\end{equation*}
The Riemann sum
\begin{align*}
\frac{1}{\lfloor \theta k_n \rfloor}
\sum_{j = 1}^{\lfloor \theta k_n \rfloor} 
\left[ \bigg( \frac{\delta k_n + \theta k_n + 1}{\delta k_n + j } \bigg)^\rho - 1 \right] 
& \rightarrow 
\int_0^1 \bigg(\frac{\delta/\theta + 1}{\delta/\theta + x} \bigg)^\rho dx - 1 \\
& = 
\begin{cases}
\frac{1 + (\theta/\delta)\rho- (\theta/\delta + 1)^\rho}{(\theta/\delta) (1-\rho)}, &\quad \delta > 0, \\
\frac{\rho}{1- \rho} , &\quad \delta = 0.
\end{cases}
\end{align*}
The error between the Riemann sum and the limit can be bounded by
\begin{equation*}
\frac{1}{\lfloor \theta k_n \rfloor} \bigg[ 1 - \bigg(\frac{\delta/\theta + 1}{\delta/\theta} \bigg)^\rho \bigg] \le \frac{1}{\lfloor \theta k_n \rfloor}.
\end{equation*}
Since $A$ is regular varying with index $\rho$,
\begin{equation*}
\RN{2}  = A\Big(\frac{n}{ k_n }\Big) \frac{1}{(\delta + \theta)^\rho} \frac{1 + (\theta/\delta)\rho- (\theta/\delta + 1)^\rho}{(\theta/\delta) (1-\rho)\rho}
+ o\Big(A\Big(\frac{n}{ k_n }\Big)\Big),
\end{equation*}
where $A(n / k_n) \sim O(1/ \sqrt{k_n})$ for $\lambda >0$ in Eq. \ref{cond lambda} and  $A(n / k_n) \sim o(1/ \sqrt{k_n})$ for $\lambda = 0$. Therefore Part \textbf{(a)} follows.

To show Part \textbf{(b)}, we have from Eq. \ref{eq: second order cond} and the fact that $ A $ is regular-varying with index $ \rho $,
\begin{equation}
\sqrt{k_n}  A\Big(\frac{n}{ \delta k_n +  \theta k_n }\Big)  \rightarrow \frac{\lambda} {(\delta+\theta)^{\rho} }
\end{equation}
and thus 
\begin{equation}
\sqrt{k_n} \RN{2} \rightarrow \lambda b_{\delta, \rho}(\theta) 
\end{equation}
and
\begin{equation*}
\sqrt{k_n} \Big(H_{k,n}(\cdot; \delta) - \alpha^{-1} g_\delta(\cdot)\Big) - \lambda b_{\delta, \rho}(\cdot) \stackrel{d}{\rightarrow} \alpha^{-1} G_\delta(\cdot).
\end{equation*} 
The covariance function
\begin{align*}
\text{Cov}\big( G_\delta(\theta_1), G_\delta(\theta_2) \big) 
& = \frac{1}{\theta_1 \theta_2} \text{Var} \Big[\int_{\delta}^{\delta + \theta_1 \wedge \theta_2} \Big(1 - \frac{\delta}{x}\Big) dW(x)\Big] \\
& = \frac{1}{\theta_1 \theta_2} \int_{\delta}^{\delta + \theta_1 \wedge \theta_2} \Big(1 - \frac{\delta}{x}\Big)^2 dx\\
& = \begin{cases}
\frac{1}{\theta_1 \theta_2} \bigg[
\theta_1 \wedge  \theta_2 
- 2 \delta \log \Big(1 + \frac{ \theta_1 \wedge  \theta_2}{\delta} \Big) + \frac{\delta (\theta_1 \wedge \theta_2) } {\delta + (\theta_1 \wedge \theta_2)} \bigg], & \delta > 0\\
\frac{1}{\theta_1 \vee \theta_2}, &\delta = 0.
\end{cases}
\end{align*}
The covariance function of the two-parameter process $ \tilde{G}(\theta, \delta) $ can be shown similarly.

\newpage

\bigskip
\begin{center}
{\large\bf SUPPLEMENTARY MATERIAL}
\end{center}

\begin{description}

\item[R Code for simulations and real data examples:] Code for R algorithms used to produce illustrations in Sections \ref{sec: intro} and \ref{sect: conv} and estimation results in Sections \ref{sect: sim} and \ref{sect: data}. (.r files)

\item[Earthquake fatality data set:] Data set used in the illustration in Section \ref{sect: data}. (comma-separated values (CSV) file)

\item[R Shiny web application:] \url{https://jingjing.shinyapps.io/hewe2}. 
This application can be applied to the user's own data to estimate parameters with real-time computation and to interactively visualize results based on user inputs. Moreover, the user can artificially remove a number of extreme values from the data and compare estimation results before and after the removal.
\end{description}

\section*{Acknowledgements}
 The authors would like to thank Zhi-Li Zhang for providing the Google+ data. This research is funded by ARO MURI grant W911NF-12-1-0385.

\bibliographystyle{apalike}

\end{document}